\documentclass[a4paper, 12pt]{article}
\usepackage{amsmath,amssymb,latexsym, enumerate, color}
\usepackage{graphicx}
\author{ Fr\'ed\'eric H\'erau\\\small Laboratoire de
  Math\'ematiques\\\small Universit\'e de Reims\\\small Moulin de la
  Housse B.P.  1039\\\small 51687 Reims cedex 2, France\\\small and
  FRE 3111 CNRS\\\small herau@univ-reims.fr \and Michael
  Hitrik\footnote{The partial support by the NSF under grant
    DMS-0653275 is gratefully acknowledged.}
  \\\small Department of Mathematics \\\small University of California
  \\\small Los Angeles \\\small CA 90095-1555, USA\\\small
  hitrik@math.ucla.edu \and
  Johannes Sj\"ostrand\footnote{A part of this work was carried out at
    the Centre Interfacultaire Bernoulli (CIB), Lausanne, and we acknowledge
  support from CIB and NSF.}\\\small IMB, Universit\'e de Bourgogne\\
  \small 9, Av. A. Savary, BP 47870\\
  \small FR-21078 Dijon C\'edex\\
  \small  and UMR 5584, CNRS\\
  \small johannes.sjostrand@u-bourgogne.fr} \date{} \title{Tunnel effect and
  symmetries for Kramers Fokker-Planck type operators\footnote{{Ce travail a b\'en\'efici\'e d'une aide de l'Agence Nationale de la Recherche
portant la r\'ef\'erence ANR-08-BLAN-0228-01}}\\
\textit{\large Dedicated to Louis Boutet de Monvel}
}

\def\ccc{{\cal C}}\def\ddd{{\cal D}}
 
\def\iii{{\cal I}} \def\kkk{{\cal K}}\def\lll{{\cal L}}
 \def\ooo{{\cal O}}

\def\uuu{{\cal U}}\def\vvv{{\cal V}}
  
\def\R{\mathbb R}

\def\D{\partial}\def\eps{\varepsilon}
\def\norm#1{\Vert#1 \Vert}

\def\abs#1{\left\vert#1\right\vert}
\def\set#1{\left\{#1\right\}}
\def\seq#1{\left<#1\right>}
\def\sep#1{\left(#1\right)}
\def\Re{{\mathrm Re\,}} \def\Im{{\mathrm Im\,}}
\def\defegal{\stackrel{\textrm{ def}}{=}}

\newcommand{\preuve}[1][\!\!]{\noindent{\it Proof #1. \ \ }}

\def\supp{{\mathrm{supp\,}}}
\def\tA{{  A^{\mathrm {t}}}}
\newcommand\tk{{\kappa ^{\mathrm{t}}}}

\def\remark{\refstepcounter{theorem}\bigskip\noindent\bf Remark \thetheorem\rm\ }
\newcommand{\proof}[1][\!\!]{\noindent{\it Proof #1. \ \ }}

\def\un{ 1 \hspace{-0.9mm} {\rm l}}

\newtheorem{dref}{Definition}[section] \newtheorem{lemma}[dref]{Lemma}
\newtheorem{theo}[dref]{Theorem} \newtheorem{prop}[dref]{Proposition}
\def\remark{\refstepcounter{dref}\bigskip\noindent\bf Remark \thedref.\rm\ }
\newtheorem{cor}[dref]{Corollary}
\newtheorem{hyp}[dref]{Hypothesis}
\renewcommand{\thedref}{\thesection.\arabic{dref}}

\newcommand{\ekv}[2]{\begin{equation}\label{#1}#2\end{equation}}

\newcommand{\no}[1]{(\ref{#1})} 
\begin{document}

\maketitle
\begin{abstract}
  We study operators of Kramers-Fokker-Planck type in the
  semiclassical limit, assuming that the exponent of the associated
  Maxwellian is a Morse function with a finite number $n_0$ of local
  minima. Under suitable additional assumptions, we show that the
  first $n_0$ eigenvalues are real and exponentially small, and
  establish the complete semiclassical asymptotics for these
  eigenvalues.  \medskip \par \centerline{\bf R\'esum\'e}

Nous \'etudions des op\'erateurs de type Kramers-Fokker-Planck dans
 la limite semi-classique quand l'exposant du
 maxwellien associ\'e est une fonction de Morse avec un nombre fini
 $n_0$ de minima locaux. Sous des hypoth\`eses suppl\'ementaires
 convenables, nous montrons que les premi\`eres 
 $n_0$ valeurs propres sont r\'eelles et exponentiellement petites et
 nous \'etablissons leur asymptotique semi-classique compl\`ete.
\end{abstract}
\tableofcontents

\section{Introduction}\label{in}
\setcounter{equation}{0}

In this article we shall continue the study that we started in
\cite{KFP2} of the exponentially small eigenvalues of
the Kramers-Fokker-Planck operator
\ekv{in.0}
{P=y\cdot h\partial _x-V'(x)\cdot h\partial _y+\frac{\gamma }{2}(y-h\partial _y)\cdot (y+h\partial _y)
\hbox{ on }{\R}^{2d}={\R}^d_x\times {\R}^d_y,\ \gamma >0,}
and similar operators.  Suppose that the potential $V\in C^\infty
(\mathbb{R}^d;\mathbb{R})$ is a Morse function on $\R^d$ such that
\begin{equation}
\label{eq011}
\D^{\alpha}V ={\cal O}(1),\quad \abs{\alpha}\geq 2,
\  \ \ \ \ \textrm{and} \ \ \ \
\abs{\nabla V}\geq 1/C,\quad \textrm{for}\quad \abs{x}\geq C>0.
\end{equation}

{ Then (as we shall review) $P$ is maximally accretive and has a unique
closed extension from ${\cal S}(\mathbb{R}^{2n})$ to an unbounded
operator: $L^2\to L^2$ that we shall also denote by $P$. The spectrum
$\sigma (P)$ is contained in the closed right half-plane. Assume for
simplicity that 
\begin{equation}\label{in.1}
V(x)\to +\infty ,\ |x|\to \infty .
\end{equation}
Then the Maxwellian $e^{-(y^2/2+V(x))/h}$ belongs to the kernel of
$P$, so $0\in \sigma (P)$. In \cite{HeSjSt05} the eigenvalues in any
band $0\le \Re z\le Ch$ were determined in the limit $h\to 0$ modulo
${\cal O}(h^\infty )$, for every fixed $C$. They are of the form $\mu
h+o(h)$ where $\mu \in \mathbb{C}$ are the eigenvalues of the
quadratic approximations of $P$ with $h=1$ at the various critical
points of $V$. These $\mu $ values are explicitly known and belong to
a cone $|\Im \mu |\le {\cal O}(\Re \mu )$. The $o(h)$ terms have
complete asymptotic expansions in powers of $h$ (and we need
fractional powers when certain multiplicities are present). 

Assume that $V$ has $n_0$ local minima, ${\bf m}_1,{\bf m}_2,...,{\bf m}_{n_0}$ and
$n_1$ critical points of index 1 that we shall call saddle
points. Then precisely $n_0$ of the above eigenvalues are $o(h)$
(i.e. with $\mu =0$) and they are actually ${\cal O}(h^\infty )$ (as
can be understood intuitively by using truncations of the Maxwellian
near the local minima as exponentially accurate quasimodes).

As we shall review below, it follows from the analysis in \cite{KFP2}
that these $n_0$ eigenvalues are actually exponentially small. In that
paper we were able to establish the exponential decay rate and a full
asymptotic expansion of the prefactor for the smallest non-vanishing
eigenvalue when $n_0=2$. 

In this paper we treat the case of general $n_0$ and our results are
similar to those for the Witten Laplacian by Bovier, Eckhoff, Gayrard,
Klein \cite{BEGK}, \cite{BGK1}, obtained with probabilistic methods and
the ones by Helffer, Klein, Nier \cite{HKN}, with simplifications by
Nier \cite{Nier04} and le Peutrec \cite{Lep08} with full asymptotic
expansions, based on the WKB-analysis and Agmon estimates in
\cite{HeSj85}.

In order to state one of our results completely we describe first a
geometric framework slightly generalizing the procedure in \cite{HKN},
\cite{Nier04}: Let ${\bf m}_1$ be a global minimum of $V$ and put
$E_{{\bf m}_1}=\mathbb{R}^{d}$. For $\sigma \gg 1$, $V^{-1}(]-\infty
,\sigma [)$ is a connected relatively compact open subset of
$\mathbb{R}^d$. When we decrease $\sigma $, $V^{-1}(]-\infty ,\sigma
[)$ remains connected until we reach a critical value $\sigma _2$ where
one of the following happens:

\medskip \par\noindent a) $\sigma _2=\phi ({\bf m}_1)$ and $V(]-\infty
,\sigma _2[)$ is empty. The procedure then stops.

\smallskip \par\noindent b) $V^{-1}(]-\infty ,\sigma _2[)$ is the
finite union of several disjoint components, $E_1\cup E_2\cup ...\cup
E_N$, {$N\ge 2$,} where the labelling is chosen so that ${\bf m}_1\in E_1$. For
$k\ge 2$, let ${\bf m}_k\in E_k$ be a minimum of ${{V}_\vert}_{E_k}$
and write $E_{{\bf m}_2}=E_2$, ..., $E_{{\bf m}_N}=E_N$. Notice that none of the
closures of the $E_{j}$ can be disjoint from the union of the other
closures, and that the intersections of their boundaries are
finite unions of saddle points. $\sigma _2$ is the common value of $V$
at those saddle points. For $2\le k\le N$, we put $\sigma
(E_{{\bf m}_k})=\sigma _2$.

\medskip In case b) we pick (successively) each of $E_j$ and consider
$E_j\cap V^{-1}(]-\infty ,\sigma [)$ with $\sigma $ decreasing from
$\sigma _2$ until one of the scenarios a) or b) appears. In case a) we
stop (with that component $E_j$) and in case b) (say for $\sigma
=\sigma _3<\sigma _2$) we get a finite union of connected
components. Choose a global minimum for each of the new components
except for the one which contains the already selected minimum
${\bf m}_j$. We continue in this way until all the local minima have been
recovered. Then for each local minimum we have an associated connected
component $E_{\bf m}$ of $V^{-1}(]-\infty ,\sigma ({\bf m})[)$. We put 
$$
S({\bf m})=\sigma ({\bf m})-\phi ({\bf m}).
$$ 
See Section \ref{lbl} for a more detailed description of the
procedure. 

\par One of our main results is 

\begin{theo}\label{in1}
Under the above assumptions, the $n_0$ exponentially small eigenvalues
$\mu _1,...,\mu _{n_0}$ are all real and can be labelled in
such a way that 
$$
\mu _k\asymp he^{-S({\bf m}_k)/h},
$$
uniformly when $h\to 0$, and with the convention that $S({\bf m}_1)=+\infty
$, $\mu _1=0$.
\end{theo}

This result is valid also for the Witten Laplacian $-\Delta
_V=d_V^*d_V$, where $d_V:=e^{-V/h}hd e^{V/h}$ and extends those of
\cite{BGK1}, \cite{HKN}, in the sense that no generic assumption is
done on the separation of critical values. 

{Under a generic assumption in the spirit of  \cite{BGK1}, \cite{HKN}, we have
full asymptotic expansions for each of the eigenvalues $\mu _k$.}

{
\begin{theo}\label{in2}
Suppose that $V$ satisfies the same hypotheses as in the preceding theorem,
and  assume in addition the following :

For every critical component $E_k$ created in the above procedure, 
\begin{enumerate}[i)]
\item there is a unique point
 ${\bf m}_k$ in $E_k$ where $\inf_{E_k}V$ is attained, so $E_k=E_{{\bf
     m}_k}$
\item  except in the case when $E_k=\mathbb{R}^d$, there is a unique
  (saddle) point $s_j=s_{j(k)}$ in $\partial E_k$ which is also on the
  boundary of another component of $V^{-1}(]-\infty ,\sigma
  (E_k)[)$. In particular we always have $N=2$ in case b) of the above procedure. 

\end{enumerate}
 Then we have the following asymptotic expansion
\begin{equation} \label{maineq}
\mu_k = h l_k(h)   e^{- 2S_k/h} \ \text{ with } \ \ l_k(h)
  = l_{k,0} + l_{k,1} h + ...,l_{k,0}>0 ,
\end{equation}
where
$$
l_{k,0}=\frac{\widehat{\lambda }_1({\bf s}_j)}{\pi
}\left(\frac{\det V''({\bf m}_k)}{-\det V''({\bf s}_j)}\right)^{\frac{1}{2}}.
$$
Here  $-\widehat{\lambda }_1({\bf
  s}_j)$ is the unique negative eigenvalue of the block matrix
$\left(\begin{matrix}0 &1\\V''({\bf s}_j) &\gamma \end{matrix}\right)$.
\end{theo}}

See Theorem \ref{theogene} and Section \ref{explt}. Under a weaker
generic assumption (see Theorem \ref{fa1}) we get this for the lowest
non-vanishing eigenvalue $\mu _2$. This last result is also valid for
the ordinary Witten Laplacian and seems to be new in that degree of
generality.}

The proof is very much based on the analysis in the paper \cite{KFP2},
where we were able to treat the case of one or two minima. As there we
shall also use in an essential way the supersymmetric structure of the
Kramers-Fokker-Planck operator, due to Bismut \cite{Bi05} and
Tailleur, Tanase-Nicola, Kurchan \cite{TaTaKu06}, see also
\cite{Leb04}.

What made it possible to go beyond the case of two wells was the observation that we have an
additional (generalized PT) symmetry, namely if we conjugate the
Kramers-Fokker-Planck operator $P$ by the unitary and self-adjoint operator $u(x,y)\mapsto u(x,-y)$, then we get the adjoint
$P^*$. This extra symmetry makes it possible to introduce a Hermitian form on $L^2$ for which $P$ is formally self-adjoint and
the restriction of this form to the spectral subspace corresponding to the $n_0$ lowest eigenvalues is positive definite, hence an inner
product. Consequently the restriction of $P$ is self-adjoint and the
$n_0$ lowest eigenvalues are real.

More generally, this extra symmetry entails that the supersymmetric approach followed in \cite{KFP2} reduces to a
self-adjoint problem, very close to that of the Witten Laplacian, when we restrict the attention to the exponentially small
eigenvalues. In particular, we can follow the work by Helffer-Klein-Nier \cite{HKN}, who
adapted probabilistic ideas in the case of the Witten Laplacian.
Le Peutrec \cite{Lep08} simplified some of the linear algebra in \cite{HKN}, and we simplify that even further.

The structure of the paper is as follows. In Section \ref{re}, we
review some of the basic analysis of~\cite{KFP2}, including a
description of the supersymmetric formalism. The generalized PT
symmetry is then explored in Section \ref{pt}, first in the original
Kramers-Fokker-Planck case, and then in the general supersymmetric
setting. In particular, we show here that the first $n_0$ eigenvalues
are real. In Section \ref{lbl} we prepare for the multiple well
analysis, by introducing a suitable labelling of the minima and the
associated saddle point values, as well as a system of quasimodes
adapted to the minima. {Sections \ref{genc} and \ref{explt} are devoted to the analysis
of the exponentially small eigenvalues in the generic case. The
general case is studied in Section \ref{gl}. The main results of this
work are Theorem \ref{theogene}, Proposition \ref{endestim}, Theorem
\ref{gl4} and Theorem \ref{fa1}.} 

It is natural here to acknowledge a
  contribution of
  L.~Boutet de Monvel. One of us had the priviledge to listen to his
  lectures at the Mittag Leffler institute in 1974 about complexes of
  pseudodifferential operators and related operators with double
  characteristics, cf \cite{Bo74}. This type of problems
  continues to be important, and the present paper deals
  once more with that situation.

\section{Review of some results from \cite{KFP2}}\label{re}
\setcounter{equation}{0}

The purpose of this section is to review some basic results from
\cite{KFP2}. We refer to that paper for proofs and more details and we
give proofs only for slightly new variants.
As in \cite{KFP2} some of the analysis works also on compact manifolds
but currently we can go all the way only on $\mathbb{R}^n$ and we
restrict the attention to that case from the beginning.
\subsection{The general case}\label{reg}

On $\mathbb{R}^n$ we consider a second order differential operator
\ekv{re.1} { \begin{split} P=&\sum_{j,k=1}^n hD_{x_j}\circ
  b_{j,k}(x)\circ hD_{x_k}+{1\over 2}\sum_{j=1}^n (c_j(x)h\partial
  _{x_j}+h\partial _{x_j}\circ c_j(x))+p_0(x)\\ =&P_2+iP_1+P_0,\quad
  D_{x_j}={1\over i}{\partial \over \partial x_j},\end{split}} where the coefficients
$b_{j,k}$, $c_j$, $p_0$ are assumed to be smooth and real, with
$b_{j,k}=b_{k,j}$. To $P$ we associate the symbol in the semi-classical
sense,
\ekv{re.2} {p(x,\xi )=p_2(x,\xi )+ip_1(x,\xi )+p_0(x),}
\ekv{re.3} {p_2(x,\xi )=\sum_{j,k=1}^n b_{j,k}(x)\xi _j\xi _k,\
  p_1(x,\xi )=\sum_{j=1}^n c_j(x)\xi _j,} so that $p_j(x,\xi )$ is a
real-valued polynomial in $\xi $, homogeneous of degree
$j$. (It is well-defined on $T^*\mathbb{R}^n$ and coincides with the Weyl symbol
mod ${\cal O}(h^2)$ locally uniformly.) We assume that
\ekv{re.4}
{p_2(x,\xi )\ge 0,\ p_0(x)\ge 0.}

\par We impose the following growth conditions
at infinity:
\ekv{re.5}
{
\partial _x^\alpha b_{j,k}(x)={\cal O}(1),\ \vert \alpha \vert \ge 0,
}
\ekv{re.6}
{
\partial _x^\alpha c_j(x)={\cal O}(1),\ \vert \alpha \vert \ge 1,
}
\ekv{re.7}
{
\partial _x^\alpha p_0(x)={\cal O}(1),\ \vert \alpha \vert \ge 2.
}

Put
\ekv{re.7.1}
{\nu (x,\partial _x)=\sum_1^n c_j(x)\partial _{x_j}.}
\par Let $f(t)\in C^\infty ([0,\infty [;[0,3/2])$ be an increasing function
with $f(t)=t$ on $[0,1]$, $f(t)=3/2$ on $[2,\infty [$, $f(t)\le t$. Put
$f_\epsilon (t)=\epsilon f(t/\epsilon )$, and introduce the
time $T_0$ average of $f_\epsilon \circ p_0$ along the integral
curves of $\nu $,
\ekv{re.8}{\langle f_\epsilon \circ p_0\rangle _{T_0}={1\over
T_0}\int_{-T_0/2}^{T_0/2}f_\epsilon \circ p_0\circ \exp (t\nu )dt.
}

By an averaging (realized by means of weak exponential weights) we
showed in \cite{KFP2} the following result:
\begin{prop}\label{re1}
Let $P$ be of the form {\rm \no{re.1}}, where $b_{j,k},c_j,p_0$ are smooth and real
and satisfy {\rm \no{re.2}}--{\rm \no{re.7}}. Define $\langle f_\epsilon \circ p_0
\rangle_{T_0}$ as in {\rm \no{re.8}}. Let $C_0>$ be sufficiently large.
Then for every $C>0$, put $\epsilon=Mh$ with $M>0$ sufficiently large. Then there exists $\widetilde{C}>0$ such that
\ekv{re.9}{\Vert (p_0+h)^{1\over 2}u\Vert
  \le \widetilde{C}(\Vert (p_0+h)^{-{1\over 2}}(P-z)u\Vert +h^{1\over
    2}\Vert u\Vert _{\{ \langle f_{\epsilon} \circ p_0\rangle
    _{T_0}\le \frac{2\epsilon }{C_0} \}}),}
for $u\in {\cal S}$, $\Re z\le Ch$.
\end{prop}

Combining this with a very simple and direct a  priori estimate, we also
deduced that
\ekv{re.10}{
\Vert B^{1\over 2}hDu\Vert ^2 \le \Vert (p_0+h)^{-{1\over 2}}(P-z)u\Vert
\Vert (p_0+h)^{1\over 2}u\Vert +C\Vert h^{1\over 2}u\Vert ^2.
}

Using elementary coercivity estimates and pseudodifferential
machinery, we got in~\cite{KFP2},
\begin{prop}\label{re2}
$${\cal R}(P-z)=L^2,\ \Re z<0.$$
Here $P:L^2\to L^2$ denotes the graph closure of $P: {\cal S}(\mathbb{R}^n) \rightarrow {\cal S}(\mathbb{R}^n)$.
\end{prop}

\begin{cor}\label{re3}
The maximal closed extension $P_{\rm max}$ of $P$ (with domain given by
$\{ u\in L^2;\, Pu\in L^2\}$ coincides with the graph closure
(the minimal closed extension), already introduced.
\end{cor}
Thus $P$ is maximally accretive. See \cite{HeNi04}, \cite{HeNi05} for
earlier and closely related results in this direction.

\medskip
\par
We shall now discuss some weighted estimates for $P-z$, leading to simplifications and improvements in Section 3 of~\cite{KFP2}. These improvements
will be used later on in this section.

\par
Let $\lambda  =\lambda (\rho )\in C^\infty (\mathbb{R}^{2n};]0,+\infty[)$ satisfy the bounds,
\begin{equation}
\label{weest.2}
\lambda ,\, \lambda ^{-1}={\cal O}(\langle \rho \rangle ^{N_0}),\
{\langle \rho \rangle =(1+|\rho |^2)^{1/2},}
\end{equation}
for some fixed $N_0\ge 0$ and assume that
\begin{equation}\label{weest.3}
\lambda ={\cal O}_1(\lambda ),
\end{equation}
in the sense that
\begin{equation}\label{weest.3'}
\partial _{\rho }^\alpha \lambda ={\cal O}(\lambda \langle \rho
\rangle^{-|\alpha |}),\ \alpha \in \mathbb{N}^{2n}.
\end{equation}

\par
In the proof of Proposition \ref{re2}, given in \cite{KFP2}, we checked that
$$
\lambda =\lambda _{\epsilon ,N}=\langle \epsilon \rho \rangle^N
$$
satisfies these assumptions uniformly for $0<\epsilon \le 1$, when
$N\in \mathbb{R}$ is fixed. The discussion there for such particular
$\lambda =\lambda _{\epsilon ,N}$ generalizes to $\lambda $ in
(\ref{weest.2}), (\ref{weest.3}) and we have with $\Lambda
:=\mathrm{Op}_ h(\lambda )$ (when $h$ is sufficiently small):
\begin{itemize}
\item The symbol of $\Lambda ^{-1}$ is equal to $\lambda ^{-1}+{\cal
    O}_1(h^2\lambda ^{-1}\langle \rho \rangle^{-2})$.
\item The symbol of $[P,\Lambda ]$ is of the form
  $\frac{h}{i}\{p,\lambda  \}+{\cal O}_0(h^2\lambda )$, where we write
  $a={\cal O}_0(m)$ if $\partial _{\rho }^\alpha a={\cal O}(m)$ for
  every $\alpha \in \mathbb{N}^{2n}$.
\item Here,
$$
\{ p,\lambda \}=\{ p_2,\lambda \}+{\cal O}_0(\lambda )=-\partial
_xp_2\cdot \partial _\xi \lambda +{\cal O}_0(\lambda )={\cal
  O}_0(\lambda \langle \rho \rangle ),
$$
where the term ${\cal O}_0(\lambda \langle \rho \rangle)$ is real-valued.
\end{itemize}
Combining these facts with $h$-pseudodifferential calculus, we see that
\begin{itemize}
\item The symbol of $[P,\Lambda ]\Lambda ^{-1}$ is equal to
  $\frac{h}{i}\frac{\{ p,\lambda \}}{\lambda }+{\cal O}_0(h^2)$.
\end{itemize}

Before continuing the main discussion, we shall give a simplification of the
main step of the proof of Proposition \ref{re2}, which is to
establish:
\begin{equation}\label{weest.4}\begin{split}
&\exists z\in \mathbb{C},\ \Re z<0,\hbox{ such that if }u\in L^2\\
&\hbox{and }(P-z)u=0,\hbox{ then }u=0.\end{split}
\end{equation}
To see this, we choose $\lambda =\lambda _{\epsilon ,N}$ with
$0<\epsilon \ll 1$ and with $N\le -2$ fixed. Then,
$$
0=(\Lambda (P-z) u|\Lambda u)=(\Lambda (P-z)\Lambda
^{-1}\Lambda u|\Lambda u)=((P-z-[P,\Lambda ]\Lambda ^{-1})\Lambda
u|\Lambda u).
$$
Here we take the real part and use that
$\mathrm{Op}_h(\frac{h}{i}\frac{\{p,\lambda \}}{\lambda })$ is
formally skew-adjoint, to get
$$
0=((P_2+P_0-\Re z+\mathrm{Op}_h({\cal O}_0(h^2)))\Lambda u|\Lambda
u)\ge (-\Re z+{\cal O}(h^2))\Vert \Lambda u\Vert^2.
$$
Assuming that$-\Re z\gg h^2$, we conclude that $\Lambda u=0$ and hence
that $u=0$.\hfill{$\Box$}

\medskip
The improvement in comparison with \cite{KFP2} is that we do not only
consider $\Re ((P-z)u|u)$ but modify this expression with weight
factors, while the weight factors in \cite{KFP2} enter at a later
stage.

\par If we let $u\in L^2$ and drop the assumption that $(P-z)u=0$, the
same calculations and Cauchy-Schwarz give
\begin{equation}\label{weest.5}
(-\Re z+{\cal O}(h^2))\Vert \Lambda u\Vert\le \Vert \Lambda (P-z)u\Vert,
\end{equation}
for $\Lambda =\Lambda _{\epsilon ,N}$ with $N\le -2$, and this remains
true for more general $\Lambda =\mathrm{Op}_h(\lambda )$ as in
(\ref{weest.2}), (\ref{weest.3}), provided that $\lambda ={\cal
  O}_1(\langle \rho \rangle^{-2})$.
\begin{prop}\label{weest1}
Let $\Lambda =\mathrm{Op}(\lambda )$ with $\lambda $ as in
{\rm (\ref{weest.2})}, {\rm (\ref{weest.3})}. Then there exists a constant $C>0$
such that if $\Re z\le -Ch^2$, $(P-z)u=v$, $u,v\in L^2$, $\Lambda v\in
L^2$, then $\Lambda u\in L^2$ and {\rm (\ref{weest.5})} holds.
\end{prop}
\par\noindent {\it Proof.}
Choose $-N\ge N_0+2$. Then (\ref{weest.5}) holds with $\Lambda $
replaced by $\Lambda _{\epsilon ,N}\Lambda $:
\begin{equation}
(-\Re z+{\cal O}(h^2))\Vert \Lambda _{\epsilon ,N}\Lambda u\Vert \le
\Vert \Lambda _{\epsilon ,N}\Lambda (P-z)u\Vert,
\end{equation}
uniformly with respect to $\epsilon $. Letting $\epsilon $ tend to
$0$, we get $\Vert \Lambda u\Vert <+\infty $ together with
(\ref{weest.5}).\hfill{$\Box$}

\medskip The proposition can be generalized by letting $u,v\in {\cal
  S}'$, with $\Lambda _0u,\, \Lambda _0v\in L^2$, for some fixed
$\Lambda _0=\mathrm{Op}_h(\lambda _0)$ with $\lambda _0$ as in
(\ref{weest.2}), (\ref{weest.3}).

\par
We shall next recall the dynamical assumptions introduced in~\cite{KFP2}. Consider the non-negative symbol
\ekv{geo.1}
{
\widetilde{p}(x,\xi )=p_0(x)+{p_2(x,\xi )\over \langle \xi \rangle ^2}.
}

\begin{hyp}\label{ny0}\rm
Assume
\ekv{ny.1}
{\hbox{The set }\{x \in \mathbb{R}^n;\, p_0(x)=0,\ \nu (x,\partial _x)=0\} \hbox{ is
finite }=\{ x^1,...,x^N\} .}
\end{hyp}
Let $\rho _j=(x^j,0)$ and introduce the critical set
\ekv{ny.2}
{
{\cal C}=\{ \rho _1,...,\rho _N\}
.}
Notice that $p_1,p_0,p_2,\widetilde{p}$ vanish to second order at each
$\rho _j$.
For functions $q$ on the cotangent space, we generalize the earlier
definition of time $T_0$ average and put
$$
\langle q\rangle _{T_0}={1\over T_0}\int_{-T_0/2}^{T_0/2}q\circ \exp
(tH_{p_1}) dt.
$$

\par We introduce the following dynamical conditions where $T_0>0$ is fixed:
\begin{hyp}\label{ny2}\rm
\ekv{ny.17}
{
\hbox{Near each }\rho _j\hbox{ we have }\langle \widetilde{p}\rangle
_{T_0}\ge {1\over C}\vert \rho -\rho _j\vert ^2,
}
\ekv{ny.18}
{
\hbox{In any set }\vert x \vert \le C,\ {\rm dist\,}(\rho ,{\cal C})\ge
{1\over C},\hbox{ we have } \langle \widetilde{p}\rangle _{T_0}(\rho )\ge {1\over
\widetilde{C}(C)},\ \widetilde{C}(C)>0.
}

\ekv{ny.20}
{\begin{split}
&\forall\hbox{ neighborhood }U\hbox{ of }\pi _x{\cal C},\, \exists C>0\, \hbox{ such that }\forall \
x\in \mathbb{R}^n \setminus U,
\\
&{\rm meas\,}(\{ t\in [-{T_0\over 2},{T_0\over 2}];\, p_0(\exp t\nu
(x))\ge {1\over C}\})\ge {1\over C}. \end{split}
}
\end{hyp}

We know that $P$ has no spectrum in the
open left half-plane. By using an elaborate method of microlocal
exponential (but bounded!) weights we showed in \cite{KFP2} the
following result.
\begin{prop}\label{app2}
For every constant $B>0$ there is a constant $D>0$ such that $P$ has
no spectrum in
\ekv{app.12}
{
\{ z\in \mathbb{C};\, \Re z<Bh,\, |\Im z|>Dh\}
}
when $h>0$ is small enough. Moreover $\Vert (P-z)^{-1}\Vert ={\cal O}_B(h^{-1} )$ for $z$ in the set {\rm \no{app.12}}.
\end{prop}

\par Let $\rho _j\in{\cal C}$ and let $F_{\rho _j}$ be the matrix of the
linearization of $H_p$ at $\rho _j$ (the so called fundamental matrix of $p$
at the doubly characteristic point $\rho _j$). According to (\ref{ny.17}), the time average of the quadratic approximation of $\widetilde{p}$
at $\rho_j$ along the Hamilton flow of the quadratic approximation of $p_1$ at $\rho_j$ is elliptic and takes its values in a closed angle
contained in the union of $\{ 0\}$ and the open right half plane. We could therefore apply a classical result to see that
the eigenvalues of $F_{\rho _j}$ are of the form $\pm \lambda _{j,k}$, $1\le
k\le n$, when repeated with their multiplicity, with $\Im \lambda
_{j,k}>0$.

\par Put
$$
q(x,\xi )=-p(x,i\xi )=p_2+p_1-p_0.
$$
Let $F_q$, $F_p$ be the fundamental matrices of $q$, $p$ at one of the
critical points $\rho _j\in{\cal C}$. Since
$$
H_q(x,\xi )={1\over i}(p'_\xi (x,\eta )\cdot {\partial \over \partial x}-
p'_x (x,\eta )\cdot {\partial \over \partial \eta  }),\hbox{ with }\eta
=i\xi ,
$$
$F_q$ and ${1\over i}F_p$ have the same eigenvalues; $\pm {1\over
i}\lambda _k$, $k=1,...,n$ ($j$ being fixed) where $\Re ({1\over i}\lambda
_k)>0$. Now $q$ is real-valued and we can apply the stable manifold theorem
to see that the $H_q$-flow has
a stable outgoing manifold $\Lambda _+$ passing through $\rho _j$ such
that $T_{\rho _j}\Lambda _+^\mathbb{C}$ is spanned by the generalized
eigenvectors corresponding to $+{1\over i}\lambda _k$, $k=1,...,n$. We also
know that $\Lambda _+$ is a Lagrangian manifold and that $q$ vanishes on
$\Lambda _+$.

Let $\Lambda _-$ be the stable incoming $H_q$-invariant manifold such
that $T_{\rho _j}\Lambda _-^\mathbb{C}$ is spanned by the generalized
eigenvectors of $F_q$ corresponding to $-{1\over i}\lambda _k$, $1\le
k\le n$. In \cite{KFP2}, Section 8, we established a transversality
lemma for $\Lambda _\pm$ which together with a deformation argument
led to the following result,
\begin{prop}\label{as4}
We have $\phi _+''(0)>0$, $\phi _-''(0)<0$.
\end{prop}
Here $\phi_{\pm}$ is the generating function for the Lagrangian manifold $\Lambda_{\pm}$.

Let
\ekv{as.1}
{
{\widetilde{\rm tr}\,}(p,\rho _j)={1\over i}\sum_k \lambda _{j,k}.
}
In~\cite{KFP2} we also obtained the following two results. 

\begin{theo}\label{as1}
We make the assumptions {\rm \no{re.1}}--{\rm \no{re.7}}, {\rm \no{ny.1}}, {\rm \no{ny.17}}--{\rm \no{ny.20}}, and recall the definition of ${\cal C}$ in {\rm \no{ny.2}}.
Let $B>0$. Then there exists $h_0>0$ such that
for $0<h\le h_0$, the spectrum of $P$ in $D(0,Bh)$ is discrete and the
eigenvalues are of the form
\ekv{as.8}
{
{\mu } _{j,k}(h)\sim h(\mu _{j,k,{0}}+h^{1/N_{j,k}}\mu _{j,k,1}+h^{2/N_{j,k}}\mu
_{j,k,2}+...),
}
where the $\mu _{j,k,0}$ are all the numbers in $D(0,B)$ of the form
\ekv{as.9}
{\mu _{j,k,0}={1\over i}\sum_{\ell =1}^n \nu _{j,k,\ell}\lambda _{j,\ell }
+{1\over 2}{\widetilde{\rm tr}\,}(p,\rho _j), \hbox{ with }\nu _{j,k,\ell}\in{\bf N},
}
for some $j\in\{ 1,...,N\}$, $N=\# {\cal C}$. (Possibly after
changing $B$, we may assume that $|\mu _{j,k,0}|\ne B$, $\forall j,k$.) Recall
here that  $\pm\lambda _{j,\ell}$ are the eigenvalues of $F_p$ at
$\rho _j$.
This description also
takes into account the multiplicities in the natural way. If the
coefficients $\nu_{j,k,\ell} $ in {\rm \no{as.9}} are unique, then $N_{j,k}=1$ and we
have only integer powers of $h$ in the asymptotic expansion {\rm \no{as.8}}.
\end{theo}

\begin{theo}\label{as2}
We make the same assumptions as in Theorem {\rm \ref{as1}}. For all $B$, $C>0$
there is a constant $D>0$ such that
\ekv{as.10}
{
\Vert (z-P)^{-1}\Vert \le {D\over h},\hbox{ for }z\in D(0,Bh)\hbox{ with
}{\rm dist\,}(z,\sigma (P))\ge {h\over C}.
}
\end{theo}

\par Still with $j=j_0$ fixed, let
\ekv{as.12}
{
\mu_{{0}} ={1\over i}\sum_{\ell =1}^n\nu _\ell \lambda _\ell +{1\over
2}{\widetilde{\rm tr}\,}(p,\rho_{j_0}),\ \nu _\ell\in {\bf N}
}
be a value as in \no{as.9} and assume that $\mu _0 $ is simple in the
sense that $(\nu _1,...,\nu _n)\in{\bf N}^n$ is uniquely determined by
$\mu _0$. In
particular, every $\lambda _\ell$ for which $\nu _\ell\ne 0$ is a simple
eigenvalue of $F_p$. Applying a classical construction
we got
\ekv{as.13}
{
\mu (h)\sim h(\mu_0 +h\mu _1+h^2\mu _2+...)
}
with uniquely determined coefficients $\mu _1$, $\mu _2$, ... and
\ekv{as.14}
{
a(x;h)\sim a_0(x)+ha_1(x)+...\hbox{ in }C^\infty ({\rm neigh\,}(x^{j_0})),
}
where $a_j(x)={\cal O}(\vert x-x^{j_0}\vert ^{(m-2j)_+})$, $m=\sum \nu _\ell$
and $a_0$
has a non-vanishing Taylor polynomial of order $m$, such that
\ekv{as.15}
{
(P-\mu (h))(a(x;h)e^{-\phi _+(x)/h})={\cal O}(h^\infty )e^{-\phi _+(x)/h}
}
in a neighborhood of $x^{j_0}$. Actually any neighborhood $\Omega \subset\subset
\mathbb{R}^n$ will do, provided that
\smallskip
\par\noindent 1) $\phi _+$ is well-defined in a neighborhood of
$\overline{\Omega }$.\smallskip
\par\noindent 2) ${{H_q}_\vert}_{\Lambda _+}\ne 0$ on $\overline{\Omega
}\setminus \{ x^{j_0}\}$. \smallskip
\par\noindent 3) $\Omega $ is star-shaped with respect to the point $x^{j_0}$
and the integral curves of the  vector field
$\nu _+:=\left(\pi_x\right)_*({{H_q}_\vert}_{\Lambda _+})$, where $\pi_x((x,\xi))=x$.\smallskip

\par We also know that $\mu (h)$ is equal ${\rm mod\,}{\cal
O}(h^\infty )$ to the corresponding value in \no{as.8}.

\par If $\gamma \subset D(0,B)$ is a closed
$h$-independent contour avoiding all the values $\mu _{j,k,0}$ in
\no{as.9}, and \ekv{as.16} { \pi _{h\gamma} ={1\over 2\pi
    i}\int_{h\gamma }(z-P)^{-1} dz } the corresponding spectral
projection, then, \ekv{as.17} { \Vert \pi _{h\gamma }(\chi ae^{-\phi
    _+/h})-\chi ae^{-\phi _+/h}\Vert _{L^2}={\cal O}(h^\infty ) } if
$\chi \in C_0^\infty (\Omega )$ is equal to one near $x^{j_0}$. It
follows that $\chi ae^{-\phi _+/h}$ is a linear combination of
generalized eigenfunctions of $P$ with eigenvalues inside $h\gamma $
up to an error ${\cal O}(h^\infty )$ in $L^2$-norm.

\par We next review some exponential decay results from Section 9 in
\cite{KFP2} and start with the case when ${\cal C}$ is reduced to a
single point:
\ekv{exp.1} {{\cal C}=\{ (0,0)\} .}

\par Let $\phi =\phi _+(x)\in C^\infty ({\rm neigh\,}(0;\mathbb{R}))$ be the
function
introduced following Proposition \ref{as4} so that $\Lambda _+=\Lambda _\phi $ is the
stable outgoing manifold through $(0,0)$ for the $H_q$-flow. Recall that by
Proposition \ref{as4}
\ekv{exp.5}
{
\phi ''(0)>0.
}
Moreover, there exists $G\in C^\infty ({\rm neigh\,}(0,\mathbb{R}^n);\mathbb{R})$
such that
\ekv{exp.7}
{
(\partial _\xi q)(x,\phi '(x))\cdot \partial _xG\asymp x^2,\ G(x)\asymp x^2.
}

\par Let $\Omega _G(r)=\{ x\in{\rm neigh\,}(0);\, G(x)\le r\}$ for $0<r\ll
1$. Outside the set $\Omega _G(C_0\epsilon )$, we put
\ekv{exp.8}
{\widehat{\psi }=\phi -\epsilon g(G),}
where $g=g(G)=\ln G$ for $G\ge C_0 \epsilon $, so that
$g'(G)=1/G$. Then
\ekv{exp.11}
{
q(x,\widehat{\psi }'(x))\le -{\epsilon \over C_0},\ x\in {\rm
neigh\,}(0,M)\setminus  \Omega _G(C_0\epsilon ),}
if $C_0>0$ is large enough.

We extend the definition of $\psi $ to a full
neighborhood of $x=0$, by putting
\ekv{exp.14}
{
g(G)=\ln (C_0\epsilon )+{1\over C_0\epsilon }(G-C_0\epsilon ),\hbox{ for }
0\le G\le C_0\epsilon .
}

\par Outside a small fixed neighborhood of 0 we flatten out the
weight. Let $f_\delta (t)=\delta f({t\over \delta })$ be the function
introduced before (\ref{re.8}). For some small and fixed $\delta _0>0$, we
put
\ekv{exp.12}
{
\psi =f_{\delta _0}(\widehat{\psi })=f_{\delta _0}(\phi -\epsilon g(G))
}
which is also well-defined as the constant $3\delta _0/2$ for large
$x$. Finally, by the same averaging procedure as in the proof of
Proposition \ref{re1}, we can add a term ${\cal O}(\epsilon )$ to
$\psi $, supported away from a fixed neighborhood of $0$, such that if
$\psi _\epsilon $ denotes the corresponding modification of $\psi $,
we have the apriori estimate,
\ekv{exp.26}
{
h\Vert e^{\psi _\epsilon /h}v\Vert +h^{1\over 2}\Vert B^{1\over
2}hD(e^{\psi _\epsilon /h}v)\Vert \le {\cal O}(1)\Vert e^{\psi _\epsilon
/h}(P-z)v\Vert +{\cal O}(h)\Vert e^{\psi _\epsilon /h}v\Vert _{\Omega
_G(C_0\epsilon )},
}
uniformly, for $\vert \Re z\vert \le Ch$ provided that $\epsilon =Ah$ for
$A$ large enough depending on $C$.

\par Now let $\mu (h)=\mu _{1,k}(h)$ be an eigenvalue of $P$ as in
\no{as.8}, \no{as.13} and assume that $\mu_0 $ is given by \no{as.12}
and is simple, as explained after that equation. Then $\mu (h)$ is a
simple eigenvalue of $P$ and is the only eigenvalue in some disc
$D(\mu (h),h/C_0)$. Let $u_{\rm WKB}(x;h)$ be the approximate
solution given in \no{as.14}, \no{as.15} and let $u=\pi _{h\gamma
}(\chi u_{\rm BKW})$ be the corresponding exact eigenfunction, where
$\gamma =\partial D(\mu ,{1\over 2C_0})$. Using (\ref{exp.26}) we
established the following result in \cite{KFP2}:
\begin{theo}\label{exp1}
\begin{enumerate}[a)]
\item Outside any $h$-independent neighborhood of 0, we have
$$
u,\, B^{1\over 2}hDu={\cal O}(e^{-1/(Ch)})
$$
in $L^2$-norm.
\item There exists a neighborhood $\Omega $ of $0$, where
\ekv{exp.27}
{\begin{split}
&u(x;h)=(a+r)e^{-\phi _+(x)/h},\\ &\Vert r\Vert _{L^2(\Omega )},\, \Vert
B^{1\over 2}hDr\Vert _{L^2(\Omega )}={\cal O}(h^\infty ). \end{split}
}
\end{enumerate}
\end{theo}
\begin{remark}\label{exp2} \rm
If we drop the assumption (\ref{exp.1}) and allow $N-1$ more points $\rho
_2,...,\rho _N$ in ${\cal C}$, then Theorem \ref{exp1} is still valid,
provided that all $\mu _{j,k,0}$ in \ref{as.8} with $j\ge 2$ are
different from the value $\mu_0 $, associated to $\rho _1=(0,0)$.
\end{remark}

\medskip
We now drop the assumption (\ref{exp.1}) completely, so that
${\cal C}=\{ \rho _1,...,\rho _N\}$, $\rho _j=(x^j,0)$, $\rho
_1=(0,0)$. We then have the following extension of Theorem \ref{exp1}.
\begin{theo}\label{exp3}
We make the assumptions of Theorem {\rm \ref{exp1}}, with the exception of
{\rm (\ref{exp.1})}. Let $u_{\mathrm{WKB}}(x)$ be as in that theorem and as
there, we put $u=\pi _{h\gamma }(u_{\mathrm{WKB}})$, with $C_0$ large
enough in the definition of $\gamma $ above, so that $\mu (h)$ is
the only asymptotic eigenvalue of $P$ inside $h\gamma $ that we can
associate with the critical point $\rho _1=(0,0)$. Then $u$ is not
necessarily an eigenfunction of $P$ (due to possible eigenvalues of
$P$ inside a disc $D(\mu (h),o(h))$, associated to other points in
${\cal C}$), but the conclusions a) and b) of Theorem {\rm \ref{exp1}}
remain valid.
\end{theo}

This result was not stated in \cite{KFP2}, but follows fairly directly
from Theorem \ref{exp1} and the earlier results by some easy and
standard arguments:
\begin{itemize}
\item Let us first ``eliminate'' ${\cal C}\setminus \{(0,0)\}$, by
  introducing $\widetilde{P}:=P+\sum_2^N \alpha \chi
  (\frac{x-x^j}{\alpha })$ for some small and fixed $\alpha >0$. Here $\chi$ is a standard cut-off function to a small
  neighborhood of $0$.
\item Then $\widetilde{P}$ fulfills the assumption (\ref{exp.1}), so
  the conclusions a) and b) of Theorem \ref{exp1} apply to
  $\widetilde{u}=\widetilde{\pi }_{h\gamma }(\chi u_{\mathrm{WKB}})$,
  where $\widetilde{\pi }_{h\gamma }$ is the spectral projection
  associated to $\widetilde{P},\,h\gamma $.
\item The resolvent identity
  $(z-P)^{-1}=(z-\widetilde{P})^{-1}+(z-P)^{-1}(P-\widetilde{P})(z-\widetilde{P})^{-1}$
  implies that
\ekv{exp.28}
{
u=\widetilde{u}+\frac{1}{2\pi i}\int_{h\gamma
}(z-P)^{-1}(P-\widetilde{P})(z-\widetilde{P})^{-1}\chi u_{\mathrm{WKB}}dz.
}
\item Using the a priori estimate (\ref{exp.26}), with $P$ replaced by
  $\widetilde{P}$ and $v=(z-\widetilde{P})^{-1}\chi u_{\mathrm{WKB}}$,
  together with Theorem \ref{as2} and Corollary \ref{re3},
  we first see that $(P-\widetilde{P})(z-\widetilde{P})^{-1}\chi
  u_{\mathrm{WKB}}$ is exponentially decaying in $L^2$, uniformly for
  $z\in h\gamma $. Using Theorem \ref{as2} and Corollary \ref{re3}
  once more, we next get that
  $\widehat{u}:=(z-P)^{-1}(P-\widetilde{P})(z-\widetilde{P})^{-1}\chi
  u_{WKB}$ is uniformly exponentially decreasing in ${\cal D}(P)$. In
  particular, $\Vert \widehat{u}\Vert+\Vert
  B^{\frac{1}{2}}hD\widehat{u}\Vert ={\cal O}(e^{-1/Ch})$ and we have
  the same estimate with $\widehat{u}$ replaced by the integral in
  (\ref{exp.28}).
\item Theorem \ref{exp3} now follows by combining the above facts for
the two terms in (\ref{exp.28}).

\end{itemize}

\subsection{The supersymmetric case}\label{ss}

\par The Witten approach has been
independently extended to the case of non-elliptic operators like the
Kramers-Fokker-Planck operator in \cite{TaTaKu06} (in supersymmetric
language) and in \cite{Bi05} (in terms of differential forms). See
also \cite{Leb04}. We here follow the presentation of Section 10 and 11 of
\cite{KFP2} and refer to that work for more details and proofs.

\par We start by a quick review of that in the semiclassical case,
then we establish some basic facts about the principal and
subprincipal symbols, especially at the critical points of the given
weight function. Finally we add some growth conditions and a dynamical
condition, so that the results of Subsection \ref{reg} can be applied.

\subsubsection{Generalities}\label{ssa}
Let
\ekv{ss.1}
{
A:\,(\mathbb{R}^n)^*\to \mathbb{R}^n
}
be a linear and invertible map. Then we have the real nondegenerate bilinear form
\ekv{ss.2}
{
( u\vert v) _A=\langle \wedge ^kA(u)\vert v\rangle ,\
u,v\in \wedge ^k(\mathbb{R}^n)^*.
}
If $a:\wedge ^k(\mathbb{R}^n)^*\to \wedge ^\ell (\mathbb{R}^n)^*$ is a linear map, we define
the "adjoint" $a^{A,*}:\wedge ^\ell (\mathbb{R}^n)^*\to \wedge ^k(\mathbb{R}^n)^*$ by
\ekv{ss.3}
{
( au\vert v) _A=( u\vert a^{A,*}v) _A.
}
(In the complexified case, we { take the sesquilinear
  extension of $(u\vert
v)_A$} and define $a^{A,*}$ the same
way.)

\par If $\omega $
is a one form, we have
\ekv{ss.5}
{ (\omega ^\wedge )^{A,*}=(A\omega)^\rfloor,  } where ${}^\wedge$ and ${}^\rfloor$ denote
the usual operators of (left) exterior product and contraction.

\par If $u,v$ are smooth $k$ forms with $\supp u\cap\supp v$ compact, we
define
$$
( u\vert v) _A=\int ( u(x)\vert v(x)) _{A}dx$$ and denote
by $a^{A,*}$ the formal adjoint of an operator $a:C_0^\infty (\mathbb{R}^n;\wedge
^kT^*\mathbb{R}^n)\to {\cal D}'(\mathbb{R}^n;\wedge ^\ell T^*\mathbb{R}^n)$.  We can consider
$$
\partial _{x_j}:C_0^\infty (\mathbb{R}^n;\wedge ^kT^*\mathbb{R}^n)\to
C_0^\infty (\mathbb{R}^n;\wedge ^kT^*\mathbb{R}^n),
$$
acting coefficient-wise, and a straightforward computation shows that
\ekv{ss.7}
{
(h\partial _{x_j})^{A,*}= -h\partial _{x_j}.
}

\par
Let $\phi \in C^\infty (\mathbb{R}^n;\mathbb{R})$
and introduce the Witten (de Rham) complex
\ekv{ss.8}
{
d_\phi =e^{-\phi /h}\circ hd\circ e^{\phi /h}=hd+(d\phi )^\wedge :\,
C_0^\infty (\mathbb{R}^n;\wedge ^kT^*\mathbb{R}^n)\to C_0^\infty (\mathbb{R}^n;\wedge ^{k+1}T^*\mathbb{R}^n),
}
with $d_\phi ^2=0$.
We have
\ekv{ss.9}
{d_\phi =\sum_1^n(h\partial _{x_j}+\partial _{x_j}\phi )\circ dx_j^\wedge ,}
where $h\partial _{x_j}+\partial _{x_j}\phi $ acts coefficient-wise and
commutes with $dx_j^\wedge $,
so
\ekv{ss.10}
{
d_\phi ^{A,*}=\sum_1^n (-h\partial _{x_k}+\partial _{x_k}\phi )\circ A(dx_k)^\rfloor .
}

\par
The corresponding Witten-Hodge Laplacian is given by
\ekv{ss.11}
{-\Delta _A=d_\phi ^{A,*}d_\phi +d_\phi d_\phi ^{A,*},}
and we have
\ekv{ss.12}
{\begin{split}
-\Delta _A=& \sum_{j,k}(-h\partial _{x_k}+\partial _{x_k}\phi
)A_{j,k}(h\partial _{x_j}+\partial _{x_j}\phi )
\\ &+\sum_{j,k}2h\partial _{x_j}\partial _{x_k}\phi \circ dx_j^\wedge
A(dx_k)^\rfloor . \end{split}
}

\par
Now write
\ekv{ss.13}
{
A=B+C,\ {B}^{\mathrm{t}}=B,\ {C}^{\mathrm{t}}=-C.
}
Then \no{ss.12} gives
\ekv{ss.14}
{\begin{split}
-\Delta _A = & \sum_{j,k}(-h\partial _{x_k}+\partial _{x_k}\phi
)B_{j,k}(h\partial _{x_j}+\partial _{x_j}\phi )
\\ &+\sum_{j,k}((\partial _{x_k}\phi )C_{j,k}h\partial _{x_j}+h\partial
_{x_j}\circ C_{j,k}\circ (\partial _{x_k}\phi ))
\\ &+\sum_{j,k}2h\partial _{x_j}\partial _{x_k}\phi \circ dx_j^\wedge
A(dx_k)^\rfloor . \end{split}
}
Note that the last term vanishes on $0$-forms, i.e. on scalar
functions. 
To recover the Kramers-Fokker-Planck operator (cf \cite{TaTaKu06}), replace $n$ by $2n$,
$$
A={1\over 2}\begin{pmatrix}0 & 1\cr -1 &\gamma \end{pmatrix}.
$$
Then \no{ss.14} gives for 0-forms:
\ekv{ss.15}
{
-\Delta _A^{(0)} = {\gamma \over 2}\sum_{j=1}^n (-h\partial _{y_j}+\partial
_{y_j}\phi )(h\partial _{y_j}+\partial
_{y_j}\phi )+hH_\phi ,
}
where
$$
H_\phi =\sum (\partial _{y_k}\phi \,  \partial _{x_k}-\partial _{x_k}\phi \,
\partial _{y_k})
$$
is the Hamilton field of $\phi $ with respect to  the standard symplectic form
$\sum dy_j\wedge dx_j$.

\par If we choose
\ekv{ss.16}
{
\phi (x,y)={1\over 2}y^2+V(x),
}
we get the  Kramers-Fokker-Planck operator (\ref{in.0}),
\ekv{ss.17}
{-\Delta _A^{(0)}=y\cdot h\partial _x-V'(x)\cdot h\partial _y+{\gamma
\over 2}(-h\partial _y+y)\cdot (h\partial _y+y).
}

\subsubsection{The principal symbol of the Hodge Laplacian}\label{ssb}
\par The principal symbol of $-\Delta _A$ in the
sense of $h$-differential operators
is scalar and given by
\ekv{ss.18}
{\begin{split}
p(x,\xi )=&\sum_{j,k}A_{j,k}(-i\xi _k+\partial _{x_k}\phi )
(i\xi _j+\partial _{x_j}\phi )
\\=&\sum_{j,k}B_{j,k}(\xi _j\xi _k+\partial _{x_j}\phi \,  \partial _{x_k}\phi
)+2i\sum_{j,k}C_{j,k}\partial _{x_k}\phi \,  \xi _j. \end{split}}
The corresponding real symbol $q(x,\xi )=-p(x,i\xi )$ is given by
\ekv{ss.19}
{\begin{split}
q(x,\xi )=&\sum_{j,k}A_{j,k}(\xi _k+\partial _{x_k}\phi )(\xi _j-\partial
_{x_j}\phi )
\\=&\sum_{j,k}B_{j,k}(\xi _j\xi _k-\partial _{x_j}\phi \,  \partial
_{x_k}\phi )+2\sum_{j,k}C_{j,k}\partial _{x_k}\phi \,  \xi _j.
\end{split}}
It vanishes on the two Lagrangian manifolds $\Lambda _{\pm \phi }$.

\par We define
\ekv{ss.20}
{
\nu _\pm={{H_q}_\vert}_{\Lambda _{\pm \phi }}.
}
Using $x_1,...,x_n$ as coordinates on $\Lambda _{\pm \phi }$, we get
\ekv{ss.21}
{
\nu _+=2\sum_{j,k}A_{j,k}\partial _{x_k}\phi \,  \partial _{x_j}=2A(\phi
'(x))\cdot \partial _x
}
\ekv{ss.22}
{
\nu _-=-2\sum_{j,k}A_{j,k}\partial _{x_j}\phi \, \partial _{x_k}=
-2\,{A}^{\mathrm{t}}(\phi'(x))\cdot \partial _x .
}

\par Let $x_0$ be a nondegenerate critical point of $\phi $, so that
$\Lambda _\phi $ and $\Lambda _{-\phi }$ intersect transversally at
$(x_0,0)$.  The spectrum of the linearization $F_q$ of $H_q$ at
$(x_0,0)$ is equal to the union of the spectra of the linearizations
\ekv{ss.23} {\nu _+^0=(2A\phi ''(x_0)x)\cdot \partial _x\hbox{ and }
  \nu _-^0=-(2{A}^{\mathrm{t}}\phi ''(x_0)x)\cdot \partial _x } of $\nu _+$
and $\nu _-$ respectively at $x_0$. Thus we are interested in the
eigenvalues of the matrices $A\phi ''$, ${A}^{\mathrm{t}}\phi ''$.  Here,
${A}^{\mathrm{t}}\phi ''={\phi ''}^{-1}{(A\phi '')^{\mathrm{t}}\phi ''}$ has the
same eigenvalues as $A\phi ''$, and similarly $\phi ''A$, $\phi''{A}^{\mathrm{t}}$ are isospectral to $A\phi ''$.
Thus \ekv{ss.24} {\begin{split}
  &\hbox{The eigenvalues of }F_q \hbox{ are given by }\pm 2\lambda
  _j,\\ &\hbox{where }\lambda _1,...,\lambda _n\hbox{ are the
    eigenvalues of }A\phi ''. \end{split}  } (Here the notation is different from
the one used prior to Proposition \ref{as4}.)

\par We assume from now on that
\ekv{ss.25}
{
B\ge 0.
}
In \cite{KFP2} we established the following result:

\begin{lemma}\label{ss1}
Let $\mu (x,\partial _x)=Mx\cdot \partial _x$ be a real linear vector field on
$\mathbb{R}^n$. Let $n_{\pm}\in{\bf N}$, $n_++n_-=n$. Then the following two
statements are equivalent:
\smallskip
\par\noindent {\rm (A)} M has $n_+$ eigenvalues with real part $>0$ and $n_-$ eigenvalues
with real part $<0$.
\par\noindent {\rm (B)} There exists a quadratic form $G:\mathbb{R}^n\to \mathbb{R}$
of signature $(n_+,n_-)$ and a constant $C>0$, such that
\ekv{ss.27}
{
\mu (x,\partial _x)(G)\ge {1\over C}\vert x\vert ^2,\ x\in\mathbb{R}^n.
}
\end{lemma}
Using this lemma we got the following proposition at the nondegenerate
critical point $x_0$ that we assume to be zero for notational reasons. Here $p^0$ is the quadratic approximation of $p$ at $(0,0)$.
\begin{prop}\label{ss2} {\rm a)} Assume that the matrix $A\phi ''$ of $\nu _+^0$
has $m_{\pm}$ eigenvalues with $\pm$ real part $>0$, $m_++m_-=n$. Then there
exists a real quadratic form ${\cal G}(x,\xi )$ on $\mathbb{R}^{2n}$ such that
\ekv{ss.31}
{
\Re p^0((x,\xi )+i\epsilon H_{\cal G}(x,\xi ))\ge {\epsilon \over C}\vert
(x,\xi )\vert ^2,\ (x,\xi )\in\mathbb{R}^{2n},\,\, 0<\epsilon \ll 1.
}
\smallskip
\par\noindent {\rm b)} Conversely, assume that there exists a quadratic form
${\cal G}$ such that {\rm \no{ss.31}} holds. Then $A\phi ''$ has $n_\pm$ eigenvalues
with $\pm$ real part $>0$, where $(n_+,n_-)$ is the signature of $\phi ''(0)$.
\end{prop}

\par The condition \no{ny.17} implies the existence of ${\cal G}$
as in a) of the proposition.

\subsubsection{The subprincipal symbol}\label{ssc}

\par We next look at the subprincipal term in \no{ss.12}, i.e. the
last sum in that equation. As we saw in \cite{KFP2}, it can be
rewritten as
\ekv{ss.35} { 2h\sum _j (\phi ''\circ {A}^{\mathrm{t}})(dx_j )
  ^\wedge \partial _{x_j }^\rfloor .  }

\par Now we restrict the attention to a nondegenerate critical point $x_0$
of $\phi $ and we shall compute the subprincipal symbol of $-\Delta _A$
at the corresponding doubly characteristic point $(x_0,0)$. At that
point $\phi ''\circ {A}^{\mathrm{t}}:T_{x_0}^*\mathbb{R}^n\to T_{x_0}^*\mathbb{R}^n$ is invariantly
defined and it is easy to check that \no{ss.35} is also invariantly
defined: we get the same quantity if we replace $dx_1,..,dx_n$, $\partial
_{x_1},..,\partial _{x_n}$, by $\omega _1,..,\omega _n$, $\omega
_1^*,..,\omega _n^*$, where $\omega _1,..,\omega _n$ is any basis in the
complexified cotangent space and $\omega _1^*,..,\omega _n^*$ is the dual
basis of tangent vectors for the natural bilinear pairing.

\par Assume that the equivalent conditions of Proposition \ref{ss2} hold
and denote the corresponding
eigenvalues (that are also the eigenvalues of $\phi ''\circ {A}^{\mathrm{t}}$) by $\lambda
_1,..,\lambda _n$ with $\Re \lambda _j>0$
for $1\le j\le n_+$ and with $\Re \lambda _j <0$ for $n_++1\le j\le
n=n_++n_-$. The eigenvalues of $F_p$ are then $\pm 2i\lambda _j$
(in view of \no{ss.24}
and the isospectrality of $F_p$ and $iF_q$ reviewed prior to
Proposition
\ref{as4})
, so
\ekv{ss.36}
{
\widetilde{{\rm tr\,}}F_p:={1\over i}\sum_{\mu \in \sigma (F_p)\atop \Im
\mu >0}\mu =\sum_1^{n_+}2\lambda _j-\sum_{n_++1}^n2\lambda _j .
}

\par The subprincipal symbol of the first term in \no{ss.12} (at
$(x_0,0)$) is equal to \ekv{ss.37} { \sum_{j,k}A_{j,k}{1\over 2i}\{
  -i\xi _k+\partial _{x_k}\phi ,i\xi _j+\partial _{x_j}\phi \}
  =-\sum_{j,k}A_{j,k}\phi ''_{j,k}=-\mathrm{tr\,}(A\phi '')=-\sum_1^n
  \lambda _j.  } {Simplifying the last term in
  (\ref{ss.12}), we get the full subprincipal symbol at $(x_0,0)$:
$$
S_P=-\sum_1^n\lambda _j+2\sum_j ((\phi ''\tA )(dx_j))^\wedge \partial _{x_j}^\rfloor
$$
and hence 
$$
\frac{1}{2}\widetilde{\mathrm{tr}\,}F_p+S_P=-2\sum_{n_++1}^n\lambda _j+2\sum_j ((\phi ''\tA )(dx_j))^\wedge \partial _{x_j}^\rfloor
$$}
\par The eigenvalues of $\sum_j (\phi ''\circ A^{\mathrm{t}})(dx_j)^\wedge \partial
_{x_j}^\rfloor$ on the space of $m$-forms are easily calculated, if we
replace $dx_1,..,dx_n$ by a basis of eigenvectors $\omega _1,..,\omega _n$
of $\phi ''\circ A^{\mathrm{t}}$, so that
$$
(\phi ''\circ {A}^{\mathrm{t}})(\omega _j)=\lambda _j\omega _j,
$$
and $\partial _{x_j}$ by the corresponding dual basis vectors $\omega _j^*$. (Here we assume to
start with that there are no Jordan blocks. This can be achieved by an
arbitrarily small perturbation of $A$, and we can extend the end result of our
calculation to the general case by continuity.) We get
\ekv{ss.38}
{
\sum _j (\phi ''\circ A^{\mathrm{t}})(dx_j )^\wedge \partial _{x_j
}^\rfloor=\sum_j \lambda _j\omega _j^\wedge  {\omega _j^*}^\rfloor .
}
A basis of eigenforms of this operator is given by $\omega _{j_1}\wedge
..\wedge \omega _{j_m}$, $1\le j_1<j_2<..<j_m\le n$ and the corresponding
eigenvalues are $\lambda _{j_1}+..+\lambda _{j_m}$.

\par Then the eigenvalues of
$$
{1\over 2}\widetilde{\mathrm{tr\,}}F_p+S_P,\hbox{ acting on }m\hbox{ forms }
$$
are \ekv{ss.39} 
{2(\lambda _{j_1}+..+\lambda _{j_m}-\sum_{n_++1}^n\lambda _j),\ 1\le
  j_1<..<j_m\le n.}  We conclude that if $m\ne n_-$, then all the
eigenvalues have a real part $>0$. If $m=n_-$, then precisely one
eigenvalue is equal to 0, while the others have positive real part  and the corresponding one
  dimensional kernel is spanned by $\omega _{n_++1}\wedge ... \wedge
  \omega _{n}$.

\subsubsection{A symmetry for adjoints}\label{ssd}
\par In this subsection we are concerned with symmetry relations for
the $A,*$ adjoints. If $D:L^2(\Omega ;\wedge ^kT^*\Omega )\to L^2(\Omega
;\wedge ^jT^*\Omega )$, then a simple calculation shows that
$$
D=(D^{A,*})^{{A}^{\mathrm{t}},*}.
$$
This can be applied to $-\Delta _A$ and we get
\ekv{su.6}{
(-\Delta _A)^{{A}^{\mathrm{t}},*}=-\Delta _{{A}^{\mathrm{t}}},\ (-\Delta
_{{A}^{\mathrm{t}}})^{A,*}=-\Delta _A.
}

\subsubsection{Quasimodes and spectral subspaces}\label{qms} So far,
we only developed the formal aspects of the supersymmetric
approach. Now assume,
\begin{hyp}\label{hyppsa}
$\phi :\mathbb{R}^n\to \mathbb{R}$ is smooth and satisfies
\begin{equation}
\label{eq3.01}
\partial_x^{\alpha}\phi(x)={\cal O}(1),\quad
\partial_x^{\alpha}\left(\langle{B\partial_x\phi,\partial_x\phi\rangle}\right)=
{\cal O}(1),\quad \abs{\alpha}\geq 2.
\end{equation}
Moreover, $\phi $ is a Morse function and
there exists a constant $C>0$ such that
\begin{equation} \label{morse}
\abs{\phi'(x)}\geq 1/C ,\quad \abs{x}\geq C.
\end{equation}
\end{hyp}

We also assume from now on that $A$ in (\ref{ss.13}) satisfies
\ekv{su.1}{B\ge 0.}

\par Then $P=-\Delta _A^{(0)}$ satisfies the assumptions
(\ref{re.4})--(\ref{re.7}) and Hypothesis \ref{ny0} is fulfilled with
$$
{\cal C}=\{ (x^1,0),...,(x^N,0)\} ,
$$
where $x^j$ denote the critical points of $\phi $. Indeed,
(\ref{ss.18}) can also be written
\begin{equation}
\label{eq3.1}
p(x,\xi)=\langle{B\xi|\xi}\rangle+2i\langle{C\phi'_x|\xi}\rangle+\langle{B\phi'_x|\phi'_x}\rangle,
\end{equation}
and the vector field $\nu $ becomes
$$
\nu =2C\phi '_x\cdot \frac{\partial }{\partial x}.
$$
This means that Proposition \ref{re2} and Corollary \ref{re3} apply to
$P=-\Delta _A^{(0)}$. Notice also that the principal part of $-\Delta
_A^{(\ell )}$ is scalar so that this operator differs from the tensor
product of $-\Delta _A^{(0)}$ with the identity in some
$\mathbb{C}^{N_\ell}$ up a zero order term which is ${\cal O}(h)$ in
$C_b^\infty $. Hence the above mentioned results apply to $-\Delta_A^{(\ell )}$ as well, provided that we strengthen the assumption on $z$
in Proposition \ref{re2} to $\Re z < -Ch$, for $C\geq 0$ suitably chosen.  

We now adopt the dynamical assumption, Hypothesis \ref{ny2}. Then
Propositions \ref{weest1}, \ref{app2}, \ref{as4}, Theorem \ref{as1}, \ref{as2},
\ref{exp3} apply also.

\begin{prop}\label{ir1}
For $z\in \mathbb{C}$, $\Re z\le {\cal O}(h)$, $v\in {\cal S}$, we have
\begin{equation}\label{ir.1}\begin{split}
&(-\Delta _A^{(\ell +1)}-z)^{-1}d_\phi v=d_\phi (-\Delta
_A^{(\ell )}-z)^{-1}v,\\
&\hbox{when }z\not\in \sigma (-\Delta
_A^{(\ell )})\cup \sigma (-\Delta _A^{(\ell +1)}),\end{split}
\end{equation}
\begin{equation}\label{ir.2}\begin{split}
&(-\Delta _A^{(\ell -1)}-z)^{-1}d_\phi^{A,*} v=d_\phi^{A,*} (-\Delta
_A^{(\ell )}-z)^{-1}v,\\
&\hbox{when }z\not\in \sigma (-\Delta
_A^{(\ell -1 )})\cup \sigma (-\Delta _A^{(\ell )}).\end{split}
\end{equation}
\end{prop}

\par\noindent {\it Proof.} By unique holomorphic continuation it suffices to
establish these relations for $-\Re z\gg h^2$.

\par For such values of $z$, we can apply Proposition \ref{weest1} to
see that $\Lambda (-\Delta _A^{(\ell )}-z)^{-1}v\in L^2$ with $\Lambda
=\mathrm{Op}_h(\langle \rho \rangle)$. In particular, $d_\phi (-\Delta
_A^{(\ell )}-z)^{-1}v\in L^2$. Let $u=(-\Delta _A^{(\ell )}-z)^{-1}v$,
so that $u,d_\phi u\in L^2$.

From $(-\Delta _A^{(\ell )}-z)u=v$ we get, using the intertwining
relations;
$$
d_\phi v=d_\phi (-\Delta _A^{(\ell )}-z)u=(-\Delta _A^{(\ell
  +1)}-z)d_\phi u.
$$
By the equality of the minimal and maximal closed extensions (Corollary \ref{re3}), and the fact that $d_\phi v, d_\phi u\in L^2$, we
get
$$
d_\phi u=(-\Delta _A^{(\ell +1)}-z)^{-1}d_\phi v,
$$
and we get (\ref{ir.1}). The proof of (\ref{ir.2}) is
similar.\hfill{$\Box$}

\medskip
\par
We shall now discuss the action of $d_\phi$ and $d_\phi ^{A,*}$ on generalized eigenspaces of $-\Delta_A^{(\ell )}$.
Let $\gamma \subset \mathbb{C}$ be a simple closed contour such that no eigenvalues of the quadratic approximations of the
operators $-\Delta _A^{(\ell)}$ at the critical points of $\phi$ for $h=1$, $\ell =0,...,n$, belong to $\gamma $.

\par Then $h\gamma \cap \sigma (-\Delta _A^{(\ell )})=\emptyset $ for
$0<h\ll 1$, so we can introduce the spectral projections,
$$
\Pi ^{(\ell)}=\frac{1}{2\pi i}\int_{h\gamma }(z+\Delta _A^{(\ell)})^{-1}dz,
$$
and their finite dimensional ranges,
$$
E^{(\ell)}=\Pi ^{(\ell )}(L^2),
$$
Since ${\cal S}$ is dense in $L^2$, we can replace $L^2$ by ${\cal S}$
in the definition of $E^{(\ell)}$ and Proposition \ref{ir1} tells us
that $\Pi ^{(\ell +1)}d_\phi =d_\phi \Pi ^{(\ell)}$ on
${\cal S}$. Consequently,
\begin{equation}\label{ac.1}
d_\phi :\, E^{(\ell)}\to E^{(\ell +1)}.
\end{equation}
In fact, if $u\in E^{(\ell)}$, we can write $u=\Pi ^{(\ell
  )}(\widetilde{u})$, $\widetilde{u}\in {\cal S}$, and then
$d_\phi u=d_\phi \Pi ^{(\ell)}(\widetilde{u})=\Pi ^{(\ell +1)}d_\phi
\widetilde{u}\in E^{(\ell +1)}$.

\par Similarly,
\begin{equation}\label{ac.2}
d_\phi ^{A,*}:\, E^{(\ell +1)}\to E^{(\ell)}.
\end{equation}

In what follows, we shall mainly consider $-\Delta_A^{(0)}$ and $-\Delta _A^{(1)}$. Let
$n_0$ be the number of local minima, $m_1,...,m_{n_0}$ of $\phi
$. Then, if $\chi \in C_0^\infty (\mathbb{R}^n)$ is a standard cut-off
to a small neighborhood of $0$, the functions
\ekv{qms.1}{f_k^{(0)}=h^{-n/4}\chi (x-m_k)e^{-(\phi -\phi (m_k))/h},\
  k=1,...,n_0,} are quasimodes of $-\Delta _A^{(0)}$ with eigenvalue 0
in the sense that
\ekv{qms.2}{\begin{split}&\Vert f_k^{(0)}\Vert\asymp 1,\ (f_k^{(0)}|f_\ell^{(0)})=0
  \hbox{ for }k\ne \ell,\\ &-\Delta _A^{(0)}f_k^{(0)}={\cal O}(e^{-1/Ch})\hbox{ in }L^2. \end{split}}
Using also the calculation for the subprincipal symbol above, we get,
(as in \cite{KFP2})
\begin{prop}\label{qms1}
If $C>0$ is large enough, then $-\Delta _A^{(0)}$ has precisely $n_0$
eigenvalues in $D(0,h/C)$ (for $h$ small enough) when counted with
their multiplicity. These eigenvalues are actually ${\cal O}(h^\infty )$.
\end{prop}

Let $E^{(0)}$ be the spectral subspace corresponding to the
eigenvalues in the proposition and let $\Pi ^{(0)}$ denote the
corresponding spectral projection. Then $e_k^{(0)}:=\Pi
^{(0)}f_k^{(0)}$ form a basis in $E^{(0)}$ such that
\ekv{qms.3}
{
f_k^{(0)}-e_k^{(0)},\ B^{1/2}hD(f_k^{(0)}-e_k^{(0)})={\cal O}(e^{-1/Ch})\hbox{ in
  }L^2.
}
This follows from Theorem \ref{exp3} and its proof. Later we shall
choose the quasimodes $f_k^{(0)}$ more carefully, using more refined
cut-offs.

Let $n_1$ be the number of saddle points, i.e. critical points of
index 1.
\begin{prop}\label{qms2}
If $C>0$ is large enough, then for $h>0$ sufficiently small, $-\Delta
_A^{(1)}$ has precisely $n_1$ eigenvalues in $D(0,h/C)$ (counted with
their multiplicity). The corresponding eigenvalues are ${\cal
  O}(h^\infty )$. We denote by $E^{(1)}\subset L^2$ the corresponding spectral subspace. 
\end{prop}

The first part follows from Theorem \ref{exp3} and the calculation of
the subprincipal symbol (see \cite{KFP2} for more details). That the
corresponding eigenvalues are ${\cal O}(h^\infty )$ and not just
${\cal O}(h^2)$ is also the consequence of a standard argument:

\par Let
\ekv{qms.4}
{
f_j^{(1)}=h^{-n/4}a_j(x;h)e^{-\phi _+(x)/h},
}
where
$$
a_j(x;h)\sim a_{j,0}(x)+ha_{j,1}+... \hbox{ in }C^\infty (\mathrm{neigh\,}(s_j,\mathbb{R}^n)),\quad a_{j,0}(s_j)\neq 0,
$$
and $\phi _+\asymp |x-s_j|^2$ solves the eiconal equation $q(x,\phi'_+)=0$ in a neighborhood of $s_j$. We assume that $f_j^{(1)}$ is a
quasimode so that in the sense of formal asymptotic expansions,
$$
-\Delta _A^{(1)}f_j^{(1)}=\mu _j(h)f_j^{(1)},\ \mu _j(h)={\cal
  O}(h^2).
$$
Using the intertwining relations
\ekv{qms.5}
{
\Delta _A^{(0)}d_\phi ^{A,*}=d_\phi ^{A,*}\Delta _A^{(1)},\ \Delta
_A^{(2)}d_\phi =d_\phi \Delta _A^{(1)},
}
we see that $d_\phi ^{A,*}f_j^{(1)}$ and $d_\phi f_j^{(1)}$ are
quasimodes to $-\Delta _A^{(0)}$ and $-\Delta _A^{(2)}$ respectively,
with the same eigenvalue $\mu _j=o(h)$. However, from Theorem
\ref{exp3} and the calculations of the subprincipal symbol above, we
know that these two operators cannot have any non-trivial quasimodes at
$s_j$ with an eigenvalue $o(h)$. Hence $d_\phi ^{A,*}f_j^{(1)}$ and
$d_\phi f_j^{(1)}$ vanish in the space of asymptotic WKB expressions,
and using that $-\Delta _A^{(1)}=d_\phi ^{A,*}d_\phi + d_\phi d_\phi
^{A,*}$, we see that $\mu _j={\cal O}(h^{\infty })$, as claimed.

\par Recall that $\Lambda _{\phi _+,j}$ is the stable outgoing
manifold through $(s_j,0)$ for the $H_q$-flow and that $\phi
_{+,j}''(s_j)>0$ by Proposition \ref{as4}. (Similarly we
have a stable incoming manifold $\Lambda _{\phi _{-,j}}$.) Recall also
(see Proposition \ref{ss2}) that the linearization of
${{H_q}_\vert }_{\Lambda _\phi }$ at that point has $n-1$ positive
eigenvalues and $1$ negative eigenvalue. Its matrix is $2A\phi
''(s_j)$, according to (\ref{ss.21}). We let $K_+, K_-\subset \Lambda _\phi
$ be the corresponding stable outgoing and incoming submanifolds of
dimension $k_+=n-1$ and $k_-=1$ respectively and recall that
$K_+\subset \Lambda _{\phi _+}$, $K_-\subset \Lambda _{\phi _-}$ and
$\phi -\phi (s_j)-\phi _{\pm,j}$ vanishes to the second order on $\pi
_x(K_\pm )$.  It is also clear that $\Lambda _\phi ,\Lambda _{\phi
  _\pm}$ intersect cleanly along $K_{\pm}$, so we get (cf
\cite[Section 11]{KFP2}): \ekv{2w.10.3} { \phi _+-(\phi -\phi
  (s_j))\asymp {\rm dist\,}(x,\pi _x(K_+))^2,\ \phi -\phi (s_j)-\phi
  _-\asymp {\rm dist\,}(x,\pi _x(K_-))^2.  }

We also know from Proposition \ref{ss2} and the isospectrality recalled prior to (\ref{ss.24}) 
that $\phi''(s_j){A}^{\mathrm{t}}$ (which is the linearization of
${{-\frac{1}{2}H_q}_|}_{\Lambda _{-\phi }}$, where we also notice that
$q$ vanishes on $\Lambda _{\pm \phi }$) has precisely one negative
eigenvalue, that the other eigenvalues have real parts $>0$, and that
$a_{j,0}(s_j)\ne 0$ is an eigenvector corresponding to the negative
eigenvalue. Here we recall that $\phi ''(s_j)A$ and $\phi
''(s_j){A}^{\mathrm{t}}$ are isospectral, since $(\phi
''(s_j)A)^{A,*}=\phi ''(s_j){A}^{\mathrm{t}}$. Also, if $a_{j,0}^*\ne 0$
is an eigenvector corresponding to the negative eigenvalue of $\phi
''(s_j)A$, then by (\cite[(11.27)]{KFP2}),
\ekv{ps.4} { (a_{j,0}^*|a_{j,0}(s_j))_A\ne 0.  }

Here we point out that the construction of the quasimode in
(\ref{qms.4}) starts by choosing $a_{j,0}(s_j)$ to be a non-vanishing
eigenvector of $\phi ''(s_j){A}^{\mathrm{t}}$, corresponding to the negative eigenvalue.

\section{Generalized PT symmetry and consequences}
\label{pt}
\setcounter{equation}{0}

\subsection{First remarks on the KFP case}\label{kfp}

In this subsection, as a preliminary and pedagogical step, we show the
reality of the small eigenvalues of the KFP operator
{(\ref{in.0})}, $n=2d$.
This corresponds to the general  supersymmetric case with
\ekv{kfp.3}
{
\phi (x,y)=\frac{y^2}{2}+V(x), \ \ \ \ A=\frac{1}{2}\left(\begin{array}{ccc} 0 &1\\ -1 &\gamma  \end{array}\right).
}
Here, $V$ is a smooth real-valued Morse function on ${\R}^d_x$ with $V''\in C_b^\infty $, $|V'(x)|\ge 1/C$ for $|x|\ge C$. Let 
${\bf m}_1, ..., {\bf m}_{n_0}$  be  the (non-degenerate) local minima of $V$, so that  $m_1 = ({\bf m}_1,0),..,m_{n_0} = ({\bf m}_{n_0}, 0)$ are the 
local minima of  $\phi$. In this case the spectral subspace $E^{(0)}$ is spanned by the functions $e_k^{(0)}$, $k=1,..,n_0$, where
\ekv{kfp.22}
{
e_k^{(0)}=h^{-\frac{n}{4}}e^{-\frac{1}{h}(\phi (x,y)-\phi ({\bf m}_j,0))}\chi_k (x, y)+{\cal O}(h^\infty )\hbox{ in }L^2.
}
where  $\chi_k \in C_0^\infty ({\R}^{n})$ was defined before.

\par Let $\kappa :{\R}^{n}\to {\R}^{n}$ be given by $\kappa (x,y)=(x,-y)$. Put $U_\kappa u=u\circ \kappa $, $u\in L^2({\R}^{n})$, so that 
$U_\kappa $ is unitary on $L^2({\R}^n)$ and also self-adjoint and equal to its own inverse. We introduce the Hermitian form
\ekv{kfp.4}
{
(u|v)_\kappa =(U_\kappa u|v),\ u,v\in L^2({\R}^{n}).
}

In the following our operators will be real so we can restrict the attention to real functions and differential forms. Notice that
\ekv{kfp.5}{P^*=U_\kappa ^{-1}PU_\kappa .}
Consequently,
$$
(Pu|v)_\kappa =(U_\kappa Pu|v)=(P^*U_\kappa u|v)=(U_\kappa u|Pv)=(u|Pv)_\kappa ,
$$
so $P$ is formally self-adjoint with respect to the Hermitian form $(u|v)_\kappa $.

\begin{prop}\label{kfp1}
The restriction of \  $(\cdot |\cdot)_\kappa $ to $E^{(0)}\times E^{(0)}$ is positive definite uniformly with respect to $h$, for $h$ small enough.
\end{prop}

\par\noindent {\it Proof.}
Since $\phi \circ \kappa =\phi $, we see that there exists $a_j>0$ independent of $h$, such that
$$
(e_k^{(0)}|e_k^{(0)})_\kappa =a_k+{\cal O}(h),\quad (e_k^{(0)}|e_{k'}^{(0)})_\kappa ={\cal O}(h^{\infty }) \hbox{ for } k\ne k'.
$$
Hence for $u=\sum_1^{n_0}u_k e_k^{(0)}$, we get
$$
(u|u)_\kappa =\sum_1^{n_0}(a_k+{\cal O}(h))|u_k|^2\ge \Vert u\Vert^2/C,
$$
and the proof is complete.\hfill{$\Box$}

\medskip
\par
In conclusion,
\begin{prop}
The restriction of  $P:E^{(0)}\to E^{(0)}$ is self-adjoint with respect to the inner product $(.|.)_\kappa $.
\end{prop}
In particular, the $n_0$ eigenvalues of $P$ with real part $<h/C$ are all real. Note that they also are ${\cal O}(h^\infty )$.

\begin{remark}\label{nonpar}
  The Maxwellian $e^{-(V(x)+y^2/2)/h}$ is an even eigenfunction with
  respect to $y$, associated to the eigenvalue 0. Let us show that no
  other eigenfunction in $E^{(0)}$ than multiples of the Mawxellian
  can be even.

In fact, assume the contrary, $U_\kappa u=u\ne 0$, $Pu=\mu u$, and
apply $U_\kappa $ to the eigenvalue equation. Using that $U_\kappa
P=P^*U_\kappa $, we get $P^*u=\mu u$. (Taking the differences of
the two equations, we get $(P-P^*)u=0$ so that the Hamilton field in
$P$ annihilates $u$). Taking the sum of the two equations, we get
$\frac{1}{2}(P+P^*)u=\mu u$, i.e.  $\frac{\gamma }{2}(y-h\partial
_y)\cdot (y+h\partial _y)u= \mu u$ (for almost all $x$). Since
$\mu $ is not among the eigenvalues of the harmonic oscillator
part, we conclude that $u=0$ and get a contradiction.
\end{remark}

\subsection{Extra symmetry and self-adjointness in the  general case}\label{sy}

Now we consider the general supersymmetric case and adopt the assumptions of Subsection \ref{ss}.
Let $\kappa :{\R}^n\to {\R}^n$ be linear and satisfy
\ekv{sy.1}{\kappa ^2=1.}
After a linear change of coordinates we may assume that $\kappa =\un_{\R^{n-d}}\oplus (-\un_{\R^d})$, for some $d\in\{ 1,2,...,n\}$. Also assume that
\ekv{sy.2}
{
\phi \circ \kappa =\phi .
}
Then $\kappa $ maps the critical set of
$\phi $ into itself. We assume that
\ekv{sy.2.5}
{
\kappa(x)=x \hbox{ for all critical points of }\phi \hbox{ of index 0
  or 1}.
}
On differential $k$-forms we define $U_{\kappa }=\kappa ^*$ as the pull-back in the usual sense. On $0$-forms, we get
$$
U_\kappa u=u\circ \kappa .
$$
On differential $1$-forms, we get
$
U_\kappa \omega ={\kappa}^{\mathrm{t}}(\omega \circ \kappa ),
$
where $\omega \circ  \kappa $ means the $1$-form obtained from $\omega $ by composing the coefficients with $\kappa $. More generally, for
$k$-forms, we get
$$
U_\kappa \omega =\wedge^k({\kappa}^{\mathrm{t}})(\omega \circ \kappa ).
$$
Using the fact that pullback and exterior differentiation commute, together with the invariance (\ref{sy.2}), we get
\ekv{sy.3}
{
U_\kappa ^{-1}d_\phi U_\kappa =d_\phi .
}

We next assume that
\ekv{sy.4}
{
\kappa A={A}^{\mathrm{t}}{\kappa}^{\mathrm{t}},
}
which can also be written as $\kappa {A}^{\mathrm{t}} = A {\kappa}^{\mathrm{t}}$, since $\kappa ^2=1$. Equivalently, $\kappa A$
(or $A{\kappa}^{\mathrm{t}}$) is symmetric. In the case of KFP, we have
\ekv{sy.5}
{
A=\frac{1}{2}\left(\begin{array}{ccc} 0 &1\\ -1 &\gamma  \end{array}\right),\ \phi =\frac{y^2 }{2}+V(x),\ \kappa =\left(\begin{array}{ccc}1 &0\\ 0 &-1 \end{array}\right)
}
and we see that
$$
\kappa A= \frac{1}{2}\left(\begin{array}{ccc}0 &1\\ 1 &-\gamma  \end{array}\right) .
$$
\begin{prop}\label{sy1} The bilinear form defined by
$$(u|v)_{A,\kappa }:=(U_\kappa u|v)_A$$
is a Hermitian form on the space of square integrable $k$-forms.
\end{prop}
\proof
We have
\begin{eqnarray*}
(u|v)_{A,\kappa }&=&\langle (\wedge^kA) (\wedge^k({\kappa}^{\mathrm{t}})) u\circ \kappa |v\rangle=\\
\langle (\wedge^k(A{\kappa}^{\mathrm{t}}))u|v\circ \kappa \rangle&=&\langle (\wedge^k(\kappa {A}^{\mathrm{t}}))u|v\circ \kappa \rangle=\\
\langle (\wedge ^k\kappa) (\wedge ^k {A}^{\mathrm{t}})u|v\circ \kappa \rangle
&=&\langle (\wedge ^k {A}^{\mathrm{t}})u|(\wedge^k {\kappa}^{\mathrm{t}})v\circ \kappa \rangle=\\
\langle u|(\wedge^kA) (\wedge^k{\kappa}^{\mathrm{t}})v\circ \kappa\rangle&=&
(v|u)_{A,\kappa },
\end{eqnarray*}
where complex conjugate signs are absent since we restrict the attention to real forms.\hfill{$\Box$}

\medskip
The same type of calculation shows that
$$
(u|U_\kappa v)_A=(U_\kappa u|v)_{A^{\mathrm{t}}}.
$$

\begin{prop}\label{sy2}
We have
$$
U_\kappa d_\phi ^{{A}^{\mathrm{t}},*}=d_\phi ^{A,*}U_\kappa .
$$
\end{prop}
\proof
Let $u$ be a $k$-form and $v$ be a $(k+1)$-form and consider
$$
(u|U_\kappa d_\phi ^{A^{\mathrm{t}},*}v)_A=(U_\kappa u|d_\phi ^{A^{\mathrm{t}},*}v)_{A^{\mathrm{t}}}=(d_\phi U_\kappa u|v)_{A^{\mathrm{t}}}.
$$
Similarly,
$$
(u|d_\phi ^{A,*}U_\kappa v)_A=(d_\phi u|U_\kappa v)_A=(U_\kappa d_\phi u|v)_{A^{\mathrm{t}}},
$$
and we get the proposition in view of (\ref{sy.3}).\hfill{$\Box$}

\medskip
\par It follows that
\ekv{sy.6}
{
U_\kappa \Delta_{{A}^{\mathrm{t}}}=\Delta _AU_\kappa ,
}
which implies that $\Delta_A$ is formally self-adjoint for the Hermitian form $(|)_{A,\kappa }$. In fact,
\begin{eqnarray*}
&&(\Delta_Au|v)_{A,\kappa }=(U_\kappa \Delta_Au|v)_A=(\Delta_{{A}^{\mathrm{t}}}U_\kappa u|v)_A=\\
&&(U_\kappa u|(\Delta_{{A}^{\mathrm{t}}})^{A,*}v)_A=(U_\kappa u|\Delta_Av)_A=(u|\Delta_Av)_{A,\kappa },
\end{eqnarray*}
where we used (\ref{su.6}). 

We also have
\begin{eqnarray*}
(d_\phi u|v)_{A,\kappa }=(U_\kappa d_\phi u|v)_A=(d_\phi U_\kappa u|v)_A=(U_\kappa u|d_\phi ^{A,*}v)_A=(u|d_\phi ^{A,*}v)_{A,\kappa },
\end{eqnarray*}
which shows that $d_\phi ^{A,*}$ is the adjoint of $d_\phi $ for our Hermitian product $(.|.)_{A,\kappa }$.

\subsection{Positivity questions}\label{ps}

First we shall prove that the eigenvalues are real. We know from (\ref{ss.14}) that
$$
\Re\left(-\Delta_A^{(0)}\right)\ge 0,
$$
so all eigenvalues have real part $\ge 0$. Proposition \ref{kfp1} and its proof carry over without any changes, so we
know that $-\Delta_A^{(0)}:E^{(0)}\to E^{(0)}$ is self-adjoint with respect to the {\it inner product} $(.|.)_{\kappa }=(|)_{A,\kappa }$ on
$E^{(0)}$. In particular, the eigenvalues of $-\Delta_A^{(0)}$ with real part in $[0,h/C]$ are all real and $={\cal O}(h^\infty )$.

\par We next consider $-\Delta_A^{(1)}$ and a vector in $E^{(1)}$ of
the form, 
\ekv{ps.2bis}
{
e_j^{(1)} = (h^{-n/4} a_j(x;h)+r_j)e^{-\phi _{+,j}(x)/h}\theta_j(x)+\widetilde{r}_j e^{-S_0/h},
}
where $S_0$ is a positive constant, $r_j$, $B^{1/2}hDr_j$ are ${\cal O}(h^\infty )$ in $L^2$ and similarly for $\widetilde{r}_j$, 
$\theta_j \in C^{\infty}_0({\rm neigh}(s_j))$, $\theta_j=1$ near $s_j$, 
and $ a_j(x;h)\sim a_{j,0}(x)+ha_{j,1}(x)+....$ where  $a_{j,0}(s_j)\ne 0$ is  an eigenvector corresponding to the negative
eigenvalue of $\phi ''(s_j)A^{\mathrm{t}}$. The vector $e_j^{(1)}$ is obtained by the spectral projection of the truncation of a quasimode as
in (\ref{qms.4}). As remarked at the end of Section \ref{re}, when constructing a quasimode for
$-\Delta ^{(1)}_{A^{\mathrm{t}}}$, we consider $a_{j,0}^*\neq 0$, which is an eigenvector corresponding to the negative eigenvalue of $\phi ''(s_j)A$.
Then
\ekv{ps.4bis}
{
(a_{j,0}^*|a_{j,0}(s_j))_A\ne 0.
}
Notice that we can take $a_{j,0}^*=\kappa ^{\mathrm{t}}a_{j,0}(s_j)$.

\par Now fix $j$ for a while, suppress this subscript from the
notation and write $a_0=a_{j,0}(s_j)$. We shall study the sign of
\ekv{ps.5} { (\kappa ^{\mathrm{t}}a_0|a_0)_A, } which we already know
is real and non-vanishing. Recall the relations \ekv{ps.6} { \kappa
  ^{\mathrm{t}}\phi '' \kappa =\phi '', \quad \phi ''=\phi ''(s_j), }
\ekv{ps.7} { \kappa A=A^{\mathrm{t}}\kappa ^{\mathrm{t}}, } and that
$A=B+C$ with $B=B^{\mathrm{t}}\ge 0$, $C=-C^{\mathrm{t}}$, $A$
bijective. Write 
\ekv{ps.8} { \kappa =\left(\begin{array}{ccc}1 &0\\ 0
      &-1 \end{array}\right):E_+\oplus E_-\to E_+\oplus E_- } where
$E_+\simeq{\mathbb{R}}^{n-d}$ is the kernel of $\kappa - \mathrm{Id}$
and $E_-\simeq{\mathbb{R}}^d$ is the kernel of $\kappa +
\mathrm{Id}$. Let $(\mathbb{R}^n)^*=E_+^*\oplus E_-^*$ be the dual
decomposition, so that $E_+^*\simeq({\mathbb{R}}^{n-d})^*$ is the kernel of $\kappa^{\mathrm{t}} - \mathrm{Id}$
and $E_-^*\simeq({\mathbb{R}}^d)^{*}$ is the kernel of $\kappa^{\mathrm{t}} + \mathrm{Id}$.

Then
(\ref{ps.6}), (\ref{ps.7}) say that \ekv{ps.9} { \phi
  ''=\left(\begin{array}{ccc}\Phi _1 &0\\ 0&\Phi _2 \end{array}\right) :
  \ E_+\oplus
  E_-\to E_+^*\oplus E_-^*, } \ekv{ps.10} {
  A=\left(\begin{array}{ccc}B_{11} &C_{12}\\
      -C_{12}^{\mathrm{t}} &B_{22} \end{array}\right) :\ E_+^*\oplus
  E_-^*\to E_+\oplus E_-, } where
$B_{11}, B_{22}\ge 0$. We now make a continuous deformation of $A$
into the identity in such a way that the above properties of $A$ are
preserved, deform $a_0$ correspondingly so that it remains an
eigenvector corresponding to the negative eigenvalue of $\phi
''A^{\mathrm{t}}$. (Notice here that the dynamical assumption (\ref{ny.17})
is automatically satisfied once we make $B$ positive
definite.) Then the quantity (\ref{ps.5}) remains non-vanishing and
hence of constant sign, so we have reduced the problem to that of
studying the sign of $(\kappa ^{\mathrm{t}}a_0|a_0)$ when $a_0$ is a
non-vanishing eigenvector of $\phi ''$ corresponding to its negative
eigenvalue. This gives

\begin{prop}\label{ps1} There are only two possibilities:
\begin{enumerate}[i)]
\item If $\Phi _1$ is the component in {\rm (\ref{ps.9})} that has a negative eigenvalue, then the quantity {\rm (\ref{ps.5})} is $>0$.
\item
If $\Phi _2$ is the component in {\rm (\ref{ps.9})} that has a negative eigenvalue, then the quantity {\rm (\ref{ps.5})} is $<0$. \end{enumerate}
\end{prop}
A more invariant way of describing the two cases in the proposition, is to say that we are in the first case if $\phi ''$ is
positive as a quadratic form on $E_-$ and that we are in the second case when $\phi ''$ is positive on $E_+$. In the case of KFP, we have
\ekv{ps.11}
{
\phi ''=\left(\begin{array}{ccc}V'' &0\\0 &1 \end{array}\right),\quad \kappa =\left(\begin{array}{ccc}1 &0\\ 0&-1 \end{array}\right),
}
so we are in the first case and the quantity (\ref{ps.5}) is positive.

\par Recall the final observation in Subsection \ref{sy}. Combining it with Proposition \ref{ps1} and the method of stationary phase as in the case
of $0$-forms, we get
\begin{prop} \label{ps1bis}
Assume  that at every saddle point, we are in the case {\rm 1} of Proposition {\rm \ref{ps1}}. Then the restriction of
$(\cdot |\cdot )_{A,\kappa }$ to $E^{(1)}\times E^{(1)}$ is uniformly positive definite and if we equip $E^{(0)}$, $E^{(1)}$ with the scalar products
$(\cdot |\cdot )_\kappa $ and $(\cdot |\cdot )_{A,\kappa }$, then the adjoint of $d_\phi :E^{(0)}\to E^{(1)}$ is $d_\phi ^{A,*}$.
\end{prop}
{
\begin{remark}\label{ps2}
Let us finally mention that our discussion applies 
to the case of the usual Witten complex. In that case, we let
$A=\imath:(\mathbb{R}^n)^*\to \mathbb{R}^n$ be the application whose
matrix is the identity in Euclidean coordinates. Naturally much of our
machinery is redundant in that case, but it may be of interest to point
out that our Theorems \ref{gl4}, \ref{fa1} seem to be new in that
degree of generality for the standard Witten Laplacian.  
\end{remark}}
\section{Labelling and quasimodes for multiple wells}
\label{lbl}
\setcounter{equation}{0}

\subsection{Separating saddle points and critical components}
\label{SSP_and_CC}
Let $\phi \in C^\infty (\mathbb{R}^n;\mathbb{R})$ be a Morse function
satisfying (\ref{eq3.01}) and (\ref{morse}). Also assume that
\ekv{eq1}
{\phi (x)\to +\infty ,\ x\to \infty ,}
so that in view of (\ref{morse}), we have $\phi (x)\ge |x|/C$ for
$|x|\ge C$, where $C>0$. In particular, $e^{-\phi /h}\in L^2$. {In
\cite{KFP2} we also treated the case when (\ref{eq1}) does not hold,
and we hope to extend our present result to that case in the future.}

\par
The function $\phi $ has finitely many critical points. The critical points of index 1 will be called saddle points. In what follows, we shall
be concerned with the saddle points and the local minima of $\phi$.

\par
Let $s$ be a saddle point of $\phi$ and let $B(s,r)=\{ x\in
\mathbb{R}^n;\, |x-s|<r\}$. Then for $r>0$ small enough, the set
$$
\{ x\in B(s,r);\, \phi (x)<\phi (s)\}
$$
has precisely 2 connected components, $C_j(s,r)$, $j=1,2$, with $C_j(s,\widetilde{r})\subset C_j(s,r)$, if $0<\widetilde{r}\le r$.
\begin{dref}\label{def1}
We say that $s\in \mathbb{R}^n$ is a separating saddle point (ssp) if it is a saddle
point and $C_1(s,r)$ and $C_2(s,r)$ are contained in different connected components of the set
$\{ x\in \mathbb{R}^n;\, \phi (x)<\phi (s)\}$. Let $\mathrm{SSP}$ denote the set of ssps.
\end{dref}

\par
Notice that this definition depends on the global behavior of $\phi
$. It can be localized somewhat:
\begin{prop}\label{prop1}
Let $s\in \mathbb{R}^n$ be a saddle point and let $\sigma \in ]\phi
(s),+\infty[$. Then $s$ is a ssp if and only if $C_1(s,r)$ and $C_2(s,r)$ are
contained in different components of $\{ x\in F(s,\sigma );\, \phi
(x)<\phi (s)\}$, where $F(s,\sigma )$ denotes the connected component
of $s$ in $\phi ^{-1}(]-\infty ,\sigma [)$.
\end{prop}

\begin{dref}\label{def2}
A connected component $E$ of the sublevel set $\phi ^{-1}(]-\infty ,\sigma [)$ will be
called a critical component (cc) if $\partial E\cap \mathrm{SSP}\ne
\emptyset $ or if $E=\mathbb{R}^n$.
\end{dref}

\par
We shall now describe a labelling system for the set LM of local minima of the function $\phi$ and a natural injective map from LM to the set CC
of critical components, instrumental in constructing an appropriate system of quasimodes, adapted to the local minima. When doing so, it will turn out
to be convenient to label the elements of both sets by 2-tuples of positive integers, $\mathbb{N}^2$. Let therefore $m_{1,1}$ stand
for a point of global minimum of $\phi$, arbitrarily chosen, but kept fixed in the following discussion. Associated to $m_{1,1}$ we have the
critical component $E_{1,1} = \mathbb{R}^n$ and we let the associated saddle point value be $\sigma_1  = \sigma(E_{1,1}) = +\infty$. Let next
$$
\sigma_2  = \sup \phi\left(\mathrm{SSP}\right).
$$
Then the sublevel set ${\cal L}(\sigma_2) := \{x\in \mathbb{R}^n, \phi(x)<\sigma_2\}$ is the union of its finitely many connected components
(all critical), of which there is precisely one containing the point $m_{1,1}$. The remaining connected components of the sublevel set
${\cal L}(\sigma_2)$ will be labelled as $E_{2,k}$, $1\leq k\leq N_2$, $N_2\geq 1$. Associated to each $E_{2,k}$, we let $m_{2,k}$ be a point of
global minimum of the restriction of $\phi$ to $E_{2,k}$, $1\leq k\leq N_2$.

\par
Continuing the labelling procedure, we let $\sigma_3$ be the largest number of the form $\phi(s)$, $s\in \mathrm{SSP}$, such that $\sigma_3 < \sigma_2$.
Decompose the sublevel set ${\cal L}(\sigma_3)$ into its connected components and perform the labelling as follows: We omit all
those components that contain the already labelled minima $m_{1,1}$ and $m_{2,k}$, $1\leq k \leq N_2$. Some of these components may be non-critical.
The remaining ones are critical and we label them as $E_{3,j}$,
$1\leq j\leq N_3$, $N_3\geq 1$. We then let $m_{3,j}$ be a point of global minimum of the restriction of $\phi$ to $E_{3,j}$, $1\leq j \leq N_3$.

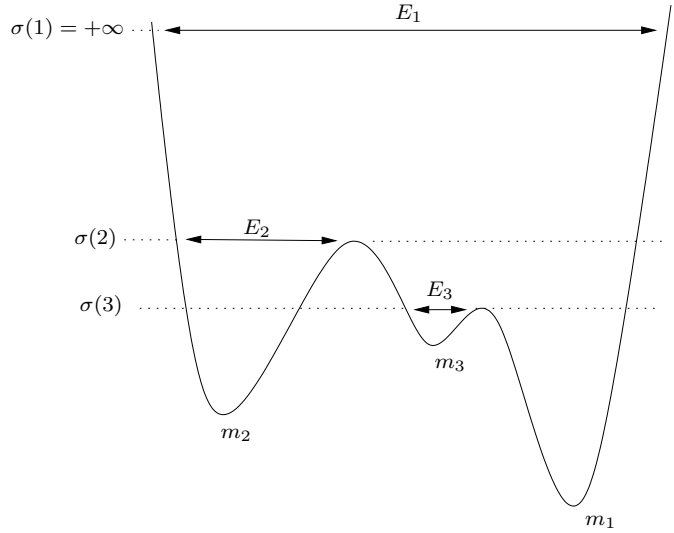
\begin{figure}[htbp]
\begin{center}

\begin{picture}(0,0)%
\includegraphics{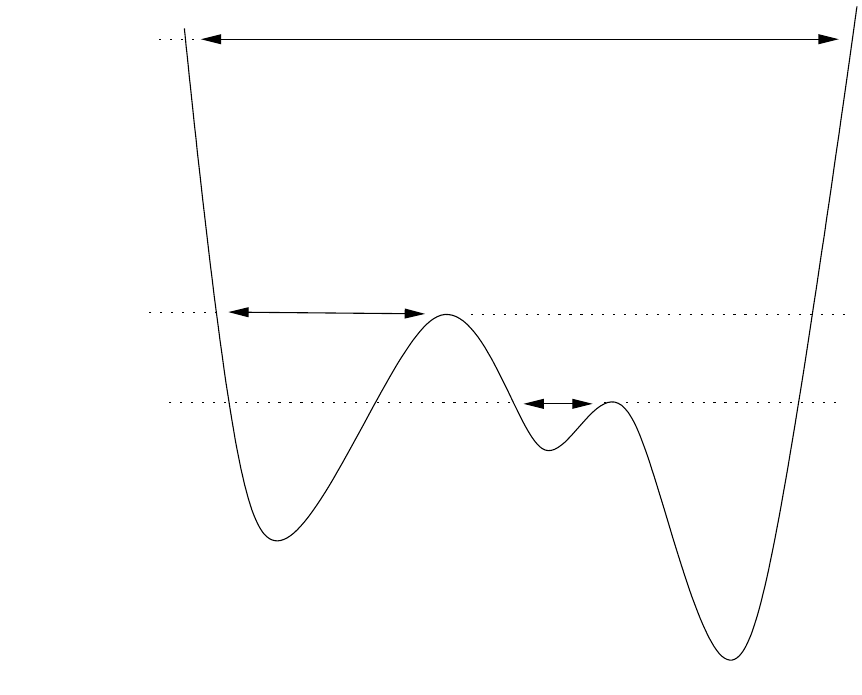}%
\end{picture}%
%
%
\setlength{\unitlength}{2763sp}%
\begingroup\makeatletter\ifx\SetFigFont\undefined%
\gdef\SetFigFont#1#2#3#4#5{%
  \reset@font\fontsize{#1}{#2pt}%
  \fontfamily{#3}\fontseries{#4}\fontshape{#5}%
  \selectfont}%
\fi\endgroup%
\begin{picture}(5889,4750)(386,-6904)
\put(5501,-6836){\makebox(0,0)[lb]{\smash{{\SetFigFont{8}{9.6}{\rmdefault}{\mddefault}{\updefault}{\color[rgb]{0,0,0}$m_1$}%
}}}}
\put(2264,-6061){\makebox(0,0)[lb]{\smash{{\SetFigFont{8}{9.6}{\rmdefault}{\mddefault}{\updefault}{\color[rgb]{0,0,0}$m_2$}%
}}}}
\put(4164,-5436){\makebox(0,0)[lb]{\smash{{\SetFigFont{8}{9.6}{\rmdefault}{\mddefault}{\updefault}{\color[rgb]{0,0,0}$m_3$}%
}}}}
\put(401,-2449){\makebox(0,0)[lb]{\smash{{\SetFigFont{8}{9.6}{\rmdefault}{\mddefault}{\updefault}{\color[rgb]{0,0,0}$\sigma(1) = +\infty$}%
}}}}
\put(3814,-2311){\makebox(0,0)[lb]{\smash{{\SetFigFont{8}{9.6}{\rmdefault}{\mddefault}{\updefault}{\color[rgb]{0,0,0}$E_1$}%
}}}}
\put(2464,-4236){\makebox(0,0)[lb]{\smash{{\SetFigFont{8}{9.6}{\rmdefault}{\mddefault}{\updefault}{\color[rgb]{0,0,0}$E_2$}%
}}}}
\put(4089,-4786){\makebox(0,0)[lb]{\smash{{\SetFigFont{8}{9.6}{\rmdefault}{\mddefault}{\updefault}{\color[rgb]{0,0,0}$E_3$}%
}}}}
\put(964,-4324){\makebox(0,0)[lb]{\smash{{\SetFigFont{8}{9.6}{\rmdefault}{\mddefault}{\updefault}{\color[rgb]{0,0,0}$\sigma(2)$}%
}}}}
\put(1014,-4924){\makebox(0,0)[lb]{\smash{{\SetFigFont{8}{9.6}{\rmdefault}{\mddefault}{\updefault}{\color[rgb]{0,0,0}$\sigma(3)$}%
}}}}
\end{picture}%

\caption{Labelling in the generic case}
\label{generic}
\end{center}
\end{figure}

\par
We go on with this procedure, proceeding in the order dictated by the elements of the set $\phi(\mathrm{SSP})$, arranged in the decreasing
order, until all local minima have been enumerated. In this way, we obtain the labelling of the set LM of local minima of $\phi$ by indices of
the form $k = (k^{\sigma}, k^{cc})\in \mathbb{N}^2$, where $k^{\sigma}$ is
the index corresponding to the separating saddle point values $\sigma$, arranged in the decreasing order,
$\sigma_1 > \sigma_2 > \sigma_3 > \ldots > \sigma_{N_0}$. It is then also clear that we get an injection from LM to CC, associating to each local
minimum $m_k$, the critical component $E_k$, containing $m_k$. We equip the set of indices $k$ with the lexicographical order.

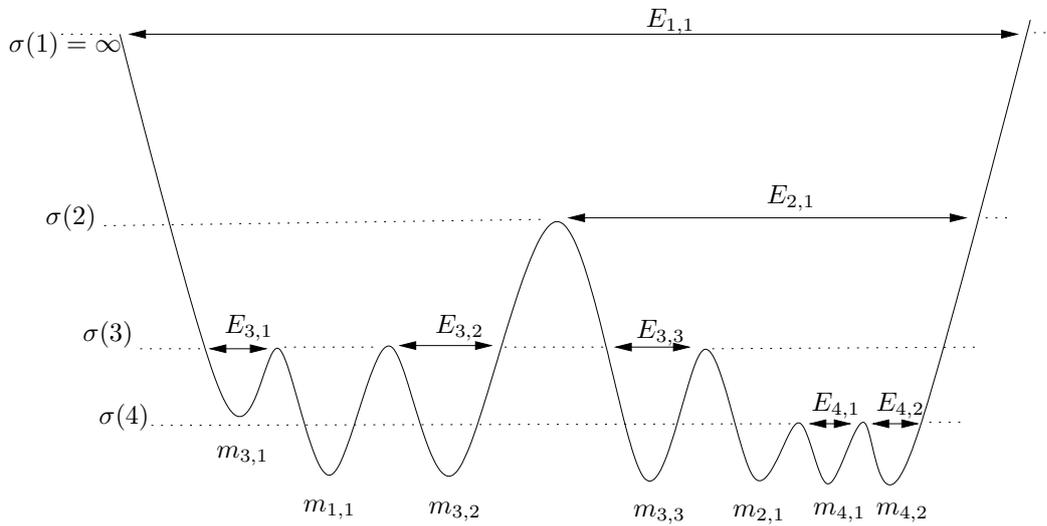
\begin{figure}[htbp]
\begin{center}

\begin{picture}(0,0)%
\includegraphics{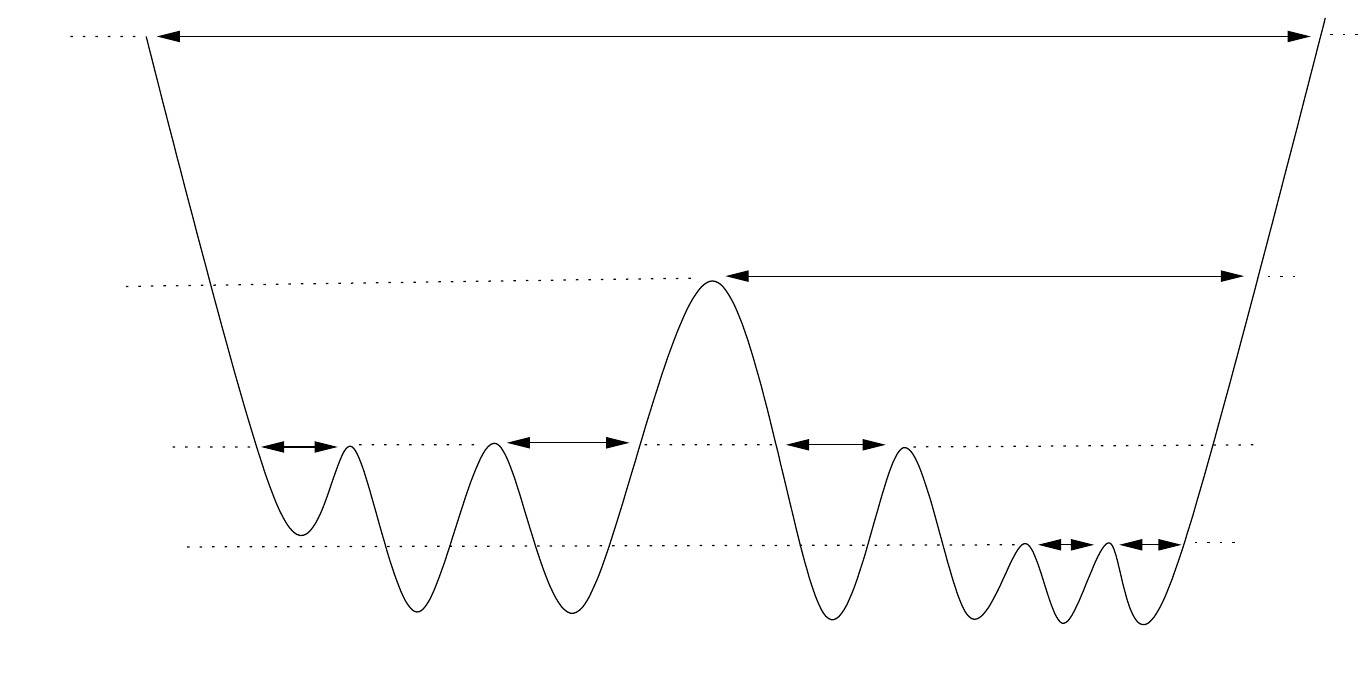}%
\end{picture}%
%
%
\setlength{\unitlength}{3158sp}%
\begingroup\makeatletter\ifx\SetFigFont\undefined%
\gdef\SetFigFont#1#2#3#4#5{%
  \reset@font\fontsize{#1}{#2pt}%
  \fontfamily{#3}\fontseries{#4}\fontshape{#5}%
  \selectfont}%
\fi\endgroup%
\begin{picture}(8202,4062)(1049,-5429)
\put(6951,-2911){\makebox(0,0)[lb]{\smash{{\SetFigFont{10}{12.0}{\rmdefault}{\mddefault}{\updefault}{\color[rgb]{0,0,0}$E_{2,1}$}%
}}}}
\put(2739,-3924){\makebox(0,0)[lb]{\smash{{\SetFigFont{10}{12.0}{\rmdefault}{\mddefault}{\updefault}{\color[rgb]{0,0,0}$E_{3,1}$}%
}}}}
\put(4376,-3924){\makebox(0,0)[lb]{\smash{{\SetFigFont{10}{12.0}{\rmdefault}{\mddefault}{\updefault}{\color[rgb]{0,0,0}$E_{3,2}$}%
}}}}
\put(5939,-3961){\makebox(0,0)[lb]{\smash{{\SetFigFont{10}{12.0}{\rmdefault}{\mddefault}{\updefault}{\color[rgb]{0,0,0}$E_{3,3}$}%
}}}}
\put(7301,-4536){\makebox(0,0)[lb]{\smash{{\SetFigFont{10}{12.0}{\rmdefault}{\mddefault}{\updefault}{\color[rgb]{0,0,0}$E_{4,1}$}%
}}}}
\put(7801,-4536){\makebox(0,0)[lb]{\smash{{\SetFigFont{10}{12.0}{\rmdefault}{\mddefault}{\updefault}{\color[rgb]{0,0,0}$E_{4,2}$}%
}}}}
\put(1064,-1724){\makebox(0,0)[lb]{\smash{{\SetFigFont{10}{12.0}{\rmdefault}{\mddefault}{\updefault}{\color[rgb]{0,0,0}$\sigma(1) = \infty$}%
}}}}
\put(1351,-3074){\makebox(0,0)[lb]{\smash{{\SetFigFont{10}{12.0}{\rmdefault}{\mddefault}{\updefault}{\color[rgb]{0,0,0}$\sigma(2)$}%
}}}}
\put(1639,-3986){\makebox(0,0)[lb]{\smash{{\SetFigFont{10}{12.0}{\rmdefault}{\mddefault}{\updefault}{\color[rgb]{0,0,0}$\sigma(3)$}%
}}}}
\put(1764,-4636){\makebox(0,0)[lb]{\smash{{\SetFigFont{10}{12.0}{\rmdefault}{\mddefault}{\updefault}{\color[rgb]{0,0,0}$\sigma(4)$}%
}}}}
\put(2676,-4899){\makebox(0,0)[lb]{\smash{{\SetFigFont{10}{12.0}{\rmdefault}{\mddefault}{\updefault}{\color[rgb]{0,0,0}$m_{3,1}$}%
}}}}
\put(3351,-5324){\makebox(0,0)[lb]{\smash{{\SetFigFont{10}{12.0}{\rmdefault}{\mddefault}{\updefault}{\color[rgb]{0,0,0}$m_{1,1}$}%
}}}}
\put(4326,-5336){\makebox(0,0)[lb]{\smash{{\SetFigFont{10}{12.0}{\rmdefault}{\mddefault}{\updefault}{\color[rgb]{0,0,0}$m_{3,2}$}%
}}}}
\put(5914,-5361){\makebox(0,0)[lb]{\smash{{\SetFigFont{10}{12.0}{\rmdefault}{\mddefault}{\updefault}{\color[rgb]{0,0,0}$m_{3,3}$}%
}}}}
\put(6689,-5361){\makebox(0,0)[lb]{\smash{{\SetFigFont{10}{12.0}{\rmdefault}{\mddefault}{\updefault}{\color[rgb]{0,0,0}$m_{2,1}$}%
}}}}
\put(7314,-5349){\makebox(0,0)[lb]{\smash{{\SetFigFont{10}{12.0}{\rmdefault}{\mddefault}{\updefault}{\color[rgb]{0,0,0}$m_{4,1}$}%
}}}}
\put(7801,-5349){\makebox(0,0)[lb]{\smash{{\SetFigFont{10}{12.0}{\rmdefault}{\mddefault}{\updefault}{\color[rgb]{0,0,0}$m_{4,2}$}%
}}}}
\put(6014,-1524){\makebox(0,0)[lb]{\smash{{\SetFigFont{10}{12.0}{\rmdefault}{\mddefault}{\updefault}{\color[rgb]{0,0,0}$E_{1,1}$}%
}}}}
\end{picture}%

\caption{Labelling in the general case}
\label{general}
\end{center}
\end{figure}

\par
Associated to the labelling procedure above is the connected graph, having a structure of a tree, or rather a root, obtained by letting $\sigma$
vary from $+\infty$ to $-\infty$, and representing the various components of ${\cal L}(\sigma)$ as points that move vertically and split at each
separating saddle point level.

\begin{figure}[htbp]
\begin{center}

\begin{picture}(0,0)%
\includegraphics{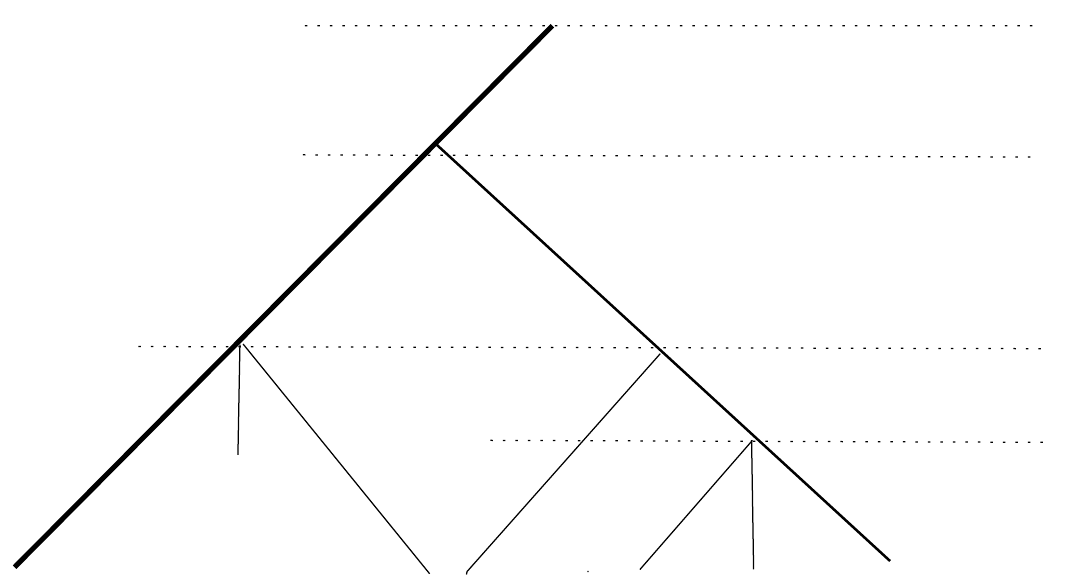}
\end{picture}%
%
%
\setlength{\unitlength}{3158sp}%
\begingroup\makeatletter\ifx\SetFigFont\undefined%
\gdef\SetFigFont#1#2#3#4#5{%
  \reset@font\fontsize{#1}{#2pt}%
  \fontfamily{#3}\fontseries{#4}\fontshape{#5}%
  \selectfont}%
\fi\endgroup%
\begin{picture}(6518,3449)(1636,-3166)
\put(4176,-199){\makebox(0,0)[lb]{\smash{{\SetFigFont{10}{12.0}{\rmdefault}{\mddefault}{\updefault}{\color[rgb]{0,0,0}$E_{1,1}$}%
}}}}
\put(4351,-2761){\makebox(0,0)[lb]{\smash{{\SetFigFont{10}{12.0}{\rmdefault}{\mddefault}{\updefault}{\color[rgb]{0,0,0}$E_{3,3}$}%
}}}}
\put(6064,-2199){\makebox(0,0)[lb]{\smash{{\SetFigFont{10}{12.0}{\rmdefault}{\mddefault}{\updefault}{\color[rgb]{0,0,0}$E_{2,1}$}%
}}}}
\put(6964,-3036){\makebox(0,0)[lb]{\smash{{\SetFigFont{10}{12.0}{\rmdefault}{\mddefault}{\updefault}{\color[rgb]{0,0,0}$E_{2,1}$}%
}}}}
\put(5164,-3049){\makebox(0,0)[lb]{\smash{{\SetFigFont{10}{12.0}{\rmdefault}{\mddefault}{\updefault}{\color[rgb]{0,0,0}$E_{4,1}$}%
}}}}
\put(6201,-3036){\makebox(0,0)[lb]{\smash{{\SetFigFont{10}{12.0}{\rmdefault}{\mddefault}{\updefault}{\color[rgb]{0,0,0}$E_{4,2}$}%
}}}}
\put(5014,-1224){\makebox(0,0)[lb]{\smash{{\SetFigFont{10}{12.0}{\rmdefault}{\mddefault}{\updefault}{\color[rgb]{0,0,0}$E_{2,1}$}%
}}}}
\put(3139,-1224){\makebox(0,0)[lb]{\smash{{\SetFigFont{10}{12.0}{\rmdefault}{\mddefault}{\updefault}{\color[rgb]{0,0,0}$E_{1,1}$}%
}}}}
\put(1651,-2724){\makebox(0,0)[lb]{\smash{{\SetFigFont{10}{12.0}{\rmdefault}{\mddefault}{\updefault}{\color[rgb]{0,0,0}$E_{1,1}$}%
}}}}
\put(2676,-2486){\makebox(0,0)[lb]{\smash{{\SetFigFont{10}{12.0}{\rmdefault}{\mddefault}{\updefault}{\color[rgb]{0,0,0}$E_{3,1}$}%
}}}}
\put(3426,-2736){\makebox(0,0)[lb]{\smash{{\SetFigFont{10}{12.0}{\rmdefault}{\mddefault}{\updefault}{\color[rgb]{0,0,0}$E_{3,2}$}%
}}}}
\put(8126,-649){\makebox(0,0)[lb]{\smash{{\SetFigFont{10}{12.0}{\rmdefault}{\mddefault}{\updefault}{\color[rgb]{0,0,0}$\sigma(2)$}%
}}}}
\put(8139,-1811){\makebox(0,0)[lb]{\smash{{\SetFigFont{10}{12.0}{\rmdefault}{\mddefault}{\updefault}{\color[rgb]{0,0,0}$\sigma(3)$}%
}}}}
\put(8126,-2349){\makebox(0,0)[lb]{\smash{{\SetFigFont{10}{12.0}{\rmdefault}{\mddefault}{\updefault}{\color[rgb]{0,0,0}$\sigma(4)$}%
}}}}
\put(8114,126){\makebox(0,0)[lb]{\smash{{\SetFigFont{10}{12.0}{\rmdefault}{\mddefault}{\updefault}{\color[rgb]{0,0,0}$\sigma(1) = \infty$}%
}}}}
\end{picture}%

\caption{Tree corresponding to the general case}
\label{generaltree}
\end{center}
\end{figure}
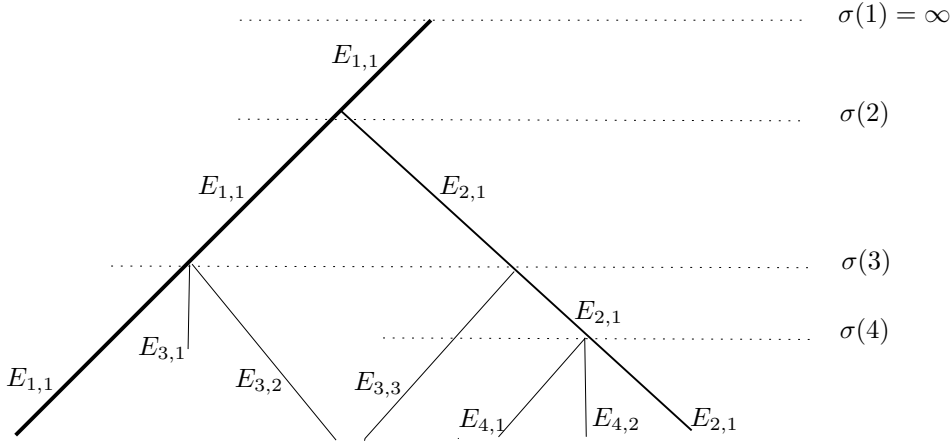

\par
We also notice that our labelling procedure has the property that if $E_{k'}$ and $E_k$ are two critical components with $E_{k'}\subset E_k$, with a proper
inclusion, then $k'>k$, and $E_k$ can be reached from $E_{k'}$ by following ascending links in the root. Given a critical component $E_k$, we let
$\sigma (k)= \sigma (E_k)$ be the corresponding saddle point value. Here, as above, we adopt the convention that
$E_{1,1} = \mathbb{R}^n$ is a critical component, and that $\sigma (1,1)=+\infty $.

\subsection{Cut-off functions and quasimodes}
\label{QM}
In this section, we shall build quasimodes adapted to the local minima of $\phi$ and the labelling of the set LM, described in the previous subsection.
When doing so, we shall first construct suitable cut-off functions, in terms of the corresponding critical components.

\par
Let $E=E(\sigma)\neq \mathbb{R}^n$ be a critical component and let $\sigma$ be
the corresponding saddle point value. Let $\chi _0\in C_0^\infty
(B(0,1),[0,1])$ be equal to 1 on $B(0,\frac{1}{2})$. For $\delta >0$
small enough, put
\ekv{eq2}
{
\widetilde{\phi }(x)=\phi (x)+\sum_{s\in \mathrm{SSP},\atop \phi
  (s)\ge \sigma}\delta ^2\chi _0\left(\frac{x-s}{\delta }\right).
}
For ${0}\le t\le \frac{\delta ^2}{C}$, with $C$ large
enough independent of $\delta $, there is a unique connected component
of $\widetilde{\phi }^{-1}(]-\infty ,\sigma + t[)$ which is $C\delta
$-close to $E(\sigma)$ in the sense that any of the two sets
is contained in the algebraic sum of the other and $B(0,C\delta)$. Let $\widetilde{E}_0$ be this component when
$t=\frac{\delta^2}{C}$. Let $f\in C^\infty (\mathbb{R};[0,1])$ be equal to 1 on
$]-\infty ,\sigma + \frac{\delta ^2}{2C}]$ and have its support in $]-\infty,\sigma+ \frac{\delta ^2}{C}[$.
The cut-off function associated to the critical component $E=E(\sigma)$ is defined as follows,
\ekv{eq3}
{
\chi _E(x)=1_{\widetilde{E}_0}(x)f(\widetilde{\phi }(x)).
}
Here $1_{\widetilde{E}_0}$ is the characteristic function of the set $\widetilde{E}_0$.

\par
Notice that outside the union of the balls $B(s,\delta )$, for $s\in \mathrm{SSP}$, $\phi (s)\ge \sigma$, we have
$\phi(x)\ge\sigma+\frac{\delta ^2}{2C}$, in the cut-off region
$\mathrm{supp\,}\nabla \chi _E$. We also notice that if $B(s,\delta )\cap
\mathrm{supp\,}\nabla \chi _E\ne \emptyset $ for some $s\in
\mathrm{SSP}$, then $\mathrm{dist\,}(s,E)\le {\cal O}(\delta )$. In
particular, $s$ can be a boundary point only of critical components
$\widetilde{E}$ which contain $E$.

\par
We shall now describe the construction of quasimodes associated to the local minima $m_k$, $k\in \mathbb{N}^2$, of $\phi$. Here we assume that the
minima have been labelled as described in the previous subsection. Let $E_k$ be the corresponding critical component containing $m_k$, as
described there.

\par
When $m_{1,1}$ is a point of global minimum of $\phi$, let $E_1 = \mathbb{R}^n$, and set
$$
f_{1,1}^{(0)}(x;h)=h^{-n/4}e^{-(\phi (x)-\phi (m_{11}))/h}.
$$
When $k\neq (1,1)$, we set, according to our labelling,
$$
f_k^{(0)}(x;h)=f_{m_k}^{(0)}(x;h)= h^{-\frac{n}{4}}\chi _{E_{k}}(x)e^{-(\phi(x)-\phi (m_k))/h},
$$
for each of the local minima $m_k$, $k\in \mathbb{N}^2$. Here $\chi_{E_k}$ has been defined as in (\ref{eq3}). Notice that then the quasimode $f_{k_2}^{(0)}$ is exponentially small near $m_{k_1}$, as soon
as $k_2>k_1$.

\par
From the properties of the cut-off functions $\chi_{E_k}$, we infer that our system of quasi-modes has the following two important properties:

\smallskip \par\noindent
1) We have $-\Delta _A^{(0)}(f_{11}^{(0)})=0$,
and for the other indices $k$, we obtain that
$$
-\Delta _A^{(0)}(f_k^{(0)})=[-\Delta _A^{(0)},\chi _{E_k}](h^{-\frac{n}{4}}e^{-(\phi -\phi (m_k))/h}).
$$
Let $\sigma(k)$ be the saddle point value associated to the critical component $E_k$ and let $S_k=\sigma(k)-\phi (m_k)$. Then we get
\begin{equation}\label{eq4}
\begin{split}
-\Delta _A^{(0)}(f_k^{(0)})=\sum_{E_k\subset \widetilde{E}\in
  \mathrm{CC}}&\sum_{s\in \mathrm{SSP}\cap \partial \widetilde{E}}1_{B(s,\delta
)}[-\Delta
_A^{(0)},\chi _{E_k}]h^{-\frac{n}{4}}e^{-(\phi -\phi (m_k))/h}\\
&+{\cal O}(h^{-N_0})e^{-(S_k+\frac{\delta ^2}{2C})/h},
\end{split}
\end{equation}
for some $N_0>0$.

\smallskip \par\noindent 2) We next claim that the quasimodes
$f_k^{(0)}$ enjoy the property of linear independence in $L^2$,
uniformly with respect to $h$. Indeed, notice first that
$\Vert{f_k^{(0)}}\Vert\asymp 1$. Consider a linear combination
$$
u=\sum u_k f_k^{(0)},\ \ u_k\in \mathbb{C},
$$
where it is understood that the summation extends over all indices
$k\in \mathbb{N}^2$, occurring in the labelling procedure of
Subsection \ref{SSP_and_CC}. For simplicity we now label them in lexicographical
order, $f_{k_1}, f_{k_2}$, ... $f_{k_{n_0}}$, where $n_0=\# \mathrm{LM}$, so
that $k_1=(1,1)$. As already noticed, if $k\neq (1,1)$, then
$f_k^{(0)}$ is exponentially small near $m_{1,1}=m_{k_1}$, while $f_{k_1}^{(0)}$
has a substantial part of its $L^2$ norm concentrated there: for every
$r>0$:
$$
\int_{B(m_{k_1},r)}|f_{k_1}^{(0)}|^2dx\ge \frac{1}{C},\quad C>0.
$$
Thus, if we multiply $u$ by $1_{B(m_{k_1},r)}\overline{f_{k_1}^{(0)}}$ and
integrate (assuming also that $r>0$ is sufficiently small), we get
$$
|(u|1_{B(m_{k_1},r)}f_{k_1}^{(0)})|\ge \frac{|u_{k_1}|}{C}-e^{-\frac{1}{Ch}}(\sum |u_k|^2)^{\frac{1}{2}}.
$$
Hence
$$
|u_{k_1}|\le C(\|u\|+e^{-\frac{1}{Ch}}\| u_\bullet\|),\ u_{\bullet
}=(u_k)\in \mathbb{R}^{\# \mathrm{LM}}.
$$
Next consider the local minimum $m_{k_2}$. We have
$\int_{B(m_{k_2},r)}|f_{k_2}^{(0)}|^2dx\ge 1/C$, while
$\int_{B(m_{k_2},r)}|f_{k_j}^{(0)}|^2dx={\cal O}(e^{-1/(Ch)})$ for
$j>2$. Taking the $L^2$-product with
$1_{B(m_{k_2},r)}f_{k_2}^{(0)}$, we get
$$
\frac{|u_{k_2}|}{C}+{\cal O}(|u_{k_1}|)\le C(\| u \|
+e^{-\frac{1}{Ch}}\|u_\bullet\|),\hbox{ so }|u_{k_2}|\le \widetilde{C}(\| u \|
+e^{-\frac{1}{Ch}}\|u_\bullet\|)
.
$$
Continuing this procedure, we get $\|u_\bullet\|\le
C(\|u\|+e^{-1/(Ch)}\|u_\bullet\|)$ and here the last term can be
absorbed, so
$$
\sum |u_k|^2\le {\cal O}(1)\| \sum u_kf_k^{(0)} \| ^2,
$$
which shows the uniform linear independence of the system of quasimodes $f_k^{(0)}$.

\medskip
\par
Let $e_k^{(0)}=\Pi^{(0)} (f_k^{(0)})$ where $\Pi^{(0)} = {\cal O}(1)$ is the
spectral projection of $-\Delta_A^{(0)}$ onto $E^{(0)}$. Arguing as in
Section 11 of \cite{KFP2} {(or as in the proof of Theorem \ref{exp3})} and using (\ref{eq4}), we first see that
$$
e_k^{(0)} = f_k^{(0)} + {\cal O}\left(h^{-N_1} e^{-(S_k - \alpha)/h}\right)\quad \mbox{in}\,\, L^2,
$$
for some small fixed $\alpha>0$ and $N_1>0$, and it is then clear that
also the system $(e_k^{(0)})$ is linearly independent in $L^2$,
uniformly with respect to $h$.

\par
In what follows, we shall need the quasimodes for the operator
$-\Delta_A^{(1)}$, associated to the saddle points of $\phi$. The
discussion here will be exactly the same as in the generic case,
described in detail in the following section.

\section{Multiple well analysis in the generic case}
\label{genc}
\setcounter{equation}{0}
In this section we shall be concerned with the generic case and shall show that the analysis of the singular values of
$d_\phi$ by Helffer, Klein and Nier in \cite{HKN} in the Witten case can be applied here, thanks to the
self-adjointness of $-\Delta^{(0)}$ in $(.|.)_{A,\kappa }$. We shall consider the case when
$e^{-\phi /h} \in L^2(\mathbb{R}^n)$ and assume from now on that we are in
the case 1 of Proposition \ref{ps1} for every saddle point $s_j$.

\subsection{The critical points and quasimodes in the generic case}
\label{labcut}
Let $\phi$ be a Morse function satisfying the same assumptions as in the beginning of Section \ref{lbl}.
The main hypothesis of this section is the following,
{\begin{hyp}\label{genc0} 
For every critical component $E_k=E_{m_k}$ as in Subsection
\ref{SSP_and_CC} we assume that
\begin{itemize}
\item ${{\phi }_\vert}_{E_k}$ has a unique point of minimum ($m_k$),
\item if $\mathrm{SSP}\cap E_k\ne \emptyset $, there is a unique ssp
  $s\in \mathrm{SSP}\cap E_k$, such
  that $\sup \phi (\mathrm{SSP}\cap E_k)=\phi (s)$ and in particular $E_k\cap
  \phi ^{-1}(]-\infty ,\phi (s)[)$ is the union of two distinct ccps.
\end{itemize}
\end{hyp}}
Combining this assumption with (\ref{sy.2.5}), we observe that $\kappa
(m_k)=m_k$ and $\kappa (s_j)=\kappa (s_j)$, for all $1\le k\le n_0$,
$1\le j\le n_1$.

In the following we let $s_0=\infty$.
The general labelling procedure described in the previous section simplifies in the generic case, since here
$N_{k^\sigma}=1$ for all $k^{\sigma}$, and hence the elements of the sets LM and CC can be labelled by positive integers.
We get the following result.

\begin{prop} \label{hyplab}
There exists an injective function
$$
\set{ 1, ..., n_0} \ni k \longmapsto j(k) \in \set{0, ...,  n_1} \ \ \
\ \hbox{with }j(1) = 0,\, \{s_{j(2)},\ldots\, ,s_{j(n_0)}\} = {\rm SSP},
$$
and a family of connected sets $E_k$, for $k \in \set{ 1, ..., n_0}$, such that the following pro\-perties hold:
\begin{enumerate}[i)]
\item We have $E_1=\mathbb{R}^n$, and $\overline{E_k}$ is compact for $k>1$. For every $k\geq 2$,
$E_k$ is the connected component containing $m_k$ in
$$
\set{ x\in \mathbb{R}^n;\, \phi(x) < \phi(s_{j(k)})},
$$
and $\phi (m_k)=\min_{E_k}\phi $. 
\item If $s_j\in E_k$ and $j=j(k')$  for some $k,k'\in \{ 1,...,n_0\}$
  and $j\in \{1,...,n_1\}$, then $k'>k$.
\end{enumerate}
\end{prop}


\par Following \cite{HKN}, we shall now introduce suitable quasimodes,
adapted to the local minima of $\phi$ and the simplified labelling, described in Proposition \ref{hyplab}. The construction can
be viewed as a special case of the general quasimode construction described in Section \ref{lbl}. In what follows, we let $\uuu$ stand for the set
of the critical points of $\phi$.

Let $\eps_1 >0$ be such that the distance between critical points is larger than $10 \eps_1$, and such
that for every critical point $c$ and $k \in \set{ 1, ... , n_0}$ we
have either $c \in \overline{E_k}$ or ${\rm dist}(c,\overline{E_k})\geq 10 \eps_1$.

When $0 < \eps_0 < \eps_1$ and $C_0 >1$ are to be defined
later, we build a family of $C_0^\infty$-cutoff functions
$\chi_{k,\eps}$, $k \in \set{1,..., n_0}$, $0<\eps < \eps_0$, in the following way: first
define $\chi_{1,\eps} = 1$. For $k \geq 2$, we consider the open
set $E_{k,\eps} := E_k \setminus \overline{B(s_{j(k)}, \eps)}$, and let
$\chi_{k,\eps}$ be a smooth function supported in $E_{k,\eps} + B(0,\eps/C_0)$ and
equal to $1$ in $E_{k,\eps} + B(0,\eps/(2C_0))$.  Then for $\eps_0$ small
enough and $C_0$ large enough, there exists $C$ such that for all
$0< \eps < \eps_0$, we have the following properties,
  \begin{enumerate}[(a)]
  \item $\chi_{k, \eps}$ is supported in $ E_k + B(0,\eps )$ and
    $\supp \chi_{k,\eps} \cap \set{ \phi < \phi(s_{j(k)})} \subset
    E_k$.
  \item The distance from any local minimum and a separating saddle point $c$ other than $s_{j(k)}$ to $\supp \nabla \chi_{k,\eps}$ is bounded from 
below by $3 \eps_1$. In addition, if $c \in \supp \chi_{k,\eps}$ then $c \in E_k$.
  \item There exists $C_\eps>0$ such that for all $x \in \supp \nabla (\chi_{k,\eps})\backslash B(s_{j(k)}, \eps)$, we have 
           $$
           \phi(s_{j(k)}) + C_\eps^{-1} \leq \phi(x) \leq \phi(s_{j(k)}) + C \eps.
           $$
         \item for any $k' \in \set{1, ..., n_0}$, if $m_{k'} \in
           \supp \chi_{k,\eps}$ then $k' \geq k$. In case $k' \neq k$
           we have
       $$
       \phi(m_{k'}) > \phi(m_k) \ \ \textrm{ and } \ \ \phi(s_{j(k')})
       < \phi(s_{j(k)}).
       $$
     \item for any $j' \in \set{ 1,..., n_1}$, if $s_j \in \supp
       \chi_{k,\eps}$ then either $j'$ does not belong to the image of
       the map $j$, or there exists $k' \geq k$ such that
       $$
       m_{k'} \in \supp  \chi_{k,\eps} \ \  \textrm{ and } \ \ j' = j(k').
       $$
       \item Inside $B(s_{j(k)}, \eps)$ we have:
       \begin{enumerate}[i)]
       \item The distance from $\supp \chi_{k,\eps} \cap B(s_{j(k)},
         \eps)$ to the projection $\pi_x(K_+)$ of the outgoing
         manifold $K_+$ is bounded from below by a constant $\delta_\eps>0$,
       \item for all $x \in B(s_{j(k)}, \eps)$ we have $\abs{ \phi(x)
           - \phi(s_{j(k)})} \leq C \eps^2 $.
                \end{enumerate}
 \end{enumerate}

\par
The quasimodes associated to the minima $m_k$ are introduced as follows,
\begin{equation} \label{fk0}
f_k^{(0)}= h^{-n/4} c_k(h)  e^{-{1\over h}(\phi (x)-\phi (m_k))}\chi _{k,\eps} (x),\quad 1\leq k \leq n_0,
\end{equation}
where $c_k(h)>0$ is a
normalization constant such that $\norm{ f_k^{(0)}}_\kappa^2 = 1$. Notice that $c_k(h) \sim c_{k,0}+ h c_{k,1}+\ldots$, with $c_{k,0}\neq 0$. 
As before, to these quasimodes we associate their projections to $E^{(0)}$,
\begin{equation}\label{ek0}
e_k^{(0)} = \Pi^{(0)}(f_k^{(0)}).
\end{equation}
In the same spirit we also define the quasimodes associated to the saddle points. For this we suppose that for all
$j \in \set{ 1, ..., n_1}$, the cutoff function  $\theta_j$  is  supported in $B(s_j, 2 \eps_1)$ and equal to $1$ in
$B(s_j,  \eps_1)$. The corresponding quasimodes and projections are defined as in (\ref{qms.4}), (\ref{ps.2bis}). We call 
\ekv{fk1}
{
f_j^{(1)}= (h^{-n/4} a_j(x,h)+r_j)e^{-\phi _{+,j}(x)/h}\theta_j(x),
}
the quasimode, where $a_j(x, h)$ is the vector already introduced, so that $\norm{ f_j^{(1)}}_{A,\kappa} = 1$.
Again we denote by $e^{(1)}_j$ the projection to $E^{(1)}$,
\begin{equation}\label{ek1}
e_j^{(1)} = \Pi^{(1)} (f_j^{(1)}).
\end{equation}
It is of the form (\ref{ps.2bis}) for a new tighter cut-off $\theta_j$. 
\subsection{Estimates for the quasimodes}\label{estqm}

We shall first work in the space $E^{(0)}$. The following result says that, modulo exponentially small errors,
the family $(e_k^{(0)})$ forms an orthonormal basis in this space.

\begin{prop}\label{estqm1} The space $E^{(0)}$ is spanned by the family $(e_k^{(0)})_{k = 1, ... , n_0}$ defined in {\rm (\ref{fk0}-\ref{ek0})}, and there exists
$\alpha>0$, independent of $\eps$, such that  for all $k,k' \in \set{1, ..., n_0}$, and $0 < \eps < \eps_0$, we have
\begin{equation}
( e_k^{(0)} | e_{k'}^{(0)})_\kappa  = \delta_{k, k'} + \ooo( e^{-\alpha/h}).
\end{equation}
\end{prop}

\preuve   We shall first repeat some arguments of Section 11 of \cite{KFP2}, essentially the same as in the proof of Theorem \ref{exp3}. Consider
$k \in \set{1, ..., n_0}$. Since
\begin{equation} \label{fk0bis}
f_k^{(0)}= h^{-n/4} c_k(h)  e^{-{1\over h}(\phi (x)-\phi (m_k))}\chi _{k,\eps} (x)
\end{equation}
where $\chi_{k,\eps}$ is the cutoff function defined above, we have
\begin{equation} \label{comhh}
\begin{split}
P (f_k^{(0)}) & = [P, \chi _{k,\eps}] \sep{  h^{-n/4} c_k(h)  e^{-{1\over h}(\phi (x)-\phi (m_k))}} \\
& = \ooo( h^{-N_0} e ^{- (S_k - C \eps)/h} ), \quad N_0>0,
\end{split}
\end{equation}
where $S_k = \phi (s_{j(k)})-\phi (m_k)$.
Here we used (d) and (f) of Subsection \ref{labcut}, where we gave estimates on the support of $\nabla \chi _{k,\eps}$.
Now we proved in Theorem 8.4 of \cite{KFP2}, and recalled in Theorem \ref{as2} that
$$
\Pi^{(0)} = \frac{1}{2 \pi i } \int_\gamma (z - P)^{-1} dz = \ooo(1),\quad
\Vert (z-P)^{-1}\Vert={\cal O}\left(\frac{1}{h}\right), z\in \gamma ,
$$
where $\gamma$ is an oriented circle of center $0$ and radius $h/C$ with $C$ large and fixed.   This implies that
\begin{equation}\label{ek0bis}
e_k^{(0)} = \Pi^{(0)} f_k^{(0)} = f_k^{(0)} + \ooo(h^{-N_1} e ^{- (S_k - C \eps)/h}) \ \  \textrm{ in } \  L^2,
\end{equation}
thanks to  the following equalities:
\begin{equation}
\begin{split}
(z - P) f_k^{(0)} &= z f_k^{(0)} - r_k \\
(z - P)^{-1} f_k^{(0)} &= \frac{1}{z}  f_k^{(0)} + (z-P)^{-1} z^{-1} r_k,
\end{split}
\end{equation}
where $r_k$ is defined by (\ref{comhh}).
Now since $f_k^{(0)}$ is normalized in $L^2$ in the sense that
  $(f_k^{(0)}|f_k^{(0)})_\kappa =1$, we get
$$
( e_k^{(0)} |  e_k^{(0)} )_\kappa = ( f_k^{(0)} |  f_k^{(0)} )_\kappa +
\ooo(h^{-N_2} e ^{- (S_k - C \eps)/h}) = 1 + \ooo(h^{-N_2} e ^{- (S_k - C \eps)/h})
$$
and for all $k$, $k' \in \set{1, n_0}$, we can write
\begin{equation} \label{efee}
( e_k^{(0)} |  e_{k'}^{(0)} )_\kappa = ( f_k^{(0)} |  f_{k'}^{(0)} )_\kappa +
\ooo(h^{-N_3} e ^{- (\min(S_k, S_{k'}) - C \eps)/h})
\end{equation}

\par
Let now $k' \neq k$, and suppose that $k'>k$ to fix the ideas. From
Subsection \ref{labcut}, we see that there are only three possible cases:
\begin{itemize}
\item $m_{k'} \in \supp \chi_{k,\eps}$ (case (d) in Subsection
  \ref{labcut}). Then
$$
 ( f_k^{(0)} |  f_{k'}^{(0)} )_\kappa = \ooo \sep{ h^{-N_4} e ^{- (\phi(m_{k'}) - \phi(m_k))/h} },
$$
since the support of $\chi_{k',\eps}$ is  included in the support of $\chi_{k,\eps}$. 
In that case we choose $ 0<  \alpha < \phi(m_{k'}) - \phi(m_k)$.
\item $m_{k} \in \supp \chi_{k',\eps}$. This case is the same as above.
\item The distance from $m_{k'}$ to $\supp \chi_{k,\eps}$ is larger
  than $3 \eps_1$ (case (b) in Subsection \ref{labcut}).
  Then 
  by construction, the supports of $\chi_{k,\eps}$ and $\chi_{k',\eps}$
  are disjoint and hence $ ( f_k^{(0)} | f_{k'}^{(0)} )_\kappa = 0$.
\end{itemize}
In all cases we can find $\alpha > 0$ independent of $\eps$ and $h$ such that
$$
 ( f_k^{(0)} |  f_{k'}^{(0)} )_\kappa = \ooo \sep{  e ^{- \alpha /h} }.
$$
Inserting this estimate in (\ref{efee}), and  possibly  shrinking $\eps_0$ and $\alpha$ so that $0<\alpha < (\min(S_k, S_{k'}) - C \eps_0)$, we
get,
\begin{equation}
( e_k^{(0)} |  e_{k'}^{(0)} )_\kappa =
\ooo( e ^{- \alpha/h}),\quad k'\neq k.
\end{equation}
The proof is complete. \hfill{$\Box$}

\begin{prop}\label{estqm2} The space $E^{(1)}$ is spanned by the family $(e_j^{(1)})_{j = 1, ... , n_1}$ defi\-ned in {\rm (\ref{fk1})}, {\rm (\ref{ek1})},
and there exists $\alpha'>0$ (independent of $\eps$) such that  for all $j,j' \in \set{1, ..., n_1}$,
\begin{equation}
(e_j^{(1)} | e_{j'}^{(1)})_{A, \kappa}  = \delta_{j, j'} + \ooo( e^{-\alpha'/h}).
\end{equation}
and $a_j(h) = a_{j,0} + h a_{j,1} + ...$.
\end{prop}

\preuve The proof of this result is simpler than that of Proposition \ref{estqm1}, thanks to the localization properties of
the cut-off functions $\theta_j$. Recall that in (\ref{fk1}) we defined
the functions
\ekv{fk1bis}
{
f_j^{(1)}= (h^{-n/4}a_j(x;h)+r_j(x))e^{-\phi _{+,j}(x)/h}\theta_j(x),
}
where $a_j(x,h) =  a_{j,0} + h a_{j, 1} + ...$ and $\norm{ f_j^{(1)}}_{A,\kappa} = 1$.
Again we denote by $e^{(1)}_j$ its projection onto $E^{(1)}$. Being in the case i) of Proposition \ref{ps1}, we can
choose the coefficient $a_{j,0}$ such that
\ekv{ps.ter}
{
(\kappa ^{\mathrm{t}}a_{j,0}(s_j)|a_{j,0}(s_j))_A = 1.
}
It follows from \cite{KFP2} that
\begin{equation*}
P^{(1)} (f_j^{(1)}) = \ooo( h^{-N_1} e ^{-\alpha' /h} ), \ \alpha '>0,\, N_1>0.
\end{equation*}
Here we may also recall that $\phi _{+,j}(x)\asymp |x-s_j|^2$ near $s_j$. Now using again Theorem 8.4 of \cite{KFP2} in the matrix case, 
or Theorem \ref{exp3}, we have
$$
\Pi^{(1)} = \frac{1}{2 \pi i } \int_\gamma (z - P^{(1)})^{-1} dz = \ooo(1),
$$
where we recall that $\gamma$ is an oriented circle of center $0$ and radius $h/C$ with $C$ large and fixed.
This implies that for some $N_2>0$,
\begin{equation} \label{ek0bibiss}
e_j^{(1)} = \Pi^{(1)} f_j^{(1)} = f_j^{(1)} + \ooo(h^{-N_2} e ^{- \alpha'/h}) \ \  \textrm{ in } \  L^2,
\end{equation}

Since the functions   $f_j^{(1)}$ are $L^2$-normalized (for the
$A,\kappa $-inner product),  we get for all $j$, $j' \in
\set{ 1, ..., n_1}$,
$$
( e_j^{(1)} | e_{j'}^{(1)} )_{A,\kappa } = ( f_j^{(1)} | f_{j'}^{(1)} )_{A,\kappa} +
\ooo(h^{-N_3} e ^{- \alpha'/h}) = \delta_{j,j'} + \ooo(h^{-N_3} e ^{- \alpha'/h}),
$$
where we also used that the supports of the functions $\theta_j$ are disjoint, since the distance between different critical points is
larger than $10 \eps_1$. Replacing $\alpha'$ by $\alpha'/2$ completes the proof. \hfill{$\Box$}


\subsection{Singular values}
Recall that in Section \ref{re}, as a consequence of Proposition \ref{ir1}, we showed that $d_{\phi}: E^{(0)}\rightarrow E^{(1)}$.
Following \cite{HKN}, we shall now estimate the singular values of this map. The first step is the
following proposition.

\begin{prop} \label{sing1}
For all $ k \in \set{ 1, ..., n_0}$,  $j \in \set{ 1, ..., n_1}$ and $0<\eps<\eps_0$,  we have
\begin{equation}
\begin{split}
& ( f_j^{(1)} | d_\phi f_k^{(0)} )_{A, \kappa} = 0 \ \ \textrm{ if } \ \ j \neq j(k), \\
& ( f_{j(k)}^{(1)} | d_\phi f_k^{(0)} )_{A, \kappa} = h^{1/2} b_k(h) e^{-S_k/h},
\end{split}
\end{equation}
where $b_k(h) \sim b_{k,0} + h b_{k,1} + ...$ . {Here $b_{k,0}$ is $\ne
0$ and will be determined in Proposition \ref{endestim}.} 
\end{prop}

\preuve
 For the first result we notice that for all  $k  \in \set{ 1, ..., n_0} $,
 \begin{equation} \label{comhhh}
\begin{split}
d_\phi f_k^{(0)} & =  ( h d^{(0)} \chi _{k,\eps}(x)) \sep{  h^{-n/4} c_k(h)  e^{-{1\over h}(\phi (x)-\phi (m_k))}}
\end{split}
\end{equation}
is supported in $\supp \nabla \chi _{k,\eps}$ which is disjoint from
the support of $\theta_j$ for all $j \neq j(k)$ according to point (c)
in Subsection \ref{labcut}.  As a consequence, the forms $ f_j^{(1)}$ and
$d_\phi f_k^{(0)}$ also have disjoint support, and we get the first result, $(f_j^{(1)} | d_\phi f_k^{(0)} )_{A, \kappa} = 0$.

\par Now we consider the case when $j= j(k)$. We shall use Gaussian type integrals in transversal directions, as in Section 11 of
\cite{KFP2}. We first notice that the domain of integration in the inner product can be reduced to $B(s_{j}, \eps)$: from
(\ref{comhhh}) and (\ref{fk1}) we have
\begin{equation}
\begin{split}
( f_{j}^{(1)} | d_\phi f_k^{(0)} )_{A, \kappa} = h^{1-n/2} c_k(h)
\int_{B(s_{j}, 2\eps_1)} &  \sep{\sep{ \kappa^{\mathrm{t}} (a_{j}(\kappa(x),h)) | d^{(0)} \chi
    _{k,\eps}(x)}_A+{\cal O}(h^\infty ) }  \\
 & \ \ \  \times \theta_j \circ \kappa(x)
e^{-(\phi _{+,j} \circ \kappa (x) + \phi(x) - \phi(m_k))/h} dx.
\end{split}
\end{equation}
Now we observe that for any  $ x \in B(s_{j}, 2\eps_1) \setminus  B(s_{j}, \eps )  $ we have
by (c) of Subsection \ref{labcut}  that
 if $x \in \supp  \nabla \chi_{k,\eps}$  then
 $$
           \phi(s_{j(k)}) + C_\eps^{-1} \leq \phi(x) \leq \phi(s_{j(k)}) + C \eps.
 $$
which implies
$$
\phi _{+,j} \circ \kappa (x) + \phi(x) - \phi(m_k) \geq  \phi(x) - \phi(m_k) = S_k + \phi(x) - \phi(s_{j(k)}) \geq S_k +  C_\eps^{-1}.
$$
We can therefore restrict the domain of integration to the smaller ball $B(s_j,\eps)$, modulo an error term as follows,
\begin{equation} \label{step1}
\begin{split}
( f_{j}^{(1)} | d_\phi f_k^{(0)} )_{A, \kappa} = &  h^{1-n/2} c_k(h)
\int_{B(s_{j}, \eps)}   \sep{\sep{ \kappa ^{\mathrm{t}}(a_{j})(\kappa(x),h) | d^{(0)} \chi
    _{k,\eps}(x)}_A +{\cal O}(h^\infty )} \\
 & \ \ \ \ \ \ \ \ \  \times \theta_j \circ \kappa(x)
e^{-(\phi _{+,j} \circ \kappa (x) + \phi(x) - \phi(m_k))/h} dx \\
& + \ooo ( h^{-N} e^{ -( S_k +  C_\eps^{-1})/h}).
\end{split}
\end{equation}

We shall need the following result.

\begin{lemma} \label{geom}
In $B(s_{j(k)}, \eps)$, we have  $ \phi_{+,j} \circ \kappa (x) + (\phi(x) -\phi (s_j))\asymp {\rm dist\,}(x,\pi _x(K_-))^2$.
\end{lemma}

\proof Referring to the end of Section \ref{re} and especially (\ref{2w.10.3}), we first recall that
\ekv{2w.10.3bb}
{
\phi _+  (x) - (\phi(x) -\phi (s_j))\asymp {\rm dist\,}(x,\pi _x(K_+))^2,
}
and
\ekv{2w.10.3bb-}
{
-\phi _- (x) + (\phi(x) -\phi (s_j))\asymp  {\rm dist\,}(x,\pi _x(K_-))^2,
}
Now recall the definition of the adjoint operator $-\Delta _{A^{\mathrm{t}}}=(-\Delta _A)^{A^{\mathrm{t}},*}$, from (\ref{su.6}).
Its principal symbol is $p_2-ip_1+p_0=p(x,-\xi )=:\check{p}(x,\xi )$, the
corresponding real "$q$"-symbol is $\check{q}(x,\xi )=q(x,-\xi )$, and
since this is the same type of operator,  the geometric discussion
above applies (notice in particular that Proposition \ref{ps1bis}
applies). We also notice for future reference that the hypotheses \ref{ny0}, \ref{ny2}
remain valid for $-\Delta_{A^{\mathrm{t}}}$. Moreover, if we are in the case 1 of Proposition
\ref{ps1} for $-\Delta _A$, then we are so also for $-\Delta
_{A^{\mathrm{t}}}$. In particular we can introduce
$\Lambda _{\phi _+^*}$, $\Lambda _{\phi _-^*}$  as the stable outgoing and
incoming $H_{\check{q}}$-invariant Lagrangian manifolds through $s_j$
and $K_\pm^*\subset \Lambda _\phi $ as the stable outgoing and incoming manifolds for
${{H_{\check{q}}}_\vert}_{\Lambda _\phi }$ (noting that $\check{q}=0$ on
$\Lambda _\phi $). Again we have  ${\rm dim\,}K_+^*=n-1$, ${\rm dim\,}K_-^*=1$ and
$$
\phi _+^*-(\phi -\phi (s_j))\asymp {\rm dist\,}(x,\pi _x(K_+^*))^2,
$$
and
\ekv{2w.10.4bb}
{
\ - \phi_-^* + (\phi -\phi (s_j))\asymp {\rm dist\,}(x,\pi _x(K_-^*))^2.
}

\par In view of the general relation
\ekv{2w.10.5}
{
J_*(H_q)=-H_{\check{q}},\hbox{ where }J:(x,\xi )\mapsto (x,-\xi ),
}
we see that $\Lambda _{\phi _-^*}=J(\Lambda _{\phi _+})$, $\Lambda _{\phi
_+^*}=J(\Lambda _{\phi _-})$, or simply
\ekv{2w.10.6}
{
\phi _-^*=-\phi_+,\ \phi _+^*=-\phi_- .
}

From the relation (\ref{sy.6}),
$$
-\Delta _A\circ U_\kappa =U_\kappa \circ (-\Delta _{A^{\mathrm{t}}}),
$$
and Egorov's theorem, we see that $\check{p}=p\circ {\cal K}$ and
hence $\check{q}=q\circ {\cal K}$, where
${\cal K}:(y,\eta )\mapsto (\kappa (y),(\kappa ^{\mathrm{t}})^{-1}\eta
)$ is the natural canonical transformation which lifts $\kappa $ to
the cotangent space. Since $\phi $ is invariant under composition with
$\kappa $, we have ${\cal K}:\Lambda _\phi \to \Lambda _\phi $ and it
follows that $K_\pm^*={\cal K}(K_\pm)$ and in particular $\pi
_x(K_\pm^*)=\kappa (\pi _x(K_\pm))$. Combining (\ref{2w.10.4bb}),
(\ref{2w.10.6}), we get
$$
\phi _+(x)+(\phi (x)-\phi (s_j))\asymp \mathrm{dist\,}(x,\pi
_x(K_-^*))^2\asymp \mathrm{dist\,}(x,\kappa (\pi _x(K_-)))^2,
$$
and replacing $x$ by $\kappa (x)$, we get the lemma. \hfill{$\Box$}

\medskip\par
We summarize for further use some of the microlocal results obtained during the proof of the preceding lemma:

\begin{lemma} \label{sumjk} Let $\kkk : (y, \eta) \longmapsto
  (\kappa(y), (\kappa^{\mathrm{t}})^{-1} \eta)$ be the lifting of
  $\kappa$ and $J : (x,\xi) \longmapsto (x,-\xi)$. We have
\begin{equation} \label{jjjj}
\Lambda _{\phi _-^*}=J(\Lambda _{\phi _+}), \ \ \ \Lambda _{\phi
_+^*}=J(\Lambda _{\phi _-}) \ \ \ \textrm{ and } \
\phi _\mp^*=-\phi_\pm,
\end{equation}
and also
\begin{equation} \label{kkkk}
\Lambda _{\phi _\pm^*}=\kkk(\Lambda _{\phi _\pm}), \ \  \pi
_x(K_\pm^*)=\kappa (\pi _x(K_\pm)),
\ \ \ \textrm{ and } \ \phi _\pm^*= \phi _\pm \circ \kappa
\end{equation}
\end{lemma}

\noindent \it End of the proof of Proposition \ref{sing1}. \rm \  Now we are able to continue the computation in (\ref{step1}). Using Lemma \ref{geom}, we get that the exponential term in the integral satisfies
$$
\phi_{+,j} \circ \kappa (x) + (\phi(x) -\phi (m_k)) - S_k \asymp  {\rm dist\,}(x,\pi _x(K_-))^2,
 $$
on $B(s_j, \eps) \cap \supp  \nabla \chi_{k,\eps}$.
{By} evaluating the
Gaussian integral in the directions transversal to $\pi _x(K_-)$, we get
\begin{equation} \label{step2}
\begin{split}
( f_{j}^{(1)} | d_\phi f_k^{(0)} )_{A, \kappa} = &  h^{1/2} {\widetilde{b}}_k(h) e^{- S_k/h}
 + \ooo ( h^{-N} e^{ -( S_k +  C_\eps^{-1})}) \\
=  & h^{1/2} b_k(h) e^{- S_k/h}.
\end{split}
\end{equation}
{It is also quite clear that $b_k$ is elliptic and a more
  detailed computation of $b_{k,0}$ will be given in Section \ref{explt}}.\hfill{$\Box$}

\medskip
\par
We are now able to compute the matrix of $d_\phi$ from $E^{(0)}$ to $E^{(1)}$ with respect to  the bases $(e^{(0)}_k)$ and $(e^{(1)}_j)$.
This result is contained in the following proposition.

\begin{prop} \label{sing2} There exists $\alpha'' >0$ such that if
  $\eps_0$ is small enough then for all $ k \in \set{ 2, ...,
    n_0}$, $j \in \set{ 1, ..., n_1}$ and $0<\eps<\eps_0$, we have
\begin{equation}
\begin{split}
& ( e_j^{(1)} | d_\phi e_k^{(0)} )_{A, \kappa} =  \ooo \sep{ e^{-(S_k + \alpha'')/h}} \ \ \textrm{ if } \ \ j \neq j(k), \\
\textrm{ and } \ \ \ & ( e_{j(k)}^{(1)} | d_\phi e_k^{(0)} )_{A, \kappa} = h^{1/2} b_k(h) e^{-S_k/h},
\end{split}
\end{equation}
where $b_k(h) \sim b_{k,0} + h b_{k,1} + ...$ is
  elliptic for $k\ne 1$, and $( e_j^{(1)} | d_\phi e_k^{(0)} )_{A, \kappa} = 0$ when $k=1$.
\end{prop}

\preuve Let $k$ and $j$ be integers as above. The case $k= 1$ is clear
since $ d_\phi e_1^{(0)} =0$. Consider now the case when $k \geq 2$:
using the self-adjointness and intertwining properties of the
projectors $\Pi^{(l)}$,  $l=0,\, 1$, which follow from Proposition \ref{ir1}, we check that
\begin{equation}
\begin{split}
( e_{j}^{(1)} | d_\phi e_k^{(0)} )_{A, \kappa}
 & = ( e_{j}^{(1)} | d_\phi \Pi^{(0)} f_k^{(0)} )_{A, \kappa}
  = ( e_{j}^{(1)} |\Pi^{(1)} d_\phi  f_k^{(0)} )_{A, \kappa} \\
  & = (\Pi^{(1)} e_{j}^{(1)} | d_\phi  f_k^{(0)} )_{A, \kappa}
   = ( e_{j}^{(1)} | d_\phi  f_k^{(0)} )_{A, \kappa} \\
 & \ \ \ \ \ \ \ \ \  = ( f_{j}^{(1)} | d_\phi  f_k^{(0)} )_{A, \kappa}  +    ( e_{j}^{(1)} - f_{j}^{(1)} | d_\phi  f_k^{(0)} )_{A, \kappa}.
 \end{split}
 \end{equation}
 According to Proposition \ref{sing1}, the result of the proposition will follow if we show that  for any $k$ and $j$, we have
 \begin{equation} \label{remain}
  ( e_{j}^{(1)} - f_{j}^{(1)} | d_\phi  f_k^{(0)} )_{A, \kappa}  = \ooo \sep{ e^{-(S_k + \alpha'')/h}}.
  \end{equation}
Let us therefore prove (\ref{remain}). By Cauchy--Schwarz,
\begin{equation} \label{remain2}
 \abs{  ( e_{j}^{(1)} - f_{j}^{(1)} | d_\phi  f_k^{(0)} )_{A, \kappa} } \leq C(A)\norm{  e_{j}^{(1)} - f_{j}^{(1)}}_{L^2}  \norm{ d_\phi  f_k^{(0)}}_{L^2} .
 \end{equation}
Now a computation very similar to the one made in the proof of (\ref{comhh}) and based on the support of $\nabla \chi_{k,\eps}$ gives that
$$
d_\phi  f_k^{(0)} = \ooo \sep{ h^{-N_0} e^{-(S_k - C \eps)/h} }
$$
where $C$ is the constant appearing in (c) and (f) of Subsection \ref{labcut}. On the other hand  we  have from (\ref{ek0bibiss}) that
$$
e_{j}^{(1)} - f_{j}^{(1)} = \ooo ( h^{-N_2} e^{- \alpha' /h} )\hbox{
  in }L^2,
$$
where $\alpha'>0$ is independent of $\eps$. Consequently, we can write
that
$$
 \abs{ ( e_{j}^{(1)} - f_{j}^{(1)} | d_\phi  f_k^{(0)} )_{A, \kappa} }
 = \ooo \sep{ h^{-(N_0 -N_2)} e^{-(S_k - C \eps + \alpha')/h}}.
 $$
 Taking $\eps_1$ small enough (e.g. such that $C \eps_1 \leq \alpha'/3$)  first and then posing $\alpha'' = \alpha'/2$, gives the estimate (\ref{remain}) and the proof of Proposition \ref{sing2} is complete. \hfill{$\Box$}

\subsection{Main result in the generic case}
We can now give the complete asymptotics for the exponentially small
eigenvalues in the generic case. Let us first prove the following result.

\begin{lemma}\label{genc1}
Consider the $n_1\times n_0$ matrix
$$
R = (r_{j,k}) := ( (e_j^{(1)} | d_\phi e_k^{(0)})_{A,\kappa})_{j \in \set{1,..,n_1}, \ k \in \set{1,..,n_0}},
$$
where we recall that
there exists $\alpha''>0$ such that for all $ k \in \set{ 2, ..., n_0}$, and $j \in \set{ 1, ..., n_1}$, 
$$
r_{j(k), k} = h^{1/2} b_k(h)e^{-S_k/h}, \ \ r_{j,k} =
\ooo(  e^{-(S_k + \alpha'')/h} ), \textrm{ when $j \neq j(k)$},
$$
while $r_{j,1} = 0$ for all $j$. Set also $r_{j(1), 1} = 0$.  Then
there exists $\eta >0$ such that the singular values
{$\nu_k(R)$}  of $R$ , {enumerated in a a suitable order,} satisfy
 $$
\nu_k(R) = |r_{j(k), k}|( 1 + \ooo(e^{-\eta/h}),\ 1\le k\le n_0.
$$
\end{lemma}

\preuve Since the first column of $R$ consists of zeros, we know that
${\nu _1=}0$ is a singular value of $R$ and
we can study the reduced matrix $R'$ with entries $ r_{j, k}' = r_{j, k+1}$,
where $k\geq 1$. We shall use that there is only one dominant
term in each column of $R'$.  Define the $(n_0-1)\times (n_0-1)$
diagonal matrix $D$ as follows,
$$
D = \textrm{diag}  ( h^{1/2} b_{k+1}(h) e^{-S_{k+1}/h}, \ k =1, ..., n_0-1).
$$
Notice that $D$ is invertible, thanks to the ellipticity of $b_k$, and that $\nu _k(D)=|r_{j(k),k}|$.
Define the characteristic matrix of $R'$ to have the entry $1$ at each dominant term (columnwise), and $0$ elsewhere:
$$
U = ({\delta_{j, j(k+1)}}),
$$
where $\delta $ is the Kronecker symbol. Then there is a constant
$\eta >0$ such that $ R' D^{-1} = U + \ooo( e^{-\eta/h}) $,
 \begin{equation} \label{MUD}
 R' = (U + \ooo( e^{-\eta/h})) D.
 \end{equation}
The Ky Fan inequalities therefore give
\begin{equation} \label{fan}
\nu_k(R') \leq (1 + \ooo( e^{-\eta/h})) \nu_k(D).
\end{equation}
To get the opposite estimate, notice that $U$ is isometric, $U^*U=1$, and write
$$
U^*R' = (1+ {\cal O}(e^{-\eta/h}))D,
$$
$$
D  = (1 + {\cal O}(e^{-\eta/h}))U^*R',
$$
to get 
$$
\nu _k(D)\le (1+{\cal O}(e^{-\eta /h}))\nu _k(R').
$$
\hfill{$\Box$}

\medskip\par
One of the main results of this paper is the following theorem.
 \begin{theo}
\label{theogene} {In addition to the general assumptions, we adopt the
Hypothesis \ref{genc0} and} 
   assume that we are in case i) of Proposition {\rm \ref{ps1}} at every
   {separating saddle point}. The exponentially small
   eigenvalues of $P$ are real and given by
 $$
 \mu_k= (h |b_k(h)|^2 +{\cal O}(h^\infty)) e^{-2S_k/h}, \ \ \ \  k = 1, ..., n_0,
 $$
 where $S_1 = + \infty$ (and $\mu_1 = 0$), by convention. {Here
 $b_k(h)\sim b_{k,0}+hb_{k,1}+...$, where $b_{k,0}\ne 0$ will be
 studied in more detail in Section \ref{explt}.}
 \end{theo}

{\proof
According to the propositions \ref{estqm1} and \ref{estqm2} the bases
$(e_k^{(0)})$ and $(e_j^{(1)})$ are orthonormal up to exponentially
small errors in $E^{(0)}$ and
$E^{(1)}$ respectively, for the $A,\kappa $ scalar products. Let
$(\widetilde{e}_k^{(0)})$ and $(\widetilde{e}_j^{(1)})$ be the
corresponding orthonormalizations (obtained by taking square roots of
the Gramians), which differ from the original bases by exponentially
small recombinations. Then with respect to the new bases, the matrix
of $d_\phi $ is $\widetilde{R}=(1+{\cal O}(e^{-\alpha /h}))R(1+{\cal
  O}(e^{-\alpha /h}))$ and from Lemma \ref{genc1} and the Ky Fan
inequlities (which will be used in a more essential way in the proof of
Theorem \ref{gl4}) we see that the conclusion of that lemma is also
valid for $\widetilde{R}$. The matrix of the restriction of $P$ to
$E^{(0)}$ with respect to the basis $(\widetilde{e}_k^{(0)})$ is
$\widetilde{R}^*\widetilde{R}$ and the theorem follows.
\hfill{$\Box$}}

\section{Explicit computation of the leading tunneling coefficient}
\label{explt}
\setcounter{equation}{0} The aim of this section to compute the
dominant term { in the amplitude} of the exponentially small eigenvalues. For this we
shall estimate the coefficient $b_k$ appearing in Theorem
\ref{theogene}. We follow essentially the proof given in
\cite{HeSj85}, but point out that here we are in a non-selfadjoint
situation. The main result is given in Proposition \ref{endestim} at
the end of the section.

\subsection{Geometric preliminaries} \label{geoo}

We recall some points concerning the study of the operator $-\Delta_A^{(l)}$. The principal symbol $p(x,\xi)$ of $-\Delta_A^{(l)}$ does not depend on $l$ and if $q(x,\xi) = - p(x, i \xi)$ then from
Subsection \ref{ssb}  we have
$$
q(x,\xi) = \seq{ A(\xi + d \phi(x))| (\xi - d \phi(x))} = ( (\xi + d \phi(x))| (\xi - d \phi(x)))_A,
$$
where we recall that $A : (\R^n)^* \rightarrow \R^n$, and where
$(\cdot |\cdot )_A$ was introduced in (\ref{ss.2}). Near a critical
point $(s_j,0)$ of index $1$, we also defined the two Lagrangian
manifolds $\Lambda_\pm = \set{ (x,\xi), \ \xi = d\phi_\pm(x)}$, on
which $q$ is equal zero. This gives the following   two eikonal equations in a neighborhood of $s_j$:
\begin{equation} \label{eik1}
0 = q(x, d \phi_{\pm, j}(x)) = ( d\phi_{\pm, j}(x) +d \phi(x)| d\phi_{\pm, j}(x) -d \phi(x))_A.
\end{equation}
We omit the index $j$ when no confusion is possible, and define in a
neighborhood of the saddle point $s_j$,
\begin{equation} \label{defg}
\left\{
\begin{array}{l}
g_+ = \phi_+ -(\phi -\phi(s_j)) \\
g_- = -\phi_- +(\phi -\phi(s_j))
\end{array} \right.
\end{equation}
Notice that $g_-$ is  the phase appearing in the Laplace integral
defining $b_k$ in formula \ref{step1} up to a constant; indeed in a neighborhood of $s_j$
 we have according to Lemma \ref{sumjk}
\begin{multline}
\phi_+ \circ \kappa + (\phi -\phi(m_k)) = \phi_+^* + (\phi -\phi(s_j)) + (\phi(s_j) - \phi(m_k)) \\ = -\phi_- + (\phi -\phi(s_j)) + (\phi(s_j) - \phi(m_k)) = g_- + (\phi(s_j) - \phi(m_k)).
\end{multline}
Using again Lemma \ref{sumjk} we also define
\begin{equation} \label{sstar}
\left\{
\begin{array}{l}
g_+^* = \phi_+^* -(\phi -\phi(s_j)) \\
g_-^* = -\phi_-^* +(\phi -\phi(s_j)) = \phi_+ + (\phi -\phi(s_j))
\end{array} \right.
\end{equation}
and with these notations the eikonal equation for $\phi_-$ can be written near $s_j$ in the following form :
\begin{equation}\label{geo0.1}
( d g_+^*(x)| d g_-(x))_A = 0,
\end{equation}
according to the following direct computation and (\ref{eik1})
\begin{equation}
\begin{split}
( d g_+^*(x)| d g_-(x))_A & = ( d \phi_+^*(x) -
d \phi(x)| -d \phi_-(x) + d \phi(x) )_A  \\
& =  ( d \phi_-(x) +
d \phi(x)| d \phi_-(x) - d \phi(x))_A = 0.
\end{split}
\end{equation}
Similarly using $\phi_+$, we get
\begin{equation}
( d g_-^*(x)| d g_+(x))_A = 0.
\end{equation}
 We now recall other properties of the functions $g_\pm$ and $g_\pm^*$, which are essentially reformulations of (\ref{2w.10.3}):
 In a neighborhood of $s_j$, we have
 \begin{equation} \label{dist}
 g_\pm(x) \asymp d( x, \pi_x(K_\pm))^2 \ \ \textrm{ and } \ \ g_\pm^*(x) \asymp d( x, \pi_x(K_\pm^*))^2
 \end{equation}

{Recall some properties of the Hamilton 
fields $\nu_\pm = \left. H_q \right|_{\Lambda_{\pm \phi}}$ near $s_j$
that we identify with their $x$-space projections.}
We already used the fact that
$$
\forall x \in \pi_x(K_\pm), \ \ \  \nu_+(x) = 2 A d\phi(x) \in T_x(\pi_x(K_\pm),
$$
and that it vanishes at $s_j$. Similarly, 
$$
\nu_- (x) = -2 A^{\mathrm{t}} d \phi(x)
$$
so that $\nu_- (x) = -\nu_+^*(x)$ and
\begin{equation} \label{bienutile}
\forall x \in \pi_x(K_\pm^*), \ \ \  \nu_-(x) = -2 \tA d\phi(x) \in T_x(\pi_x(K_\pm^*),
\end{equation}
where the $K_\pm^*$ are associated to the operator $-\Delta_{ \tA}$
(see the end of Section \ref{re}).

\par
When looking for an accurate expression of the eigenmodes of $P = -\Delta_A $ on $l$ forms we need to know precisely the conjugate operator
$$
P_+ = e^{\phi_+/h} P e^{-\phi_+/h},
$$
in a neighborhood of a saddle point $s_j$ (recall that we omit the index $j$ when no confusion is possible).
\begin{lemma} \label{pplus}
In a neighborhood of $s_j$ we have on $l$-forms
\begin{equation}
\begin{split}
P_+ = hd (hd)^{A,*} +  (hd)^{A,*} hd - &  \overbrace{( dg_-^*(x)| d g_+(x))_A}^{=0} \\
& + h \lll_{Ad(\phi + \phi_+)} + h \lll^{A, *}_{ \tA d(\phi - \phi_+)},
\end{split}
\end{equation}
where $\lll_\nu$ denotes the Lie derivative in the $\nu$ direction.
\end{lemma}
\proof
We first recall that $P = d_\phi d_\phi^{A,*} + d_\phi^{A,*} d_\phi$ on  forms (of arbitrary order) and we get
\begin{equation}
\begin{split}
P_+ = &   e^{\phi_+/h} (d_\phi d_\phi^{A,*} + d_\phi^{A,*} d_\phi)  e^{-\phi_+/h} \\
= &   \sep{ e^{\phi_+/h} d_\phi  e^{-\phi_+/h}}\sep{   e^{\phi_+/h} d_\phi^{A,*}  e^{-\phi_+/h}} +  \sep{ e^{\phi_+/h} d_\phi^{A,*}  e^{-\phi_+/h}} \sep{   e^{\phi_+/h} d_\phi  e^{-\phi_+/h} }\\
= &  \sep{ e^{\phi_+/h} d_\phi  e^{-\phi_+/h}}\sep{   e^{-\phi_+/h} d_\phi  e^{\phi_+/h}}^{A,*} +  \sep{ e^{-\phi_+/h} d_\phi  e^{\phi_+/h}}^{A,*} \sep{   e^{\phi_+/h} d_\phi  e^{-\phi_+/h} } \\
= &  \sep{ hd + d(\phi-\phi_+)^\wedge}\sep{ hd + d(\phi+\phi_+)^\wedge}^{A,*} \\
 & +\sep{ hd + d(\phi+\phi_+)^\wedge}^{A,*} \sep{ hd + 
d(\phi-\phi_+)^\wedge}.
\end{split}
\end{equation}
For a one form $\omega$  we have $ (\omega^\wedge )^{A,*} = (A \omega)^\rfloor$ so
\begin{multline} \label{ppluseq}
P_+ =    hd (hd)^{A,*} + (hd)^{A,*} hd \\
 + d(\phi-\phi_+)^\wedge(A d(\phi+ \phi_+))^\rfloor + (A d(\phi+ \phi_+))^\rfloor d(\phi-\phi_+)^\wedge \\
 + hd (A d(\phi+ \phi_+))^\rfloor + (A d(\phi+ \phi_+))^\rfloor hd  \\
 + d(\phi-\phi_+)^\wedge (hd)^{A,*} + (hd)^{A,*}d(\phi-\phi_+)^\wedge.
\end{multline}
Using the formula $\seq{\omega, \nu} = \nu^\rfloor \omega^\wedge +\omega^\wedge \nu^\rfloor $ for a vector field $\nu$ and a one-form $\omega$, we get
\begin{multline} \label{eikpp}
d(\phi-\phi_+)^\wedge(A d(\phi+ \phi_+))^\rfloor + (A d(\phi+ \phi_+))^\rfloor d(\phi-\phi_+)^\wedge \\
= ( d(\phi+ \phi_+)| d(\phi-\phi_+))_A   =  -( dg_-^*(x)| d g_+(x))_A = 0.
\end{multline}
Using also Cartan's formula $ \lll_\nu \omega = d( \nu^\rfloor \omega) + \nu^\rfloor d\omega$
with $\nu = A d(\phi+ \phi_+)$ we can write
\begin{equation} \label{liepp}
h \lll_{Ad(\phi + \phi_+)} = hd (A d(\phi+ \phi_+))^\rfloor + (A d(\phi+ \phi_+))^\rfloor hd
\end{equation}
and using also
  $$
   (\nu^\rfloor)^{A,*} =( \tA {}^{-1} \nu)^\wedge,
   $$
for a vector $\nu$,   we get
\begin{multline} \label{liestar}
h \lll_{ \tA d(\phi - \phi_+)}^{A,*} = \sep{  hd ( \tA d(\phi - \phi_+))^\rfloor} ^{A,*} + \sep{ ( \tA d(\phi- \phi_+))^\rfloor hd }^{A,*} \\
= \sep{ ( \tA d(\phi- \phi_+))^\rfloor} ^{A,*} (hd)^{A,*} + (hd)^{A,*} \sep{ ( \tA d(\phi- \phi_+))^\rfloor} ^{A,*}  \\
= d(\phi - \phi_+)^\wedge  (hd)^{A,*} + (hd)^{A,*} d(\phi -
\phi_+)^\wedge .
\end{multline}
Putting together (\ref{eikpp}-\ref{liepp}-\ref{liestar}) in (\ref{ppluseq}), and using the fact that the eikonal equation is satisfied by $\phi_+$,  we get the lemma.
\hfill{$\Box$}

\medskip
It is convenient to introduce adapted coordinates in a neighborhood
of the saddle point $s_j$. We begin with defining a $\ccc^\infty$ function $y_1$ by
$$
y_1 = \pm \sqrt{g_+(x)}
$$
with the sign $\pm$ depending on a choice of a side of the (local) hypersurface $\pi_x(K_+)= \set{ g_+(x) = 0 }$. We impose for example the sign $+$ on the side containing $m_k$ for $j=j(k)$ defined according to the injection defined in the preceding section. Of course we have
\begin{equation} \label{yun}
y_1^2 = g_+ = \phi_+-(\phi-\phi(s_j)).
\end{equation}
We also introduce $y_1^* = y_1 \circ \kappa$, so that according to Lemma (\ref{sumjk})
\begin{equation*} 
(y_1^*)^2 = (y_1 \circ \kappa)^2 = \phi_+^*-(\phi-\phi(s_j))=  g_+^*,
\end{equation*}
and we  notice that
$$
A d y_1^* \in T \pi_x(K_-) \ \ \ \textrm{ on } \ \pi_x(K_-).
$$
Now we consider the restriction of $g_-$ to $\pi_x(K_+)$. Since
$$
g_-(.)  \asymp d(., \pi_x(K_-))^2  \textrm{ near } s_j,
 $$
 we can apply the Morse lemma to $g_-$ on $\pi_x(K_+)$ and get smooth
 functions $z_2, ..., z_n$ on $\pi_x(K_+)$ such that,
$$
g_- = z_2^2 + ... +  z_n^2  \textrm{ on } \ \pi_x(K_+).
$$
{E}xtend the functions $z_k$ to a whole neighborhood of $s_j$ using the Hamilton Jacobi equation
$$
\nabla_{A dy_1^*} z_k = 0,
$$
or equivalently
$$
\seq{ d y_1^*, d z_k}_A = 0.
$$
Since
$ \seq{ d y_1^*, d g_-}_A = 0 $ {by (\ref{geo0.1})}, the function  $g_-$ does not depend on the $y_1$ coordinate and we have
$$
g_- = z_2^2 + ... +  z_n^2  \textrm{ in a whole neighboorhod of } s_j.
$$
We  introduce the variables $z_k^* = z_k \circ \kappa$, and we also have
$$
g_-^* = (z_2^*)^2 + ... +  (z_n^*)^2  \textrm{ in a neighboorhod of } s_j.
$$
Note that  according to (\ref{defg}) and (\ref{sstar}),
\begin{equation} \label{expphiplus}
\begin{split}
& \phi_+ = \frac{1}{2} ( g_+ + g_-^*) = \frac{1}{2} \sep{  y_1^2 + (z_2^*)^2 + ... + (z_n^*)^2 }, \\
& \phi- \phi(s_j) = \frac{1}{2} ( - g_+ + g_-^*) = \frac{1}{2} \sep{  - y_1^2 + (z_2^*)^2 + ... + (z_n^*)^2 }
\end{split}
\end{equation}
Note that $\phi \circ \kappa = \phi$ {gives}
$$
- y_1^2 + (z_2^*)^2 + ... + (z_n^*)^2 = - (y_1^*)^2 + z_2^2 + ... + z_n^2,
 $$
so that
\begin{equation}  \label{sumsim}
y_1^2 + z_2^2 + ... + z_n^2 = (y_1^*)^2 + (z_2^*)^2 + ... + (z_n^*)^2.
\end{equation}

\subsection{A twisted Hodge operator}

We introduce now a special Hodge operator.  We use  the notation of
\cite{HeSj85}. We recall first the definition of the usual Hodge
* operator, denoted $*$ or $\iii$ (and $*$ or $\iii_l$ when restricted
on $l$-forms):  For any real $l$-forms $\omega$ and $\mu$,
$$
(\omega | \mu) d x= \omega \wedge (*\mu) = \omega \wedge ( \iii_l
\mu),\quad dx=dx_1\wedge ...\wedge dx_n,
$$
where $(\cdot |\cdot )$ denotes the scalar product on $l$-forms inherited from the identification of $T \R^n$ and $T^* \R^n$ via the standard euclidian metric. We denote by
$$
\imath : T^* \R^n \longrightarrow T \R^n
$$
 the corresponding application whose matrix is the identity in euclidian coordinates. This corresponds to a choice of an additional structure in our problem.

 In fact the right object is the following Hodge operator denoted $*^A$
 or $\iii^A$, and whose restriction on $l$ forms will be denoted again
 by $*^A$ or $\iii_l^A$: For any real $l$-forms {$\omega$} and $\mu$,
 it is defined by
$$
(\omega| \mu )_A d x = \omega \wedge (*^A\mu) = \omega \wedge ( \iii_l^A \mu).
$$
 We give now some properties of the twisted Hodge operators $*^A$ and $*^{\tA}$.

\begin{lemma} \label{defs}
 We have on $l$-forms
 $$
 i) \ \ *^A = * \circ \wedge^l( \imath^{-1} \tA), \ \ \ \ \  \ \ ii) \ \ ( *^A)^{-1} = (-1)^{l(n-l)} \wedge^{l} ( \imath^{-1} \tA)^{-1} *,
 $$
 where $*$ denotes $\iii_l$ in the first  equality and $\iii_{n-l}$ in
 the second one. In addition if $\omega$ is an $l$-form and $\mu$ an
 $(n\! - \!l)$-form we have 
 $$
iii) \ \  (*^A \omega) \wedge  (*^A \mu) = (\det (\imath^{-1}\tA)) \omega \wedge \mu.
 $$
 \end{lemma}

\preuve
Let $\omega$ and $\mu$ be two $l$-forms. Then we can write
\begin{multline}
\omega \wedge *^A \mu = ( \omega| \mu)_A dx= \seq{ \wedge^l A \omega| \mu} dx = \seq{ \omega| \wedge^l \tA \mu} dx \\
 = \sep{ \omega| \wedge^l (\imath^{-1}\tA) \mu} dx = \omega \wedge  (*  \wedge^l (\imath^{-1}\tA) \mu )
 \end{multline}
 This proves i).

\par
The invertibility of $*^A$ and formula ii) are direct consequences of
i) and the fact that $\iii_l \circ \iii_{n-l} = (-1)^{l(n-l)}$ for the
usual Hodge operator. In order to prove iii) we recall first that
 {\begin{equation*}\begin{split}
 &(* \omega) \wedge  (* \mu) =(*\omega |\mu )dx=(\mu |*\omega )dx\\
&=\mu \wedge (**\omega )=(-1)^{l(n-l)}\mu \wedge \omega =\omega \wedge
\mu .
\end{split}
 \end{equation*}
 Using property i) yields the result since
 $$
((\wedge^l(\imath ^{-1}A))\wedge (\wedge^{n-l}(\imath
^{-1}A))=\wedge^n(\imath ^{-1}A)=\det (\imath ^{-1}A)\,\mathrm{Id} .
 $$}
 The proof is complete.
 \hfill{$\Box$}

\medskip\par
The second lemma is devoted to the commutation properties of $*^A$ with the Lie derivatives.

\begin{lemma} \label{starlie}
Let $\nu$ be a smooth vector field. We have
$$
*^A \lll_\nu^{A,*}  = -  \lll_\nu *^A \ \ \textrm{ and } \ \ \ \lll_\nu^{A,*} *^{\tA} = - *^{ \tA} \lll_\nu,
$$
\end{lemma}

\preuve{
Consider the identity
$$
{\cal L}_\nu ((u(x)|v(x))_Adx)={\cal L}_\nu (u\wedge (*^Av)).
$$
{ If $u,v$ have compact support, then the left hand side
is the Lie derivative of a compactly supported differential form and
has therefore the integral zero. On the other hand it takes the form
$({\cal L}_\nu u|v)_Adx-(u|Mv)dx$ for some first order differential
operator $M$, so by integration, we recognize that $M$ is the $A*$
adjoint of ${\cal L}_\nu $ and }
the left hand side becomes 
$$
({\cal L}_\nu u|v)_Adx-(u|{\cal L}_\nu ^{A,*}v)_Adx.
$$
The right hand side is equal to
$$({\cal L}_\nu u)\wedge (*^Av)+u\wedge ({\cal L}_\nu *^Av).$$
The first terms in the two expressions coincide and hence the second
terms also. This leads to
$$
-u\wedge (*^A{\cal L}_\nu ^{A,*}v)
=u\wedge ({\cal L}_\nu *^Av) 
$$
and the first identity in the lemma follows.}

\par {
In order to prove the second identity we shall first prove the general
identity of independent interest:
\begin{equation}\label{**}
*^A*^{\tA}=*^{\tA}*^A=(-1)^{l(n-l)}\det (\imath^{-1}\tA)\,\mathrm{Id}.
\end{equation}
This follows from the computation
\begin{equation*}
\begin{split}
&\det (\imath ^{-1}\tA)\omega \wedge \mu =(*^A\omega )\wedge (*^A\mu )=
(*^A\omega |\mu )_Adx\\
&=(\mu |*^A\omega )_{\tA}dx=\mu \wedge (*^{\tA}*^A\omega
)=(-1)^{(n-l)l}(*^{\tA}*^A\omega )\wedge \mu .
\end{split}
\end{equation*}
}

The second identity in the lemma now follows by applying
  $*^{\tA}$ to the right and to the left in the first one and using
  (\ref{**}).   \hfill{$\Box$}

\medskip\par
We also need a relation between $*^A$ and $\kappa$:

\begin{lemma} \label{propstar}
We have
$
\kappa^* *^{ \tA} = \det \kappa *^A \kappa^*.
$
\end{lemma}

\preuve
{The statement in the lemma is equivalent to the statement that 
$$\kappa ^*\omega \wedge \kappa^**^{\tA}\mu =(\det \kappa )\kappa
^*\omega \wedge*^A\kappa ^*\mu   $$ for all $n-l$-forms $\mu $ and all
$l$ forms $\omega $. Here the left hand side is equal to
$$
\kappa ^*(\omega \wedge *^{\tA}\mu )=\kappa ^*((\omega |\mu
)_{\tA}dx)=(\wedge ^l\tA
\omega \circ \kappa |\mu \circ \kappa )(\det \kappa ) dx.$$}

\par {The right hand side is equal to $(\det \kappa )(\kappa ^*\omega
|\kappa ^*\mu )_A dx$ so we only have to identify the scalar products
in the two expressions:
\begin{equation*}\begin{split}
(\kappa ^*\omega |\kappa ^*\mu )_A=(\wedge^lA\wedge^l\tk
\omega \circ \kappa |\wedge^l\tk \mu \circ \kappa )=(\wedge^l (\kappa
A\tk )\omega \circ \kappa |\mu \circ \kappa ).
\end{split}\end{equation*}
Here $\kappa A\tk=\tA (\tk )^2=\tA$, so the two scalar products are
equal.}  \hfill{$\Box$}

\subsection{Expressions for the quasimodes $f^{(1)}_j$ }\label{expr}

In this subsection we compute the leading amplitude of
the quasimode $f_j^{(1)}$ in (\ref{fk1})  on the manifolds $\pi_x(K_+)$ and 
$\pi_x(K_-^*)$, associated to the saddle point $s_j$.

As a warm up, we shall first show that 
\begin{equation}\label{expr.1}
*^{\tA}(dz_2^*\wedge ...\wedge dz_n^*) =f dy_1
\end{equation}
in a neighborhood of $s_j$, where $f$ is smooth and non-vanishing. In
fact, composing the orthogonality relation $(dy_1^*|dz_k)_A=0$ with
$\kappa $, it is not difficult to see that we also have
\begin{equation}\label{expr.2}
(dy_1|dz_k^*)_{\tA}=0.
\end{equation}
If $\omega $ is an $n-1$ form, we have
$$
(\omega |dz_2^*\wedge ...\wedge dz_n^*)_{\tA }dx=\omega \wedge *^{\tA}(dz_2^*\wedge ...\wedge dz_n^*).
$$
From (\ref{expr.2}) it follows that the left hand side vanishes as
soon as $\omega $ can be written as the exterior product of $dy_1$ and
an $n-2$ form, and then it is easy to see that (\ref{expr.1}) holds. 
$f(s_j)$ will be computed below.

With these notations we have the following result:
\begin{prop} \label{fj}
In a neighborhood of $s_j$ we have
$$
f_j^{(1)}(x) = \sep{ h^{-n/4} a_j(x,h) + \ooo(h^\infty)} e^{-\phi_{j,+}(x)/h},
$$
where we recall that $\phi_{j,+} = (y_1^2 + (z_2^*)^2 + ... + (z_n^*)^2)/2$ and where
$$
a_j(x,h) = a_{j,0}(x) + h a_{j,1}(x) + ...,
$$
\begin{equation} \label{aj}
\begin{split}
  \textrm{ with } \ \ \ & \forall x \in \pi_x(K_+), \ \ \ \ a_{j,0}(x) = \widetilde{\alpha}_{j} dy_1 \\
  \textrm{ and } \ \ \ & \forall x \in \pi_x(K_-^*), \ \ \ \
  a_{j,0}(x) = \widehat{\alpha}_{j} *^{ \tA}(d z_2^* \wedge ... \wedge
  d z_n^*).
\end{split}
\end{equation}
Here $\widehat{\alpha}_{j,0}$ and $\widetilde{\alpha}_{j,0}$ are
non-vanishing constants.
\end{prop}

\preuve { We already know that $a_{j,0}(s_j)$ has to belong to the
kernel of $\frac{1}{2}\widetilde{\mathrm{tr}}F_p+S_P$ at $(s_j,0)$ and
it follows from the discussion at the end of Subsubsection \ref{ssc}
(or from \cite{KFP2}) that $a_{j,0}(s_j)$ is an eigenvector in the
negative eigenspace of $\phi ''(s_j)\tA $. (See also (\ref{ps.4}).)
However, $dy_1(s_j)$ is such a vector since it is orthogonal to
$T_{s_j}(\pi _xK_+)$, the spectral subspace corresponding to the
eigenvalues with real part $>0$ of $A\phi ''(s_j)=(\phi ''(s_j)\tA
)^{\mathrm{t}}$. (Cf (\ref{ss.24}).) Thus we know from the start that
$a_{j,0}(s_j)$ is given by any of the two equivalent expressions in
(\ref{aj}).}

{Mimicking the proof given in \cite{HeSj85}, we use the expression of
the conjugate operator $P_+$ in Lemma \ref{pplus} and have to solve
 \begin{equation} \label{transport}
 P_+ a_{j} = 0,
 \end{equation}
in the sense of formal asymptotic expansions. We get the first 
transport equation $T a_{j,0} = 0$, where
 $$
 T \defegal  \lll_{Ad(\phi + \phi_+)} +  \lll^{A, *}_{ \tA d(\phi - \phi_+)} =   \lll_{Ad g_-^*} -  \lll^{A, *}_{ \tA dg_+}
$$
where we used the definitions of $g_-^*$ and $g_+$ in (\ref{sstar}),
(\ref{defg}). The transport equation for $a_{j,0}$ reads
 \begin{equation} \label{transp}
 \sep{ \lll_{Ad g_-^*} -  \lll^{A, *}_{ \tA dg_+} } a_{j,0} = 0.
 \end{equation}}

\par {In order to prove the first identity in(\ref{aj}) we try to solve
(\ref{transp}) along $\pi _x(K_+)$ with 
\ekv{expr.3}
{
a_{j,0}(x)=\widetilde{\alpha }_j(z_2,...,z_n)dy_1+{\cal O}(y_1).
}
For general reasons, we already know that the vector field part of the
transport operator is tangent to $\pi _x(K_+)$. Using (\ref{liestar}),
we get 
$$
{\cal L}^{A,*}_{\tA g_+}a_{j,0} = (dg_+)^\wedge
d^{A,*}a_{j,0}+d^{A,*}(dg_+)^\wedge a_{j,0}. 
$$
Here $dg_+=0$ on $\pi _x(K_+)$ so the first term vanishes there. Using
the form of $a_{j,0}$ above and the fact that $dg_+=2y_1dy_1$, we get
$(dg_+)^\wedge a_{j,0}={\cal O}(y_1^2)$ and the second term also
vanishes on $\pi _x(K_+)$.}

{Thus we only have to solve along $\pi _x(K_+)$ the equation
\begin{equation} \label{transreduce}
 \lll_{Ad  g_-^*} a_{j,0}(x) = 0,
 \end{equation}
 still with $a_{j,0}$ of the form (\ref{expr.3}) and we have to check
 that $\widetilde{\alpha_j}$ is in fact constant on $\pi_x(K_+)$.
 From Cartan's formula for Lie derivatives, (\ref{transreduce}) reads
\begin{equation}
\begin{split}
  0 & =  d \sep{ \widetilde{\alpha}_j(z) (A dg_-^*)^\rfloor d y_1} +
  \sep{ (A d g_-^*)^\rfloor d\widetilde{\alpha}_j} d y_1
  -d\widetilde{\alpha }_j\underbrace{\langle Adg_-^*|dy_1\rangle }_{=0} \\
 & = d ( \widetilde{\alpha}_j \underbrace{( dg_-^*| dy_1)_A}_{=0} ) + 
( dg_-^*| d \widetilde{\alpha}_j )_A d y_1 \\
 & = ( dg_-^*| d \widetilde{\alpha}_j)_A d y_1  = 2 ( d \phi| d \widetilde{\alpha}_j)_A d y_1,
 \end{split}
 \end{equation}
where we used that $d g_-^* = 2 d \phi$ on $\pi_x(K_+)$, and also the eikonal equation. Recalling that $2 A d\phi = \nu_+ \in T (\pi_x(K_+))$ we get  that on $\pi_x(K_+)$
$$
(\ref{transreduce}) \Longleftrightarrow \nabla_{A d\phi} \widetilde{\alpha}_j = 0 \Longleftrightarrow \widetilde{\alpha}_j = \mathrm{Cte},
$$
from the standard properties of non-degenerate vector fields. This proves the first assertion in (\ref{aj}).}

\par
{For the second one the tools are essentially the same.
On $\pi_x(K_-^*) $ we look for a solution of the type
$$
a_{j,0}(x) =  \widehat{\alpha_j}(y_1^*) *^{ \tA} ( d z_2^*\wedge
... \wedge d z_n^*)+{\cal O}(z^*).
$$
and we have to  check that $\widehat{\alpha_j}$ is in fact constant on
$\pi_x(K_-^*)$.}

\par {Let us check that the first term in
  (\ref{transp}) vanishes on $\pi _x(K_-^*)$. $dg_-^*$ vanishes on
  $\pi _x(K_-^*)$, so
$$
{\cal L}_{Adg_-^*}a_{j,0}=d\circ (Adg_-^*)^{\rfloor}a_{j,0}\hbox{
  on }\pi _x(K_-^*),
$$
and it suffices to check that 
$$
d\left((Adg_-^*)^{\rfloor}*^{\tA }(dz_2^*\wedge ...\wedge
  dz_n^*)\right) =0,
$$
i.e. that 
$$
0=d((Adg_-^*)^{\rfloor} fdy_1)=d(f(dg_-^*|dy_1)_A),
$$
and this is zero, since $(dg_-^*|dy_1)_A=0$ as we have already
observed and used.}

{Thus the transport equation on $\pi _*(K_-^*)$ becomes
\begin{equation}
 \lll^{A, *}_{ \tA dg_+} *^{ \tA} \sep{  \widehat{\alpha_j}(y_1^*)
  d z_2^*\wedge ... \wedge d z_n^* } = 0.
 \end{equation}
Lemma \ref{starlie} shows that  the preceding equality is equivalent to
\begin{equation}
 \lll_{ \tA dg_+}  \widehat{\alpha_j}(y_1^*)
  ( d z_2^*\wedge ... \wedge d z_n^*) = 0.
 \end{equation}
Exactly as in the preceding case and using again the eikonal equation,  it reduces on $\pi_x(K_-^*)$  to
 $$
  2 \seq{ d \phi, d \widehat{\alpha}_j }_{ \tA } d z_2^* \wedge ... \wedge d z_n^* = 0
  $$
  and again this is satisfied  if $\widehat{\alpha}_j = \mathrm{Cte}$: Indeed it also reads
  $$
  \nabla_{ \tA d\phi} (\widehat{\alpha}_j (y_1^*))  = 0,
  $$
  and  $2 \tA d\phi = -\nu_- \in T (\pi_x(K_-^*))$ from (\ref{bienutile}). Again we used the standard properties of non-degenerate vector fields. The proof of Proposition \ref{fj} is complete.} \hfill{$\Box$}

\medskip\par
  Now we  evaluate the coefficients $\widetilde{\alpha}_j$ and $\widehat{\alpha}_j$. For this we simply write that $f_j^{(1)}$ is $(\kappa, A)$-normalized:
  \begin{equation*}
  \begin{split}
  1  & = \norm{f_j^{(1)}}^2_{\kappa, A} = \sep{ \kappa^* f_j^{(1)}| f_j^{(1)}}_A \\
     & = h^{-n/2} \int \theta_j \circ \kappa (x) \theta_j (x) \seq{\kappa^* a_j (x)| a_j (x)}_A e^{- \phi_+ \circ \kappa (x)/h} e^{-\phi_+ (x)/h} dx
     \end{split}
\end{equation*}
where $\theta_j$ is the truncation function introduced in the
preceding section. 

\par Using (\ref{expphiplus}) we check that
\begin{equation} \label{phas}
\begin{split}
\phi_+ \circ \kappa +\phi_+ & =
\frac{1}{2} \sep{  y_1^2 + (z_2^*)^2 + ... + (z_n^*)^2 } + \frac{1}{2} \sep{  (y_1^*)^2 + z_2^2 + ... + z_n^2 } \\
 & = \frac{1}{2} \sep{  y_1^2 + z_2^2 + ... + z_n^2 }  + \frac{1}{2} \sep{  (y_1^*)^2 + (z_2^*)^2 + ... + (z_n^*)^2 }  \\
 & = (y_1^*)^2 + (z_2^*)^2 + ... + (z_n^*)^2\\
& =y_1^2+z_2^2+...+z_n^2,
 \end{split}
 \end{equation}
 where we used (\ref{sumsim}) for the last two equalities.
 
As for the amplitude in the integral, we use that at $s_j$ we have two expressions for $a_j$. To leading order w.r.t. $h$, we have:
 \begin{equation} \label{ampl}
 \begin{split}
 \seq{\kappa^*a_{j, 0}| a_{j,0}}_A dx
 & =  \widetilde{\alpha}_{j} \widehat{\alpha}_{j}
(  dy_1^*|   *^ {\tA} (d z_2^* \wedge ... \wedge d z_n^*) )_A dx \\
 &=  \widetilde{\alpha}_{j} \widehat{\alpha}_{j}
 (  *^\tA (d z_2^* \wedge ... \wedge d z_n^*)| dy_1^*)_{\tA} dx \\
 & = \widetilde{\alpha}_{j} \widehat{\alpha}_{j}   *^\tA (d z_2^* \wedge ... \wedge d z_n^*) \wedge  *^\tA d y_1^* \\
 & =\pm  (\det \imath^{-1}A) \widetilde{\alpha}_{j} \widehat{\alpha}_{j} \sep{  d y_1^* \wedge d z_2^* \wedge ... \wedge d z_n^* } \\
  \end{split}
  \end{equation}
  where we used Lemma \ref{defs} {and stopped trying to
    follow up the signs} for simplicity. Putting (\ref{phas}) and (\ref{ampl}) in the
  integral, using the change of variable $\kappa$ and applying the Laplace
  method, we get
   \begin{equation*}
  \begin{split}
  1  & = \norm{f_j^{(1)}}^2_{\kappa, A} \\
     & = \pm h^{-n/2}  (\det \imath^{-1} A)
     \widetilde{\alpha}_{j} \widehat{\alpha}_{j}  \int \theta_j \circ
     \kappa \, \theta_j\  e^{-(y_1^2 + z_2^2 + ... + z_n^2)/h} dy_1
     dz_2 ... dz_n{+{\cal O}(h)} \\
    & = \pm \pi^{n/2} (\det \imath^{-1}A ) \widetilde{\alpha}_{j} \widehat{\alpha}_{j} + \ooo( h).
     \end{split}
\end{equation*}
so that
\begin{equation} \label{qsi}
 \pm 1  = \pi^{n/2}  (\det \imath^{-1} A ) \widetilde{\alpha}_{j} \widehat{\alpha}_{j}.
\end{equation}

\par We shall next compute $f(s_j)$ in (\ref{expr.1}). Equivalently,
we shall compute $g=g(s_j)$ in the relation
\ekv{cal.1}
{
*^Ady_1=gdz_2^*\wedge ...\wedge dz_n^*.
}
Indeed, if we apply $*^{\tA}$ to this and use (\ref{**}), we get 
\ekv{cal.2}
{
\det (\imath^{-1}\tA)dy_1=g*^{\tA}dz_2^*\wedge ...\wedge dz_n^*,
}
which is (\ref{expr.1}) with 
\ekv{cal.3}
{
f=\pm \frac{\det (\imath^{-1}\tA)}{g}.
}

\par (\ref{cal.1}) means that
$$
(\omega |dy_1)_Adx_1\wedge...\wedge dx_n=\pm g\omega \wedge
dz_2^*\wedge ...\wedge dz_n^*,
$$
for all 1 forms $\omega $ (and with a sign independent of $\omega
$). With $\omega =dy_1^*$, we get
\ekv{cal.4}
{
(dy_1^*|dy_1)_Adx_1\wedge ...\wedge dx_n=\pm gdy_1^*\wedge
dz_2^*\wedge ...\wedge dz_n^*.
}

\par Now use (\ref{expphiplus})
\ekv{cal.4.5}
{\phi -\phi (s_j)=\frac{1}{2}(-y_1^2+(z_2^*)^2+...+(z_n^*)^2)=\frac{1}{2}(-(y_1^*)^2+z_2^2+...+z_n^2),}
where the last equality follows from $\phi \circ \kappa =\phi $. On
$\pi _x(K_-)$ we have $z_2=...=z_n=0$ and hence
$$
\nu _+=2Ad\phi =-2A(y_1^*dy_1^*)=-2y_1^*Ady_1^*.
$$
Since $\nu _+$ is tangent to $\pi _x(K_-)$, we conclude that 
\ekv{cal.5}
{
Ady_1^*=b\frac{\partial }{\partial y_1},
}
so 
\ekv {cal.6}{\nu _+=-2by_1^*\frac{\partial }{\partial y_1}\hbox{ on
  }\pi _x(K_-).}

\par For notational reasons, we sometimes write $z_1$, $z_1^*$ instead
of $y_1$, $y_1^*$. Recall that 
$$
z_j^*(x)=z_j(\kappa (x))=\sum_{k=1}^\nu \kappa _{j,k}z_k,
$$
where $(\kappa _{j,k})$ is the matrix of $\kappa $ with respect to the
coordinates $z_1,...,z_n$. On $\pi _x(K_-)$ we get $y_1^*=\kappa
_{1,1}y_1$ and (\ref{cal.6}) becomes
\ekv{cal.7}
{
\nu _+=-2b\kappa _{1,1}y_1\frac{\partial }{\partial y_1}\hbox{ on }\pi _x(K_-).
}
On the other hand, we know that $\pi _x(K_-)$ is an eigenspace of the
linearization of $\nu _+$ with associated eigenvalue
$-2\widehat{\lambda }_1<0$, {writing $\widehat{\lambda }_1=-\lambda _1$
where $\lambda _1$ is the negative eigenvalue in (\ref{ss.24}),}
and a comparison with (\ref{cal.7}) shows that $-2b\kappa
_{1,1}=-2\widehat{\lambda }_1$, so
\ekv{cal.8}
{
b=\frac{\widehat{\lambda }_1}{\kappa _{1,1}}.
}

\par Combining this with (\ref{cal.5}), we get 
\ekv{cal.9}
{
(dy_1^*|dy_1)_A=\langle Ady_1^*|dy_1\rangle =b\underbrace{\langle
  \frac{\partial }{\partial y_1}|dy_1\rangle}_{=1}=\frac{\lambda
  _1}{\kappa _{1,1}}.
}
In Remark \ref{cal1} we will give practically computable formulas for
$\kappa _{1,1}$. Also recall that the quantity (\ref{cal.9}) is $>0$
since we assume that we are in case i) of Proposition \ref{ps1}. Hence
$\kappa _{1,1}>0$. From the last equality in (\ref{phas}) we also know
that the matrix $(\kappa _{j,k})$ is orthogonal and in particular that
$\kappa _{1,1}\le 1$:
\ekv{cal.11}
{
0<\kappa _{1,1}\le 1.
}

\par Inserting (\ref{cal.9}) in (\ref{cal.4}), we get
\ekv{cal.12}
{
g=\pm \frac{\widehat{\lambda }_1}{\kappa _{1,1}}\frac{1}{\det \frac{\partial
    z^*}{\partial x}},
}
where we recall that we write $z_1^*$ for $y_1^*$ whenever
convenient. From (\ref{phas}), we get
$$
{\frac{1}{2}} (\phi _++\phi _+\circ \kappa )''_{xx}=\left(\frac{\partial
    z^*}{\partial x}\right)^{\mathrm{t}}\frac{\partial z^*}{\partial x},
$$
so 
$$
\det \frac{\partial z^*}{\partial x}=\pm (\det {\frac{1}{2}} (\phi _++\phi _+\circ
\kappa )''_{xx})^{1/2},
$$
and (\ref{cal.12}) gives
\ekv{cal.13}{g=\pm \frac{\widehat{\lambda }_1}{\kappa _{1,1}(\det {\frac{1}{2}} (\phi _++\phi
    _+\circ \kappa )''_{xx})^{1/2}}.} From (\ref{cal.3}) we now get at
$s_j$:
\ekv{cal.14}
{
f=\pm \frac{\kappa _{1,1}}{\widehat{\lambda }_1}(\det {\frac{1}{2}} (\phi _++\phi _+\circ
\kappa )''_{xx})^{1/2}\det (\imath^{-1}\tA).
}

\begin{remark}\label{cal1}
Here is a direct way of getting $\kappa _{1,1}$: Recall that $\kappa
_{1,1}$ is equal to $y_1^*/y_1$ on $\pi _x(K_-)$, so
$$
\kappa _{1,1}^2=\frac{(y_1^*)^2}{y_1^2}\hbox{ on }\pi _x(K_-).
$$
Using (\ref{cal.4.5}), we get $(y_1^*)^2=-2(\phi -\phi (s_j)) $ on $\pi
_x(K_-)$. Similarly, by (\ref{phas}), we have $y_1^2=\phi _++\phi_+
\circ \kappa $ on $\pi _x(K_-)$, so
\ekv{cal.14.5}
{
\kappa _{1,1}^2=\frac{-(\phi -\phi (s_j))}{{\frac{1}{2}}(\phi _++\phi _+\circ
  \kappa )}\hbox{ on }\pi _x(K_-).
}
It follows from the discussion that the right hand side of
(\ref{cal.14.5}) is constant on $\pi _x(K_-)$. Thus, in order to
compute $\kappa _{1,1}$ up to the sign, it suffices to compute the
eigendirection $\pi _x(K_-)$ and the {Hessian of the} positive definite solution $\phi
_+$ of the Hamilton-Jacobi equation.
\end{remark}

\medskip
\par Combining (\ref{expr.1}) and (\ref{aj}), we get 
\ekv{cal.15}
{
\widehat{\alpha }_jf(s_j)=\widetilde{\alpha }_j,
}
where $f$ is given above, and using this and (\ref{cal.14}) in
(\ref{qsi}), we get
\ekv{cal.16}
{
1=\pm \pi ^{n/2}\det (\imath^{-1}A)\det (\imath^{-1}\tA)\frac{\kappa
  _{1,1}}{\widehat{\lambda }_1}(\det {\frac{1}{2}} (\phi _++\phi _+\circ \kappa
)''_{xx})^{1/2}\widehat{\alpha }_j^2.
}
Here we notice that $\imath : (\mathbb{R}^n)^*\to \mathbb{R}^n$ is
symmetric by its definition: $\langle \imath \omega |\mu \rangle
=(\omega |\mu )$, $\omega ,\mu \in (\mathbb{R}^n)^*$, so 
$$
\det \imath^{-1}\tA =\det (\imath^{-1}\tA)^{\mathrm{t}}=\det
A(\imath^{-1})^{\mathrm{t}}=\det A\imath^{-1}=\det \imath^{-1}A.
$$
Thus we get from (\ref{cal.16})
\ekv{cal.17}
{
\pm 1=\pi ^{n/2}(\det\imath^{-1}A)^2(\det {\frac{1}{2}} (\phi _++\phi _+\circ
\kappa )''_{xx})^{1/2}\frac{\kappa _{1,1}}{\widehat{\lambda }_1}\widehat{\alpha }_j^2,
}
where the sign to the left is the one that allows $\widehat{\alpha }_j$
to be real. Assuming that $\kappa _{1,1}>0$ (which can be arranged by
the choice of sign of $y_1^*$), we get,
\ekv{cal.18}
{
\widehat{\alpha }_j=\pm \pi ^{-\frac{n}{4}}(\det
\imath^{-1}A)^{-1}(\det {\frac{1}{2}} (\phi _++\phi _+\circ \kappa
)''_{xx})^{-\frac{1}{4}}
\left(\frac{\widehat{\lambda }_1}{\kappa _{1,1}}\right)^{\frac{1}{2}},
} 
where we are free to choose the sign (implying a choice of sign for
$\widetilde{\alpha }_j$).

\subsection{Computation of the coefficient}

Now we can compute the coefficient  of the singular matrix
defined in (\ref{step2}) according to (\ref{step1}) in the preceding section.
First recall the expressions for the quasimodes $f^{(1)}_j$ and $f_k^{(0)}$ in a neighborhood of the saddle point $s_j$ :

As for $f^{(1)}_j$ we have the following definition from Proposition \ref{fj}
 in a neighborhood of $s_j$
$$
f_j^{(1)}(x) = \sep{ h^{-n/4} a_j(x,h) + \ooo(h^\infty)} \theta_j (x) e^{-\phi_{j,+}(x)/h},
$$
where we recall that $\phi_{j,+} = (y_1^2 + (z_2^*)^2 + ... + (z_n^*)^2)/2$.
and that
\begin{equation} \label{ajbis}
\begin{split}
 \forall x \in \pi_x(K_-^*), \ \ \ \ a_{j,0}(x) = \widehat{\alpha}_{j} *^\tA(d z_2^* \wedge ... \wedge d z_n^*).
\end{split}
\end{equation}
Now we deal with the quasimode $f_k^{(0)}$ where $k = k(j)$  according to the injection defined in the preceding section. We have
\begin{equation*}
f_k^{(0)}= h^{-n/4} c_k(h)  e^{-{1\over h}(\phi (x)-\phi (m_k))}\chi _{k,\eps} (x)
\end{equation*}
where $\chi_{k,\eps}$ is the cutoff function defined in the generic
case and $c_k(h)>0$ a normalization constant.
In the following it is convenient to pose
$$
X_k = \sqrt{\det(\phi''(m_k))}
$$
and a direct application of the Laplace method gives
\begin{equation} \label{defX}
c_k = (1+{\cal O}(h))\frac{\sqrt{X_k}}{\pi^{n/4}},
\end{equation}
so that
\begin{equation*} \label{fk0ter}
f_k^{(0)}=   \frac{\sqrt{X_k}}{(\pi h)^{n/4}}(1+{\cal O}(h))  e^{{-}(\phi (x)-\phi (m_k))/h}\chi _{k,\eps} (x).
\end{equation*}
We now choose
 the behavior of the cut-off $\chi_{k,\eps}$ function near $s_j$. In addition to all the properties recalled in Subsection \ref{labcut}, we impose that in a small
neighborhood $\vvv$ of $\pi_x(K_-)  \cap \set{ y_1 >0}$,
$$
\chi_{k,\eps}(x) = \widetilde{\chi}_{k, \eps}(y_1) \ \ \textrm{ and } \ \ \
\textrm{ supp } \chi_{k,\eps} \cap \set{ y_1 \leq 0} = \emptyset.
$$
As in Subsection \ref{labcut} we impose that the function $\theta_j$ is equal to $1$ in a far larger ball of radius $\simeq \eps_1 \gg \eps$.

{
\begin{figure}[htbp]
\begin{center}
\begin{picture}(0,0)%
\includegraphics{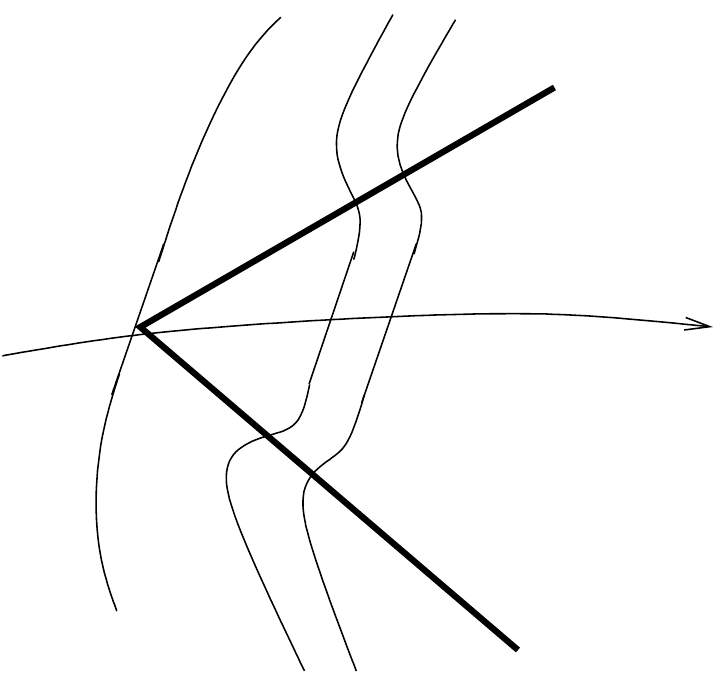}%
\end{picture}%
%
%
\setlength{\unitlength}{3947sp}%
\begingroup\makeatletter\ifx\SetFigFont\undefined%
\gdef\SetFigFont#1#2#3#4#5{%
  \reset@font\fontsize{#1}{#2pt}%
  \fontfamily{#3}\fontseries{#4}\fontshape{#5}%
  \selectfont}%
\fi\endgroup%
\begin{picture}(3439,3232)(677,-3549)
\put(3339,-2111){\makebox(0,0)[lb]{\smash{{\SetFigFont{12}{14.4}{\rmdefault}{\mddefault}{\updefault}{\color[rgb]{0,0,0}$\pi_x(K_-)$}%
}}}}
\put(1289,-474){\makebox(0,0)[lb]{\smash{{\SetFigFont{12}{14.4}{\rmdefault}{\mddefault}{\updefault}{\color[rgb]{0,0,0}$\pi_x(K_+)$}%
}}}}
\put(4101,-1749){\makebox(0,0)[lb]{\smash{{\SetFigFont{12}{14.4}{\rmdefault}{\mddefault}{\updefault}{\color[rgb]{0,0,0}$y_1$}%
}}}}
\put(1089,-1811){\makebox(0,0)[lb]{\smash{{\SetFigFont{12}{14.4}{\rmdefault}{\mddefault}{\updefault}{\color[rgb]{0,0,0}$s_j$}%
}}}}
\put(3439,-711){\makebox(0,0)[lb]{\smash{{\SetFigFont{12}{14.4}{\rmdefault}{\mddefault}{\updefault}{\color[rgb]{0,0,0}level line of  $\phi$}%
}}}}
\put(2726,-1724){\makebox(0,0)[lb]{\smash{{\SetFigFont{12}{14.4}{\rmdefault}{\mddefault}{\updefault}{\color[rgb]{0,0,0}\tiny $\chi_{\epsilon,k} = 1$}%
}}}}
\put(1739,-1774){\makebox(0,0)[lb]{\smash{{\SetFigFont{12}{14.4}{\rmdefault}{\mddefault}{\updefault}{\color[rgb]{0,0,0}\tiny $\chi_{\epsilon,k} = 0$}%
}}}}
\end{picture}%
\end{center}
\end{figure}

\fontfamily{ptm}

We can  compute the tunneling coefficient defined in (\ref{step1}):
\begin{multline} \label{step1ter}
( f_{j}^{(1)} | d_\phi f_k^{(0)} )_{A, \kappa} \\
=   h^{1-n/2} c_k(h)
 \int 
\theta_j  \circ \kappa(x) \sep{ \sep{ \kappa ^*a_{j}}(x,h) | d \chi_{k,\eps}(x)}_A
e^{-(\phi _{+,j} \circ \kappa (x) + \phi(x) - \phi(m_k))/h} dx \\
 + \ooo ( h^{-N} e^{ -( S_k +  C_\eps^{-1})/h})
\end{multline}
We use the adapted coordinates. As already noticed the phase is up to a division by $h$ equal to
\begin{equation*} 
\begin{split}
  -(\phi _{+,j} \circ \kappa (x) + \phi(x) - \phi(m_k))
& =  -(\phi _{+,j} \circ \kappa (x) + \phi(x) - \phi(s_j)) - S_k \\
& = -g_-(x) -S_k \\
& =  -(z_2^2 + ... + z_n^2) - S_k.
\end{split}
\end{equation*}
We therefore get
\begin{multline} \label{intdx}
( f_{j}^{(1)} | d_\phi f_k^{(0)} )_{A, \kappa} \\
=   h^{1-n/2} e^{-S_k/h} c_k(h)
\bigg( \int
\theta_j  \circ \kappa(x) \sep{ \sep{ \kappa ^*a_{j}}(x,h) | d \chi_{k,\eps}(x)}_A e^{-(z_2^2 + ... + z_n^2)/h} dx \bigg) \\
 + \ooo ( h^{-N} e^{ -( S_k +  C_\eps^{-1})/h})
\end{multline}
where we recall that $j = j(k)$.
On $\pi_x(K_-^*)$ we have by Proposition \ref{fj}
$$
a_j(x,h) =  \widehat{\alpha}_j  *^{ \tA} (d z_2^* \wedge ... \wedge d z_n^*) + \ooo(h).
$$
Since  $z_l^* = z_l \circ \kappa $ we can use Lemma \ref{propstar} to get on  $\pi_x(K_-)$
\begin{equation*} 
\begin{split}
 \forall x \in \pi_x(K_-), \ \ \ \ a_{j,0}^*(x) \defegal \kappa^* a_{j,0} = (\det \kappa) \widehat{\alpha}_{j} *^A(d z_2 \wedge ... \wedge d z_n).
\end{split}
\end{equation*}
From the expression of $\chi_{k,\eps}$ near $s_j$ we also have in a neighborhood of
$\pi_x(K_-)$,
$$
d \chi_{k,\eps}(x) = \widetilde{\chi}_{k,\eps}'(y_1) dy_1,
$$
 and putting these two expressions together with iii) of  Lemma \ref{defs}, we get on $\pi_x(K_-)$
\begin{multline} \label{AA}
 \sep{ ( \kappa ^*a_{j})(x,h) , d \chi_{k,\eps}(x)}_A dx \\
 =   \pm \widehat{\alpha}_j \widetilde{\chi}_{k,\eps}'(y_1) \sep{   *^{ A} ( d z_2 \wedge ... \wedge d z_n)| d y_1}_A dx+ \ooo(h) \\
  = \pm  \widehat{\alpha}_j \widetilde{\chi}_{k,\eps}'(y_1)    *^{ A} ( d z_2 \wedge ... \wedge d z_n) \wedge *^A d y_1 + \ooo(h) \\
 =  \pm \widehat{\alpha}_j (\det \imath^{-1} A) \widetilde{\chi}_{k,\eps}'(y_1) dy_1 \wedge dz_2 \wedge ... \wedge dz_n + \ooo(h)
\end{multline}
Using the Laplace method in the $z$ coordinates we therefore get from (\ref{intdx}) and (\ref{AA})
\begin{equation} \label{intdx2}
\begin{split}
& ( f_{j}^{(1)} | d_\phi f_k^{(0)} )_{A, \kappa} \\
 & \equiv \pm  (\det \imath^{-1} A)   h^{1-n/2}   e^{-S_k/h} c_k(h)  (\pi h)^{(n-1)/2}  \widehat{\alpha}_j
 \int_{\pi_x{K_-} \cap \vvv} \widetilde{\chi}_{k,\eps}'(y_1) dy_1
\\
& \equiv \pm (\det \imath^{-1} A) \pi^{(n-1)/2} h^{1/2} e^{-S_k/h} c_k(h) \widehat{\alpha}_j.
\end{split}
\end{equation}
modulo $ \ooo(h^{3/2}e^{-S_k/h})$.
Combining this with (\ref{defX}), (\ref{cal.18}), we finally get

\begin{prop} \label{endestim}
Let $s_j$ be a saddle point corresponding to a minimum $m_k$, with $j=j(k)$. Then we have
$$
( f_{j}^{(1)} | d_\phi f_k^{(0)} )_{A, \kappa} = 
\pm \left(\frac{h}{\pi
  }\right)^{\frac{1}{2}}\left(\frac{\widehat{\lambda }
    _1}{\kappa _{1,1}}\right)^{\frac{1}{2}}\frac{(\det \phi
  ''(m_k))^{\frac{1}{4}}}
{(\det {\frac{1}{2}} (\phi _++\phi _+\circ \kappa )''(s_j))^{\frac{1}{4}}}e^{-S_k/h}
+ \ooo( h^{3/2}e^{-S_k/h}),
$$
where we recall that $-\widehat{\lambda }_1$ is the negative eigenvalue of
$A\phi ''(s_j)$ and that $\kappa _{1,1}$ is the restriction of
$$
\left( \frac{-2(\phi -\phi (s_j))}{\phi _++\phi _+\circ \kappa }\right)^{\frac{1}{2}}
$$ to $\pi _x(K_-)$ (whose tangent space at $s_j$ is the corresponding eigenspace).
As a consequence the exponentially small eigenvalues of $P = -\Delta_A^{(0)}$ are real and asymptotically given by
$$
\mu_k = h \left( \frac{1}{\pi }\frac{\widehat{\lambda }
    _1}{\kappa _{1,1}}\frac{(\det \phi
  ''(m_k))^{\frac{1}{2}}}
{(\det {\frac{1}{2}} (\phi _++\phi _+\circ \kappa )''(s_j))^{\frac{1}{2}}}  + \ooo(h)\right) e^{-2 S_k/h}, \ \ \ \  k =1, ..., n_0
$$
where $S_1 = \infty$ (and $\mu_1=0$) by convention.
\end{prop}

To finish this section we shall make Proposition \ref{endestim} even
more explicit in the case of the Kramers-Fokker-Planck operator
(\ref{in.0}). Assume for simplicity that the saddle point $s_j$ is
placed at the origin $x=0$, $y=0$. Our calculations only concern the
quadratic approximation, so we may assume most of the time that $V$ is
quadratic. After an orthogonal change of variables, we may assume that 
\ekv{kf.1}
{
V(x)=-\frac{v}{2}x_1^2+\sum_2^d \frac{v_j}{2}x_j^2,
}
where $v,v_j$ are $>0$. The variables can be completely separated and
only the $x_1,y_1$ variables are of interest, so in most of the
following discussion we will assume that $d=1$ with some remarks about
the formulation of the corresponding results when $d>1$. Thus we
consider the one dimensional potential,
\ekv{kf.2}
{
V(x)=-\frac{v}{2}x^2
}
where $v>0$. Recall that
\ekv{kf.3}{
\begin{split}
p&=i(y\xi +vx\eta) +\frac{\gamma }{2}(y^2+\eta ^2),\\
q&=-p(x,y;i\xi ,i\eta )=y\xi +vx\eta+\frac{\gamma }{2}(\eta ^2-y^2),\\
\phi (x,y)&=\frac{1}{2}(y^2-vx^2). 
\end{split}
}

\par Let $\phi _+(x,y)=\frac{1}{2}(ax^2+2bxy+cy^2)$ be the unique positive
definite quadratic form which solves the eiconal equation
\ekv{kf.4}
{
q(x,y;\partial _x\phi _+,\partial _y\phi _+)=0.
}
Expanding the left hand side as a quadratic form and equating the
coefficients to zero, we get
\ekv{kf.5}
{
\frac{\gamma }{2}b^2+vb=0,\ a+vc+\gamma bc=0,\ b+\frac{\gamma
}{2}c^2=\frac{\gamma }{2}.
}

\par Choosing $b=0$ as the solution of the first equation, leads to the
two solutions $\pm \phi (x,y)$, neither of which is positive
definite. Thus we have to choose the other solution, $b=-2v/\gamma
$. Then the last equation in (\ref{kf.5}) leads to 
$c=\pm \sqrt{1+\frac{4v}{\gamma ^2}}$ and the condition that $\phi _+$
is positive definite imposes the choice of the plus sign. Finally we
determine $a$ from the middle equation and get
$a=v\sqrt{1+\frac{4v}{\gamma ^2}}$. Thus the unique solution to our
problem is 
\ekv{kf.6}
{
\phi _+=\frac{v}{2}\sqrt{1+\frac{4v}{\gamma ^2}}x^2-\frac{2v}{\gamma
}xy+\frac{1}{2}\sqrt{1+\frac{4v}{\gamma ^2}}y^2,
}
which can be seen directly to be positive definite. (In higher
dimesions, we get 
$$
\phi _+^{(d)}(x,y)=\phi _+^{(1)}(x_1,y_1)+\sum_2^d \frac{1}{2}(y_j^2+v_jx_j^2),
$$ 
differing from the expression for $\phi ^{(d)}(x,y)$ only in the
variables $x_1,y_1$, where we used the superscripts 1 and $d$ to
indicate the functions $\phi _+$ and $\phi $ in dimension 1 and $d$
respectively. In this case, $\kappa :\, (x,y)\mapsto (x,-y)$, so 
\ekv{kf.7}
{
\phi _++\phi _+\circ \kappa =v\sqrt{1+\frac{4v}{\gamma ^2}}x^2+\sqrt{1+\frac{4v}{\gamma ^2}}y^2.
}

\par Next we compute the negative eigenvalue $-2\widehat{\lambda }_1$ of $2A\phi
''(s_j)$ and the corresponding eigenspace $\pi _x(K_-)$. Since 
$$
2A\phi ''(s_j)=\left(\begin{matrix} 0 &1\\ -1 &
    \gamma \end{matrix}\right)\left(\begin{matrix}-v &0\\ 0
    &1\end{matrix}\right)=\left(\begin{matrix}0 &1\\ v &\gamma \end{matrix}\right),
$$
the eigenvalues are $\frac{\gamma }{2}\pm \sqrt{(\gamma /2)^2+v}$, so 
\ekv{kf.8}
{
\widehat{\lambda }_1=\frac{1}{2}(\sqrt{(\gamma /2)^2+v}-\gamma /2). 
}
the corresponding eigenspace $\pi _x(K_-)$ is generated by the vector
$(1,\sqrt{(\gamma /2)^2+v}-\frac{\gamma }{2})$ (In $d$ dimensions,
$\pi _x(K_-)$ is generated by a vector $(x,y)$ with $(x_1,y_1)$ equal
to the vector above and with the other components equal to 0.) 

\par $\kappa _{1,1}^2$ is the restriction of $(-2\phi )/(\phi _++\phi
_+\circ \kappa )$ to $\pi _x(K_-)$ (where by our normalization, we have
$\phi (s_j)=0$, $s_j=(0,0)$), so after some straight forward
calculations, we get
\ekv{kf.9}
{
\kappa _{1,1}=\frac{\gamma }{\sqrt{\gamma ^2+4v}}.
} 
The formulae (\ref{kf.8}) and (\ref{kf.9}) remain valid in $d$
dimensions, with $-v$ denoting the negative eigenvalue of
$V''(s_j)=V''(0)$.

\par In the first formula in Proposition \ref{endestim}, we write
\ekv{kf.10} { \frac{\det \phi ''(m_k)}{\det
      \frac{1}{2}(\phi _++\phi _+\circ \kappa )''(s_j) }
    =\frac{\det \phi ''(m_k)}{-\det \phi ''(s_j)}\times\frac {-\det \phi
      ''(s_j)}{\det \frac{1}{2}(\phi _++\phi _+\circ \kappa )''(s_j)}.
  } With $V$ as in (\ref{kf.1}), the last factor reduces to
the one-dimensional case and we get \ekv{kf.11} { \frac {-\det \phi
    ''(s_j)}{\det \frac{1}{2}(\phi _++\phi _+\circ \kappa )''(s_j)}=
  \frac{-\det \left(\begin{matrix}-v &0\\ 0
        &1\end{matrix}\right)}{\det
    \left(\begin{matrix}v\sqrt{1+\frac{4v}{\gamma ^2}}&0\\
        0&\sqrt{1+\frac{4v}{\gamma
            ^2}}\end{matrix}\right)}=\frac{1}{1+\frac{4v}{\gamma ^2}}.
}

\par Combining (\ref{kf.8})--(\ref{kf.11}) with Proposition
\ref{endestim}, we get after some straight forward reductions,
\ekv{kf.12}
{
(f_j^{(1)}|d_\phi f_k^{(0)})=\pm \left(\frac{h}{\pi
  }\right)^{\frac{1}{2}}
\left(\frac{\det V ''({\bf m}_k)}{-\det V
    ''({\bf s}_j)}\right)^{\frac{1}{4}}
\widehat{\lambda }_1^{\frac{1}{2}}e^{-\frac{S_k}{h}}+{\cal O}(h^{\frac{3}{2}}) e^{-\frac{S_k}{h}
},}
so that the exponentially small eigenvalues are given by
$$
\mu _k=\frac{h}{\pi }\left(\frac{\det V ''({\bf m}_k)}{-\det V
    ''({\bf s}_j)}\right)^{\frac{1}{2}}\widehat{\lambda
}_1e^{-\frac{S_k}{h}}+{\cal O}(h^2)e^{-\frac{S_k}{h}}.
$$

\section{Multiple well analysis in the general case}\label{gl}
\setcounter{equation}{0}
In this section we return to the general case and show how the analysis of Sections \ref{lbl} and \ref{genc} may be combined to obtain the other main result of this work,
Theorem \ref{gl4}. {In Theorem \ref{fa1} we also give the
full asymptotic expansion for the smallest non-vanishing eigenvalue
under a generic assumptions which is weaker than the one in Theorem
\ref{theogene}}. The section is concluded by a brief discussion of an explicit example of a system with three minima and three saddle points,
where the values of $\phi$ at the different minima and the saddle points may coincide, illustrating the transition from a 
degenerate case to a generic one. 

\subsection{Statement of the main result and the matrix of $d_{\phi}: E^{(0)}\rightarrow E^{(1)}$} \label{sl}
In Section \ref{lbl} we have constructed an injection from the set
LM of local minima of $\phi $ to the set CC of critical components: $m_k \mapsto E_k$, $k\in \mathbb{N}^2$.
Let $\sigma (k)= \sigma (E_k)$ be the corresponding saddle point value of
$E_k$, with the convention that $\mathbb{R}^n$ is a critical
component and that $\sigma (\mathbb{R}^n)=+\infty $. It is of the form
$E_{1,1}$ where $m_{1,1}$ is a point of global minimum of $\phi$. For $m_k \in \mathrm{LM}$, we
put \ekv{gl.2}{S_k =\sigma (k)-\phi (m_k)>0, \ (S_1=+\infty ).}
For simplicity, let us arrange the indices $k\in \mathbb{N}^2$ in a suitable order, to be chosen below, and write them as $k_1$, $k_2, \ldots k_{n_0}$,
with $k_1 = (1,1)$.

\begin{theo}\label{gl4}
  We label the local minima as above, so that $E_{1,1}=\mathbb{R}^n$, and let
  $\mu_2,...,\mu_{n_0}$ denote the $n_0-1$ non-vanishing
  eigenvalues of $-\Delta _A^{(0)}$, which are $o(h)$. For $h$ small
  enough they are real and exponentially small. More precisely, with a
  suitable labelling of the eigenvalues, we have \ekv{gl.3}{\mu_{\ell} \asymp
    he^{-2S_{k_{\ell}}/h},\quad 2\leq \ell \leq n_0.}
\end{theo}

\par
Next, we recall that in Section \ref{lbl} we have constructed a basis for the subspace $E^{(0)}$, given by $e_k^{(0)} = \Pi^{(0)} (f_k^{(0)})$,
$k\in \mathbb{N}^2$, and checked that the system $(e_k^{(0)})$ is uniformly linearly independent in $L^2$. Also, in Section \ref{genc}, we have
introduced a basis of one-forms $(e_j^{(1)})$, $1\leq j \leq n_1$, for the subspace $E^{(1)}$, associated to the saddle points, and this
construction is applicable here as well. In particular, Proposition 
\ref{estqm2} remains valid.

\par
Similarly to the generic case, we can analyze the matrix of $d_\phi :E^{(0)}\to E^{(1)}$ with respect to the bases $(e_k^{(0)})$ and
$(e_j^{(1)})$, and verify the following properties,
$$
d_\phi e_{k_1}^{(0)}=0,\quad k_1 = (1,1),
$$
while for $\ell \geq 2$, we have
$$
d_\phi e_{k_{\ell}} ^{(0)}=\sum_{j=1}^{n_1}r_{j,\ell}e_j^{(1)},
$$
where for some fixed $\alpha >0$,
\begin{enumerate}[1)]
\item $r_{j,\ell}=h^{\frac{1}{2}}b_{j,\ell}(h)e^{-S_{k_{\ell}}/h}$, when $s_j\in
  \mathrm{SSP}\cap\partial E_{k_{\ell}}$. Here
$b_{j,\ell}(h)\sim \sum_{\nu =0}^\infty b_{j,\ell}^\nu h^\nu$, $b_{j,\ell}^0\ne 0$.

\item $r_{j,\ell}=h^{\frac{1}{2}}\left( {\cal O}(e^{-(\phi (s_j)-\phi
    (m_{k_{\ell}}))/h})+{\cal O}(e^{-(S_{k_{\ell}}+ \alpha
    )/h})\right) $, when $s_j\in
  \mathrm{SSP}\cap \partial \widetilde{E}$ and $\widetilde{E}\ne E_{k_{\ell}}$
  is a critical component containing $E_{k_{\ell}}$.

\item $r_{j,\ell}={\cal O}(h^{\frac{1}{2}}e^{-(S_{k_{\ell}}+\alpha )/h})$
  otherwise, i.e. when $s_j$ is a saddle point that is not on the boundary of a critical component
  containing $E_{k_{\ell}}$.
\end{enumerate}

The matrix of $d_\phi :E^{(0)}\to E^{(1)}$ is of the form
$(0\,R')=R=(r_{j,\ell})$, where $(j,\ell)\in \{1,...,{n_1}\}\times
\{1,...,{n_0}\}$ and $0$ indicates the first vanishing column. Let $D$ be the $(n_0-1)\times (n_0-1)$
matrix $\mathrm{diag\,}(e^{-S_{k_{\ell}}/h})_{\ell=2,\ldots\, ,n_0}$, and write
\ekv{gl.7}{
R'=h^{\frac{1}{2}}\widetilde{R}D,\quad \widetilde{R}=(\widetilde{r}_{j,\ell}),}
so that $\widetilde{r}_{j,\ell}$ is equal to $b_{j,\ell}$,
${\cal O}(e^{-(\phi (s_j)-\phi (m_{k_{\ell}})-S_{k_{\ell}})/h})+{\cal O}(e^{-\alpha
/h})$, ${\cal O}(e^{-\alpha /h})$ respectively, in the three
different cases above. Notice that in the second case, $\phi
(s_j)-\phi (m_{k_{\ell}})-S_{k_{\ell}}>0$ but if we want to allow $\phi$ to vary nicely
with parameters, we have no uniform bound from below by a positive constant.

\subsection{Singular values}\label{si}
Clearly $\Vert \widetilde{R}\Vert\le {\cal O}(1)$ and we shall prove
that $\widetilde{R}$ is injective with a uniformly bounded left
inverse satisfying
\ekv{gl.8}{\Vert \widetilde{R}^{-1}\Vert \le {\cal O}(1). }

Accepting this in this subsection, we shall now estimate the singular
values of $R$. Using the standard notation, we let
$s_1(\widetilde{R})\ge ...\ge s_{n_0-1}(\widetilde{R})$ be the
singular values of $\widetilde{R}$, i.e. the eigenvalues of
$(\widetilde{R}^*\widetilde{R})^{\frac{1}{2}}$. Then in view of
(\ref{gl.8}), we have
\ekv{gl.9}{
\frac{1}{C}\le s_j(\widetilde{R})\le C. }
The Ky Fan inequalities and (\ref{gl.7})  tell us that
$$
s_{\nu} (R')\le h^{\frac{1}{2}}\Vert \widetilde{R}\Vert s_{\nu}(D),
$$
where we recall that $\Vert R\Vert={\cal O}(1)$. Without loss of generality, we
may assume that $S_{k_2} \ge S_{k_3} \ge ...\ge S_{k_{n_0}}$. Then $s_{\nu}
(D)=e^{-S_{k_{n_0+1-{\nu}}}/h}$, and we get the upper bound,
\ekv{gl.10}
{
s_{\nu}(R')\le {\cal O}(1)h^{\frac{1}{2}}e^{-S_{k_{n_0+1-{\nu}}}/h}.
}

To get a lower bound we shall use (\ref{gl.9}) and write
$$
\widetilde{R}=h^{-\frac{1}{2}}R'D^{-1}
$$
so that the Ky Fan inequalities give,
$$
\frac{1}{{\cal O}(1)}\le s_{n_0-1}(\widetilde{R})\le h^{-\frac{1}{2}}s_{\nu}(R')s_{n_0-\nu}(D^{-1}).
$$
Here $s_{n_0-\nu}(D^{-1})=e^{S_{k_{n_0+1-\nu}}/h}$ and we get
$$
s_{\nu} (R')\ge \frac{h^{1/2}}{{\cal O}(1)}e^{-S_{k_{n_0+1-\nu}}/h},\
{1\le \nu \le n_0-1.}
$$
In conclusion,
\ekv{gl.11}
{
s_{\nu} (R')\asymp h^{\frac{1}{2}}e^{-S_{k_{n_0+1-\nu}}/h}.
}

We can now end the {\bf proof of Theorem \ref{gl4}.}. We may assume
that $e_1^{(1)},...,e_{n_1}^{(1)}$ is an orthonormal basis in
$E^{(1)}$ for the scalar product $(\cdot |\cdot )_{A,\kappa }$. Let
$g_1^{(0)},...,g_{n_0}^{(0)}$ be an orthonormal basis in $E^{(0)}$ for
the scalar product $(\cdot |\cdot )_{A,\kappa }$ on functions. Then,
as we have seen, the eigenvalues $\mu_{\ell}$ are the squares of the
singular values of the matrix $\widehat{R}$ of $d_\phi :E^{(0)}\to E^{(1)}$ with
respect to the two orthonormal bases. However the uniform linear
independence of the basis $e_{k_1}^{(0)},...,e_{k_{n_0}}^{(0)}$ means that
$\widehat{R}=RU$ where $U$ and $U^{-1}$ are uniformly bounded. By the
Ky Fan inequalities we then know that the singular values of $R$ and
those of $\widehat{R}$ are pairwise of the same order of
magnitude. The theorem therefore follows from (\ref{gl.11}). \hfill{$\Box$}

\subsection{Proof of (\ref{gl.8})}\label{pfgl8}

We define the indicator matrix $\widetilde{R}_0$ of $\widetilde{R}$
by replacing $\widetilde{r}_{j,\ell}$ by $0$ in the third case in the
description of $\widetilde{R}$ above (after (\ref{gl.7})). Then,
$$
\widetilde{R}_0=(\widehat{r}_{j,\ell})_{{1}\le j\le n_1\atop
2\le \ell \le n_0},
$$
where
\begin{enumerate}[1)]
\item $\widehat{r}_{j,\ell}=\widetilde{r}_{j,\ell}=b_{j,\ell}$ when $s_j\in
  \mathrm{SSP}\cap\partial E_{k_{\ell}}$,
\item  $\widehat{r}_{j,\ell}=\widetilde{r}_{j,\ell}={\cal O}(e^{-(\phi
    (s_j)-\phi (m_{k_{\ell}})-S_{k_{\ell}} )/h}+e^{-\alpha /h})$ when $s_j\in
  \mathrm{SSP}\cap \partial \widetilde{E}$ and $E_{k_{\ell}} \subset
  \widetilde{E}\in \mathrm{CC} $, $E_{k_{\ell}}\ne \widetilde{E}$,
\item $\widehat{r}_{j,\ell}=0$ otherwise.
\end{enumerate}

Then $\widetilde{R}=\widetilde{R}_0+{\cal O}(e^{-\alpha /h})$ and in
order to prove (\ref{gl.8}) it suffices to prove an a priori estimate
\ekv{gl.12}
{
\Vert u\Vert\le {\cal O}(1)\Vert \widetilde{R}_0u\Vert,\quad u\in \mathbb{C}^{n_0-1},
}
since we then get an analogous estimate for $\widetilde{R}$.

We shall prove (\ref{gl.12}) by estimating the components of $u$ one
after the other in a suitable order and thus eliminating them
successively.

\par If SSP is empty or equivalently, if $m_{1,1}$ is the only local
minimum of $\phi $, then there is nothing to prove.

\par If $\mathrm{SSP}\ne \emptyset $, let $\sigma_{2}$ be the {\it
  smallest} value in $\phi (\mathrm{SSP})$. Here we are not using the labelling introduced in Section \ref{lbl}.
The components of $\phi^{-1}(]-\infty ,{\sigma_2}[)$ can be split into the critical components, labelled as $E_k$, with
$k=(k^{\sigma_{2}}, k^{cc})$, $1\leq k^{cc}\leq N_{k^{\sigma_{2}}}$, and the others, say,
$F_1,...,F_M$. For notational convenience, we shall relabel the critical components $E_k$ as
$\widetilde{E}_1,\ldots \widetilde{E}_N$. Here each $\overline{\widetilde{E}}_j$ intersects at least one other, $\overline{\widetilde{E}}_{k}$,
while each $\overline{F}_\ell$ is disjoint from the closure of the other components. Thus the union of the $\widetilde{E}_j$
can be grouped into say, $\widetilde{E}_1\cup...\cup \widetilde{E}_{\ell_1}$,
$\widetilde{E}_{\ell_1+1}\cup ...\cup \widetilde{E}_{\ell_1+\ell_2}$,... such that the corresponding unions of the closures are connected
and mutually disjoint. Here $\ell_1,\ell_2,...$ are all $\ge 2$.

\par In what follows, it will be useful to change the labelling of the local minima whenever we find it convenient and here even
the global minimum $m_{1,1}$ may have to change its notation. Consider one of the groups introduced above, say
$\widetilde{E}_1,...,\widetilde{E}_{\ell_1}$ and write $\ell$ instead
of $\ell_1$ for simplicity.  Let $m_k$ be the unique local minimum of
${{\phi }_\vert}_{\widetilde{E}_k}$ (recalling that $\sigma _2=\inf
\phi (\mathrm{SSP})$). Recall that in the construction of the quasimodes, we have assigned a critical component $E_k$ to
each local minimum and since $\sigma_2$ is minimal, we have
$E_k\supset \widetilde{E}_k$. The same construction also shows that
$E_k\ne \widetilde{E}_k$ for precisely one value of $k$, say for
$k=1$. Notice that the corresponding local minimum $m_1$ is also a
global minimum of the restriction of $\phi $ to $\widetilde{E}_1\cup ...\cup
\widetilde{E}_\ell$.

Since $\overline{\widetilde{E}}_1\cup ...\cup
\overline{\widetilde{E}}_\ell$ is connected, we know that
$\overline{\widetilde{E}}_1\cap \overline{\widetilde{E}}_k\ne\emptyset $
for at least one value of $k\ne 1$, say $k=2$. This intersection is a
finite set of ssps. Choose one of them and denote it by $s_2$. We
shall now estimate $u_2$. Notice that $\widehat{r}_{2,2}=b_{2,2}$ is
elliptic, while $\widehat{r}_{2,1}=0$. Since $s_2$ cannot be in the
boundary of $\widetilde{E}_k$ for any $k\ge 3$, we know (using also
the minimality of $\sigma _2$) that $\widehat{r}_{2,k}=0$ for all other
$k$ associated to any other local minimum different from $m_1$ and
$m_2$. It follows that
$$
(\widetilde{R}_0u)(2)=b_{2,2}u(2),
$$
and by the ellipticity of $b_{2,2}$, we deduce the a priori estimate
$$
|u(2)|\le {\cal O}(1)\Vert \widetilde{R}_0u\Vert .
$$

\par If $\ell \ge 3$, we may assume, after relabelling of
$m_3,...,m_\ell$, that $\overline{\widetilde{E}}_2\cap
\overline{\widetilde{E}}_3\ne \emptyset $. Let $s_3$ be a ssp in
this intersection. Then (using again the minimality of $\sigma
_2$) we have
$$
(\widetilde{R}_0u)(3)=b_{3,2}u(2)+b_{3,3}u(3),
$$
so, using the ellipticity of $b_{3,3}$, we get
$$
|u(3)|\le {\cal O}(1)(\Vert \widetilde{R}_0u\Vert +|u(2)|)\le {\cal
  O}(1)\Vert \widetilde{R}_0u\Vert.
$$

\par Continuing this way, we get
$$
|u(2)|+...+|u(\ell)|\le {\cal O}(1)\Vert \widetilde{R}_0u\Vert.
$$

Thus, in order to prove (\ref{gl.8}), it suffices to do so when
$u(2)=...=u(\ell)=0$. Indeed, for a general $u\in \mathbb{C}^{n_0-1}$,
write
$$
u=u'+\sum_{k=2}^\ell u(k)e_k,\quad u'(k)=0,\hbox{ for }2\le k\le \ell,
$$
where $e_k$ is the canonical basis vector of index $k$, and if we
assume that (\ref{gl.8}) holds when $u(2)=...=u(\ell)=0$, then
$$
\Vert u'\Vert\le {\cal O}(1)\Vert \widetilde{R}_0u'\Vert\le {\cal
  O}(1)(\Vert \widetilde{R}_0u\Vert+{\cal
  O}(1)(|u(2)|+...+|u(\ell)|))\le {\cal O}(1)\Vert
\widetilde{R}_0u\Vert .
$$

\par
In other words, we have eliminated the minima $m_2,...,m_\ell$ from
the discussion and only $m_1$ survives among the
$m_1,m_2,...,m_\ell$. We do the same for the other groups $\widetilde{E}_{\ell_1+1}\cup
...\cup \widetilde{E}_{\ell_1+\ell_2}$, ..., (if there are more than one) so that for
each group, we eliminate all the local minima except one which is also
a global minimum in the corresponding union.

If $\sigma _2$ is the only element of $\phi (\mathrm{SSP})$, then we
only have one group as above and no $F_j$. We have then eliminated all
the $u(k)$ and the proof is complete.

If $\phi (\mathrm{SSP})$ contains more elements, let $\sigma_3 > \sigma _2$ be the smallest one. Again (after a change of notation and
forgetting about the earlier $\widetilde{E}_k$ and $F_\nu$) we have
$$
\phi ^{-1}(]-\infty ,\sigma _3[)=\widetilde{E}_1\cup
...\cup\widetilde{E}_N\bigcup F_1\cup...\cup F_M,
$$
where $\widetilde{E}_\nu $ are the ccs and $F_\mu $ are the other
components. Each of the components may have several local minima but
we have already eliminated all but one, say $m$ which is also a global minimum
for that component, and $\widetilde{E}_\nu \subset E_m$ and $F_\mu
\subset E_m $ respectively. This implies that we can repeat the
procedure of elimination precisely as we did at the level $\sigma
_2$. After finitely many steps all the local minima are eliminated and
we get (\ref{gl.8}).

\subsection{Full asymptotics for the smallest non-vanishing eigenvalue}\label{fa}

We return to the general situation in Subsection \ref{sl} and label
the minima $m_{k_1}$, $m_{k_2}\ldots\, ,m_{k_{n_0}}$ so that $m_{k_1}=m_{1,1}$
and $S_{k_2}$ is maximal among all the $S_{k_j}$, $j=2,3,...,n_0$.  We
assume that there is a gap between $S_{k_2}$ and the other $S_{k_j}$:
\begin{equation}\label{fa.1}
S_{k_2}> \max_{j\ge 3}S_{k_j}=:S'.
\end{equation}
Also assume that
\begin{equation}\label{fa.2}
\partial E_{k_2} \hbox{ contains precisely one ssp.}
\end{equation}

\begin{theo}
\label{fa1}
{We assume {\rm (\ref{fa.1})}, {\rm (\ref{fa.2})}, in
addition to the general assumption of Theorem \ref{gl4}.} Then   
the smallest
non-vani\-shing ei\-gen\-value $\mu _2$ is given by
$$
\mu_2=(h|b_2(h)|^2+{\cal O}(h^\infty ))e^{-2S_2/h},
$$
where we write $S_2$ instead of $S_{k_2}$ for short.
Here $b_2$ is as in Theorem {\rm \ref{theogene}}.
\end{theo}

As before, we may assume that the basis $(e_j^{(1)})$ in $E^{(1)}$ is
orthonormal, while we have uniform linear independence for the basis
$(e_k^{(0)})$ in $E^{(0)}$. We may assume however that
\begin{equation}\label{fa.3}
\Vert e_{k_2}^{(0)}\Vert_\kappa =1.
\end{equation}
Define $R,R',\widetilde{R},D$ as in and around (\ref{gl.7}). Then 
\begin{equation}\label{fa.4}
\widetilde{R}D=\left(\begin{array}{ccc}b_2(h) &{\cal O}(e^{-\alpha
      /h})\\ {\cal O}(e^{-\alpha /h}) &W
 \end{array}\right) \left(\begin{array}{ccc}e^{-S_2/h}&0\\0
        &\widehat{D} \end{array}\right) ,
\end{equation}
where $\widehat{D}=\mathrm{diag\,}(e^{-S_j/h})_{3\le j\le n_0}$ and we
write $S_j$ instead of $S_{k_j}$ for simplicity.

Here, we know that $b_2$ is elliptic of order 0 and that
$\left(\begin{array}{ccc}b_2&0\\0 &W \end{array}\right)$ has a
uniformly bounded left inverse. Consequently,
\begin{equation}\label{fa.4.5}
\frac{1}{{\cal O}(1)}\le W^*W\le {\cal O}(1).
\end{equation}
Theorem \ref{fa1} follows from:

\begin{prop}
\label{fa2}
The smallest singular value of $\widetilde{R}D:\mathbb{R}^{n_0-1}\to
\mathbb{R}^{n_1}$ is equal to $|b_2|e^{-S_2/h}+{\cal
  O}(e^{-(S_2+\alpha )/h})$ for some $\alpha >0$, independent of $h$,
when $\mathbb{R}^{n_1}$ is equipped with the standard Hilbert norm and
$\mathbb{R}^{n_0-1}$ is equipped with the norm $[\cdot ]$, defined by
$$
[x]^2=\sum_{2\le j,k\le n_0}x_jx_k(e_j^{(0)}|e_k^{(k)})_\kappa .
$$
\end{prop}

We know that $[\cdot ]$ is uniformly equivalent to the standard norm
and that $[e_2]=1$, where $e_2=(1,0,...,0)\in \mathbb{R}^{n_0-1}$. The
square of the singular value is equal to
\begin{equation}\label{fa.5}
\inf _{[x]=1}\Vert \widetilde{R}Dx\Vert ^2.
\end{equation}

Eliminating the terms ${\cal O}(e^{-\alpha /h})$ in (\ref{fa.4}) will
modify the quantity (\ref{fa.5}) by a factor $1+{\cal O}(e^{-\alpha
  /h})$ and if we do that elimination, then (\ref{fa.5}) becomes
\begin{equation}\label{fa.6}
\inf_{[x]=1} |b_2|^2 e^{-2S_2/h}x_2^2+\Vert W\widehat{D}x'\Vert ^2,\ x=(x_2,x').
\end{equation}
Taking $x_2=1$, $x'=0$, we see that this infimum is $\le
|b_2|^2e^{-2S_2/h}$, and we shall see that it is very close to this
value. Let $x$ be a point (with $[x]=1$) where the infimum is
attained. Then $\Vert W\widehat{D}x'\Vert = {\cal O}(e^{-S_2/h})$ and
from (\ref{fa.4.5}) and the fact that $\widehat{D}^{-1}={\cal
  O}(e^{S'/h})$, we conclude that $\Vert x'\Vert ={\cal
  O}(e^{-(S_2-S')/h})$. It is then clear that $x_2=\pm 1+{\cal
  O}(e^{-(S_2-S')/h})$ and the infimum in (\ref{fa.6}) is
$$
\ge |b_2|^2e^{-2S_2/h}(1-{\cal O}(e^{-(S_2-S')/h}))\ge |b_2|^2e^{-2S_2/h}-{\cal O}(e^{-(S'+S_2)/h}).
$$
This gives Proposition \ref{fa2} and Theorem \ref{fa1}.

\subsection{An example with three minima and three saddle points}
We shall finally briefly discuss an example illustrating
the generic and general cases, for which computations can be made by
hand. We study the case of three saddle points and three local minima
as illustrated by Figure \ref{3min3sp}. It is natural to denote simply by $m_k$,
$k=1, 2, 3$ the local minima and $s_j$, $j=1, 2, 3$ the saddle points.
The \it generic case \rm corresponds to the case when Hypothesis \ref{genc0}
is satisfied. An interesting situation, entering in the previous
\it general case\rm, {appears when $\phi $ depends smoothly on some
parameter such that we may have
\begin{equation} \label{limitcase}
\phi(m_1) = \phi(m_2) = \phi(m_3), \textrm{ and } \phi(s_1) = \phi(s_2) = \phi(s_3).
\end{equation}
We will restrict the attention to a small neighborhood of a parameter
value for which (\ref{limitcase}) holds.}

\begin{figure}[htbp]
\begin{center}

\begin{picture}(0,0)%
\includegraphics{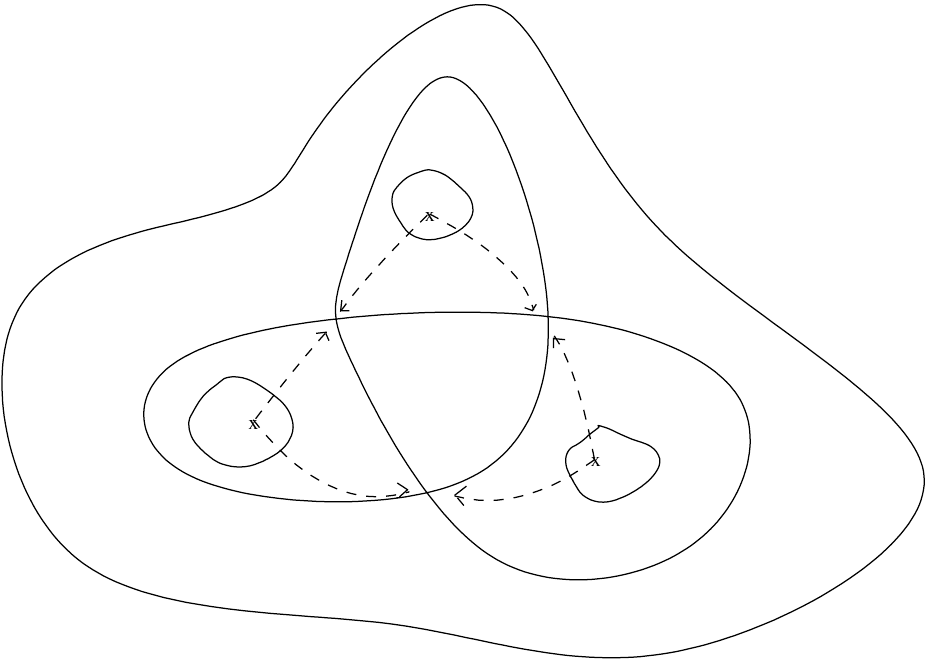}%
\end{picture}%
\setlength{\unitlength}{3158sp}%
\begingroup\makeatletter\ifx\SetFigFont\undefined%
\gdef\SetFigFont#1#2#3#4#5{%
  \reset@font\fontsize{#1}{#2pt}%
  \fontfamily{#3}\fontseries{#4}\fontshape{#5}%
  \selectfont}%
\fi\endgroup%
\begin{picture}(5558,3943)(330,-4270)
\put(1553,-2901){\makebox(0,0)[lb]{\smash{{\SetFigFont{6}{7.2}{\rmdefault}{\mddefault}{\updefault}{\color[rgb]{0,0,0}\small $m_2$}%
}}}}
\put(3995,-3155){\makebox(0,0)[lb]{\smash{{\SetFigFont{6}{7.2}{\rmdefault}{\mddefault}{\updefault}{\color[rgb]{0,0,0}\small $m_3$}%
}}}}
\put(2852,-1525){\makebox(0,0)[lb]{\smash{{\SetFigFont{6}{7.2}{\rmdefault}{\mddefault}{\updefault}{\color[rgb]{0,0,0}\small $m_1$}%
}}}}
\put(3155,-1997){\makebox(0,0)[lb]{\smash{{\SetFigFont{6}{7.2}{\rmdefault}{\mddefault}{\updefault}{\color[rgb]{0,0,0}\small $-b_1$}%
}}}}
\put(3868,-2611){\makebox(0,0)[lb]{\smash{{\SetFigFont{6}{7.2}{\rmdefault}{\mddefault}{\updefault}{\color[rgb]{0,0,0}\small $a_3$}%
}}}}
\put(3303,-3359){\makebox(0,0)[lb]{\smash{{\SetFigFont{6}{7.2}{\rmdefault}{\mddefault}{\updefault}{\color[rgb]{0,0,0}\small $-b_3$}%
}}}}
\put(2301,-3126){\makebox(0,0)[lb]{\smash{{\SetFigFont{6}{7.2}{\rmdefault}{\mddefault}{\updefault}{\color[rgb]{0,0,0}\small $a_2$}%
}}}}
\put(1906,-2428){\makebox(0,0)[lb]{\smash{{\SetFigFont{6}{7.2}{\rmdefault}{\mddefault}{\updefault}{\color[rgb]{0,0,0}\small $-b_2$}%
}}}}
\put(2619,-1941){\makebox(0,0)[lb]{\smash{{\SetFigFont{6}{7.2}{\rmdefault}{\mddefault}{\updefault}{\color[rgb]{0,0,0}\small $a_1$}%
}}}}
\put(3719,-2160){\makebox(0,0)[lb]{\smash{{\SetFigFont{6}{7.2}{\rmdefault}{\mddefault}{\updefault}{\color[rgb]{0,0,0}\small $s_2$}%
}}}}
\put(2167,-2117){\makebox(0,0)[lb]{\smash{{\SetFigFont{6}{7.2}{\rmdefault}{\mddefault}{\updefault}{\color[rgb]{0,0,0}\small $s_3$}%
}}}}
\put(2725,-3409){\makebox(0,0)[lb]{\smash{{\SetFigFont{6}{7.2}{\rmdefault}{\mddefault}{\updefault}{\color[rgb]{0,0,0}\small $s_1$}%
}}}}
\end{picture}%

\caption{Case of 3 saddle points at the same level }
\label{3min3sp}
\end{center}
\end{figure}
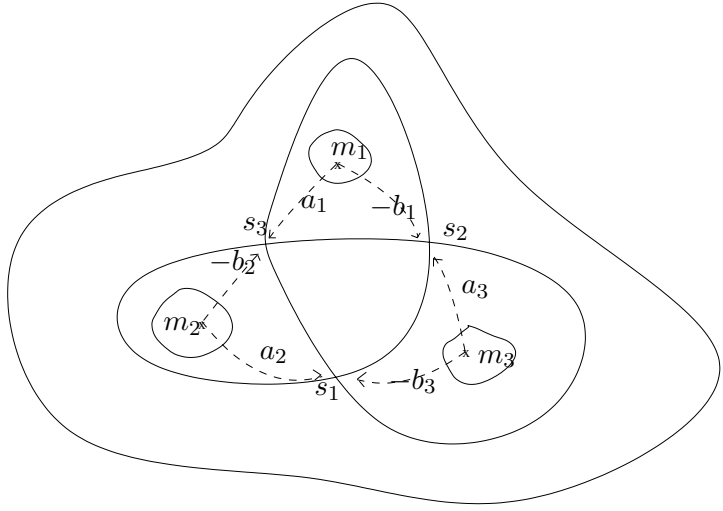

In this situation we can describe in some detail some facts about the restriction of the twisted Laplacian
$-\Delta_A^{(0)}$ to the space $E^{(0)}$ corresponding to the exponentially small eigenvalues.  We know from the previous sections
that this space is of dimension 3, since this is the number of local
minima. The question is to study the eigenvalues in more detail, and we
already know that one of them is $0$ and that they are all real. 

\par {It will be convenient to use a different system of
  critical components and define $E(m_k)$ to be the connected
  component of $\phi ^{-1}(]-\infty ,\min_{j\ne k}\phi (s_j)[)$ that
  contains $m_k$ and to define $\chi _k$ correspondingly as in Section
  \ref{lbl}. The eigenspace $E^{(0)}$ is generated by
0-forms of the form
$$
e_k^{(0)} =  h^{-n/4} c_k(h)\chi_k(x) e^{-(\phi-\phi(m_k))/h} + \ooo(e^{-1/Ch}),
$$
where $c_k(h)\sim c_{k,0}+hc_{k,1}+...$ is a normalization constant
with $c_{k,0}>0$.}
Similarly we know that  1-forms generating  the 3-dimensional eigenspace $E^{(1)}$ of $-\Delta^{(1)}_A$ associated with
exponentially small eigenvalues can be chosen of the following form:
$$
e_j^{(1)}=(h^{-n/4}a_j(x;h) +r_j(x))\theta_j(x)e^{-\phi _+(x)/h} + \ooo( e^{-1/Ch}),
$$
where $a_j(x;h)\sim a_{j,0}(x)+ha_{j,1}+... \hbox{ in }C^\infty
(\mathrm{neigh\,}(s_j,\mathbb{R}^n))$ and $r_j$ are as in
(\ref{ps.2bis}). Both families $(e^{(0)}_k)$ and $(e^{(1)}_j)$ are
bases {that we can assume to be orthonormal for the $A,\kappa $ scalar
  products.}

\par {The matrix $R$ of $d_\phi :E^{(0)}\to E^{(1)}$ with
  respect to the bases $(e_k^{(0)})$ and $(e_j^{(1)})$ takes the form
$$
R=h^{1/2}(\rho _{j,k}e^{-\phi (s_j)/h}e^{\phi (m_k)/h}),
$$
where $|\rho _{j,k}|\asymp 1$ for $j\ne k$, and $={\cal O}(1)e^{-\frac{1}{Ch}}$
when $j=k$, uniformly for the parameter in a small neighborhood of the
value where (\ref{limitcase}) holds. The Maxwellian $e^{-\phi /h}$ can
be written
$$
e^{-\phi /h}=h^{\frac{n}{4}}\sum_k \frac{(1+{\cal
    O}(e^{-\frac{1}{Ch}}))}{c_k(h)}e^{-\frac{\phi (m_k)}{h}}e_k^{(0)},
$$
and is annihated by $d_\phi $, so we know that the coefficent vector
belongs to the kernel of $R$. It follows that 
$$
\sum_{k;\, k\ne j}\rho _{j,k}\frac{1}{c_k}={\cal O}(e^{-\frac{1}{Ch}}),
$$
so 
$$
\rho _{j,k}=\rho _{j,k}^0 +{\cal O}(e^{-\frac{1}{Ch}}),
$$
where 
$$
\rho _{j,k}^0=\begin{cases}0, \ k=j, \\d_jc_k,\ k=j+1,\\ -d_jc_k,\
  k=j+2\end{cases},\quad d_j\asymp 1,\ c_k\asymp 1
$$
with the cyclic convention for the indices. Introducing 
$$
\mu _k=c_ke^{\phi (m_k)/h},\ \sigma _j=d_je^{-\phi (s_j)/h},
$$
we get 
$$
R=h^{1/2}R_0+e^{-\frac{1}{Ch}}({\cal O}(1)e^{-(\phi (s_j)-\phi (m_k))/h}))_{j,k},
$$
where 
$$
R_0 =
\left(
\begin{matrix}
0 & a_2 & -b_3 \\
-b_1 & 0 & a_3 \\
a_1 & -b_2 & 0
\end{matrix}\right),
$$
$$
a_1 = \sigma_3\mu_1, \ \ a_2 = \sigma_1\mu_2, \ \ \  a_3 = \sigma_2\mu_3,
$$
$$
b_1 = \sigma_2\mu_1, \ \ b_2 = \sigma_3\mu_2, \ \ \  b_3 = \sigma_1\mu_3.
$$
}

We see that
{$R_0 \left(\begin{matrix}
\mu_1^{-1} \\ \mu_2^{-1} \\ \mu_3^{-1} \end{matrix}\right) =0$}
reflecting the fact that $d_\phi (e^{-\phi/h}) = 0$. This implies that
$\det R_0 = 0$. We have
$$
R_0^* R_0 =
\left(\begin{matrix}
0 & -b_1 & a_1 \\
a_2 & 0 & -b_2 \\
-b_3 & a_3 & 0
\end{matrix}\right)
\left(\begin{matrix}
0 & a_2 &-b_3 \\
-b_1 & 0 & a_3 \\
a_1 & -b_2 & 0
\end{matrix}\right) =
\left(\begin{matrix}
a_1^2 + b_1^2 & -a_1b_2 &-b_1a_3 \\
-b_2a_1 & a_2^2+b_2^2 & -a_2b_3 \\
-a_3b_1 & -b_3a_2 & a_3^2+b_3^2
\end{matrix}\right).
$$
Since $\det (R_0^* R_0 ) = 0$ we get that for all $\lambda \in \R$,
\begin{equation} \label{dett}
\begin{split}
& \det( \lambda - R_0^* R_0) = \lambda^3 -
(a_1^2 + a_2^2 + a_3^2 + b_1^2 + b_2^2 + b_3^2)\lambda^2 \\ & \ \ \ \ \ \ +
\underbrace{\big[ (a_2^2 a_3^2  + b_2^2 b_3^2 + b_2^2 a_3^2) + (a_3^2 a_1^2  + b_3^2 b_1^2 + b_3^2 a_1^2) +  (a_1^2 a_2^2  + b_1^2 b_2^2 + b_1^2 a_2^2)\big]}_{\defegal D} \lambda.
\end{split}
\end{equation}
To shorten the notation we let $\alpha_j= a_j^2$, $\beta_j= b_j^2$ and $\gamma_j = \alpha_j + \beta_j$ for $j= 1, 2, 3$.  Then one of
the eigenvalues of $R_0^* R_0$ is $0$ and the other two are positive and given by
$$
\lambda_\pm = \frac{\gamma_1 + \gamma_2 + \gamma_3}{2} \pm \sep{ \sep{ \frac{\gamma_1 + \gamma_2 + \gamma_3}{2}}^2 - D}^{1/2}.
$$
Notice that we already know that the eigenvalues are real nonnegative thanks to Section \ref{sy}. In particular we have
\begin{equation} \label{ineqjoh}
0 < D \leq \sep{\frac{\gamma_1 + \gamma_2 + \gamma_3}{2}}^2.
\end{equation}
We can localize the eigenvalues a little more: we get for example that
$$
0 < \lambda_- \leq \frac{\gamma_1+\gamma_2 + \gamma_3}{2} \leq \lambda_+ < \gamma_1+\gamma_2 + \gamma_3
$$
 and also from the fact that that $\lambda_+ \lambda_- = D$, we get
 $$
 \frac{D}{\gamma_1 + \gamma_2 + \gamma_3} < \lambda_- \leq \frac{2D}{\gamma_1 + \gamma_2 + \gamma_3}.
 $$
Now we can study the case when (\ref{limitcase}) occurs, as illustrated by Figure \ref{3min3sp}.
This case is included in the general situation and no longer in the generic one, since in particular
3 critical connected components appear when we are at the critical level (see Figure \ref{3min3sp}).
If we assume in addition that $b_1= b_2 = b_3$ and $a_1 = a_2 = a_3$, which corresponds roughly to the rotation invariant case,
then we get $\alpha_j = \beta_j =: \alpha$, $\gamma_j = 2\alpha$ for $j=1, 2,3$, and $D= 9 \alpha^2 = \sep{ \frac{\gamma_1+\gamma_2 + \gamma_3}{2}}^2$,
so that the eigenvalues $\lambda_+$ and $\lambda_-$ of $R_0$ are equal.
 $$
 \lambda_+ = \lambda_- = 3\alpha \approx h e^{-2(\phi(s) -\phi(m))/h},
 $$
where $s$ (respectively $m$) is any of the saddle points (respectively minima).

\remark The right-hand side inequality in  (\ref{ineqjoh}) also reads
\begin{equation*}
\begin{split}
0 \leq & \  (\gamma_1 + \gamma_2 + \gamma_3)^2 - 4D \\
=  & \ (\alpha_1 + \alpha_2 + \alpha_3 + \beta_1+\beta_2 + \beta_3)^2 \\
 & \ \ \ \  - 4 ( \alpha_2 \alpha_3  + \beta_2  \beta_3 +  \beta_2 \alpha_3) -4 (\alpha_3 \alpha_1  +  \beta_3  \beta_1 +  \beta_3 \alpha_1) -4  (\alpha_1 \alpha_2  +  \beta_1  \beta_2 +  \beta_1 \alpha_2) \big) \\
=  & \  ( \alpha_1 - \beta_3)^2  + (\alpha_2 - \beta_1)^2 + (\alpha_3 - \beta_2)^2  \\
& \ \ \ \  - 2( \alpha_1 - \beta_3)(\alpha_2 - \beta_1)  - 2( \alpha_1 - \beta_3)(\alpha_3 - \beta_2)
-2 (\alpha_2 - \beta_1)(\alpha_3 - \beta_2) 
 \end{split}
 \end{equation*}
Let us denote by $\ddd$ this last expression. From the general study we know that $0$ is an eigenvalue and this implies 
\begin{equation} \label{equality}
\alpha_1 \alpha_2 \alpha_3  = \beta_1  \beta_2  \beta_3
\end{equation}  
since $\alpha_j =a_j^2$, $\beta_j = b_j^2$, and 
\begin{multline}
\det( \lambda - R_0^* R_0) = \lambda^3 -
(a_1^2 + a_2^2 + a_3^2 + b_1^2 + b_2^2 + b_3^2)\lambda^2 \\
+ \big[ (a_2^2 a_3^2  + b_2^2 b_3^2 + b_2^2 a_3^2) + (a_3^2 a_1^2  + b_3^2 b_1^2 + b_3^2 a_1^2) +  (a_1^2 a_2^2  + b_1^2 b_2^2 + b_1^2 a_2^2)\big]  \lambda \\
- (a_1a_2a_3 -b_1b_2b_3)^2.
\end{multline}
Just notice that without the assumption (\ref{equality}), $\ddd$ may be negative: this can happen for example if we suppose
$  \alpha_1 - \beta_3 = \alpha_2 - \beta_1 = \alpha_3 - \beta_2 \defegal \delta >0 $ 
since we get then 
$$
\ddd = - 3 \delta^2 <0.
$$
Note that in this case $\alpha_1\alpha_2\alpha_3 > \beta_1\beta_2\beta_3$ by direct computation so that of course (\ref{equality}) is not satisfied.

\end{document}